\documentclass[11pt]{article}

\usepackage{fullpage}

\usepackage{amsmath}

\usepackage[font={footnotesize}]{caption}

\usepackage[pdftex]{graphicx}

\usepackage[numbers,sort&compress,super]{natbib}
\newcommand{\BibTeX}{\textsc{Bib}\TeX}       % corect BibTeX appearance
\newcommand\T{\rule{0pt}{3ex}}
\newcommand\B{\rule[-1.8ex]{0pt}{0pt}}

\title{Maximum Fidelity}

\author{Ali Kinkhabwala\\
		Department of Systemic Cell Biology\\
		Max Planck Institute of Molecular Physiology\\
	   	Otto-Hahn-Str.~11\\
	   	44227 Dortmund, Germany\\
		email: ali.kinkhabwala@mpi-dortmund.mpg.de
}

\begin{document}

% generate the title page from the info in the headers above
\maketitle

% 200 words max Abstract
\abstract{
The most fundamental problem in statistics is the inference of an unknown probability distribution from a finite number of samples.  For a specific observed data set, answers to the following questions would be desirable:  (1) \textit{Estimation}:  Which candidate distribution provides the \textit{best} fit to the observed data?, (2) \textit{Goodness-of-fit}: How concordant is this distribution with the observed data?, and (3) \textit{Uncertainty}:  How concordant are other candidate distributions with the observed data?  A simple unified approach for univariate data that addresses these traditionally distinct statistical notions is presented called ``maximum fidelity''.  Maximum fidelity is a frequentist approach (in the strict manner of Peirce and Wilson) that is fundamentally based on model \textit{concordance} with the observed data. The fidelity statistic is a general information measure based on the coordinate-independent cumulative distribution and critical yet previously neglected symmetry considerations. A highly accurate gamma-function approximation for the distribution of the fidelity under the null hypothesis (valid for any number of data points) allows direct conversion of fidelity to absolute model concordance ($p$ value), permitting immediate computation of the concordance landscape over the entire model parameter space. Maximization of the fidelity allows identification of the most concordant model distribution, generating a method for parameter estimation. Neighboring, less concordant distributions provide the ``uncertainty'' in this estimate. Detailed comparisons with other well-established statistical approaches reveal that maximum fidelity provides an optimal approach for parameter estimation (superior to maximum likelihood) and a generally optimal approach for goodness-of-fit assessment of arbitrary models applied to any number of univariate data points distributed on the circle or the line. Extensions of this approach to binary data, binned data, and multidimensional data are described, along with improved versions of classical parametric and nonparametric statistical tests. Maximum fidelity provides a philosophically consistent, robust, and seemingly optimal foundation for statistical inference.  All findings in this manuscript are presented in an elementary way in order to be immediately accessible to researchers in the sciences, social sciences, and all other fields utilizing statistical analysis (medicine, engineering, finance, etc.).
}

\clearpage

%\listoffigures

%%%%%%%%%%%%%%%%%%%%%%%%%%%%%%%%%%%%%%%%%%%%%%%%%%%%%%%%%%%%%%
\section{Introduction}\label{sec:intro}
%%%%%%%%%%%%%%%%%%%%%%%%%%%%%%%%%%%%%%%%%%%%%%%%%%%%%%%%%%%%%%

Researchers seeking simple methods to properly analyze their data are instead confronted by a bewildering miscellany of statistical tools and approaches.
%\cite{kanji_100_2006}
This is primarily due to the existence of multiple distinct and irreconcilable philosophies of statistical inference. The research community would greatly benefit by the establishment of a single, uncontroversial, philosophically consistent, and sufficiently general approach for the optimal estimation of model parameters (and their ``uncertainty'') and the assessment of goodness-of-fit for arbitrary models applied to data. In this manuscript, I present an approach called ``maximum fidelity''  that appears to uniquely achieve these aims for the analysis of univariate data. For multivariate data, maximum fidelity allows for a proper understanding of why it is impossible to generate a single unique approach; however, after a generic choice of convention, maximum fidelity can also be extended to the analysis of multidimensional data.

The fidelity is an information-based statistic derived from the cumulative distribution (the probability integral 
transform\cite{pearson_note_1902,pearson_method_1933}), the latter only uniquely defined for univariate data. Maximization of the fidelity leads to an optimal means of parameter estimation for arbitrary models applied to univariate data. A simple conversion of the fidelity to its associated $p$ value (under the null hypothesis) allows for efficient, intuitive, and generally optimal estimation of model concordance. Similar concordance evaluation for neighboring models in the parameter space permits a notion of parameter ``uncertainty''. While maximum fidelity is only uniquely defined for univariate data (due to its reliance on the cumulative distribution), its application to higher dimensional data is also possible upon transformation of multidimensional data to univariate data through invocation of appropriate model and/or coordinate system symmetries (referred to below as ``inverse Monte Carlo"). 

\begin{figure}[!b]
   \begin{center}
\vspace{-8.5cm}
     \includegraphics*[width=6.5in]{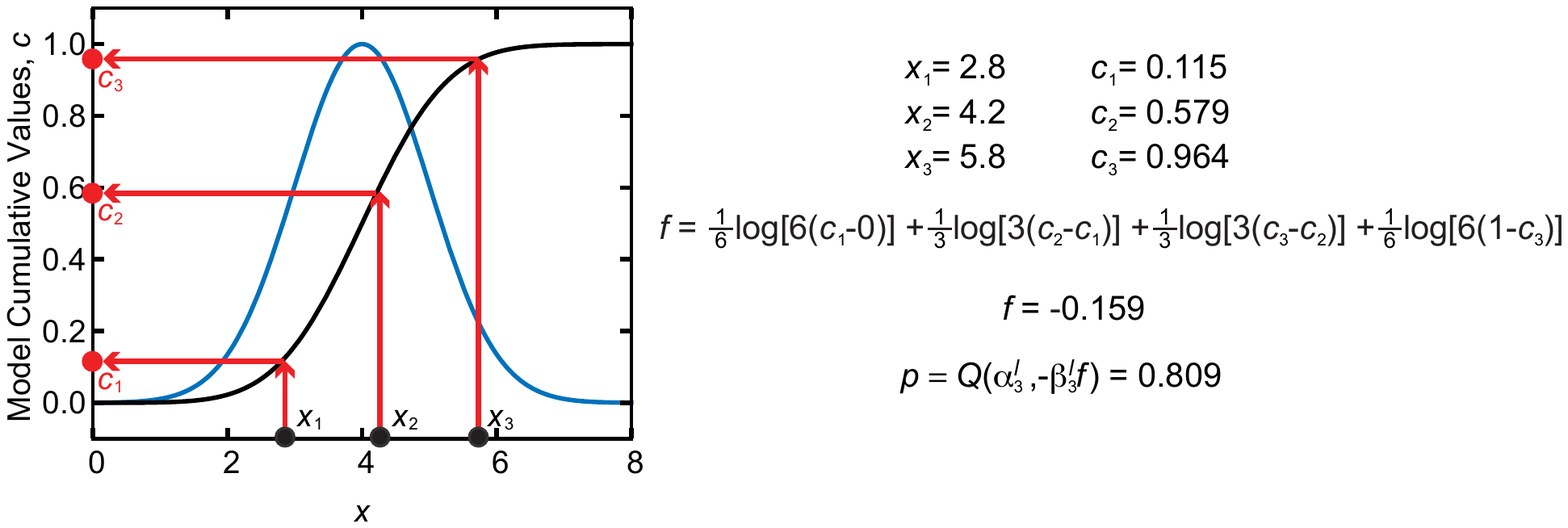}
\vspace{-10.5cm}
      \caption{Example of the use of the fidelity for testing the concordance of a hypothesized Gaussian model ($x_0=4$, $\sigma=1$) with three observed data points. Here, the peak-normalized Gaussian probability distribution (blue) is also expressed in terms of its cumulative distribution (black), which maps the points $x_i$ to their corresponding cumulative values $c_i$. The fidelity based on these particular $c_i$ values is $f=-0.159$, which can be immediately converted to its corresponding $p$ value through a highly-accurate gamma-function approximation, yielding $p=0.809$ (with the standard notion of $p>0.05$ indicating a good fit).}
      \label{fig:fidover}
   \end{center}
\end{figure}

An illustrative example of the use of the fidelity for determining the level of concordance of a Gaussian model with a particular data set is shown in Fig.~\ref{fig:fidover}. Depicted are three observed data points: $x_1$, $x_2$, and $x_3$. A Gaussian model is hypothesized to fit these data with its probability distribution shown in blue and its corresponding cumulative distribution in black. The three data points map to model-dependent cumulative values, from which the fidelity can be determined, which is the sum of the logarithm of the cumulative spacings between the points with a critical additional factor-of-two at the boundaries. Using a simple yet highly accurate gamma-function approximation of the fidelity distribuiton under the null hypothesis, the fidelity can be immediately converted to concordance ($p$ value), allowing for an absolute assessment of model goodness-of-fit to the data. While a Gaussian model was chosen for this example, maximum fidelity is equally applicable to arbitrary probability distributions.

Maximum fidelity represents a significant departure from traditional approaches to statistical inference. Many traditional approaches are ultimately premised on a logical fallacy concerning the notion of \textit{probability}, as Charles Sanders Peirce (1839-1914) already perceptively noted over a century ago:
\begin{quote}
``The relative probability of this or that arrangement of Nature is something which we should have a right to talk about if universes were as plenty as blackberries, if we could put a quantity of them in a bag, shake them well up, draw out a sample, and examine them to see what proportion of them had one arrangement and what proportion another. But, even in that case, a higher universe would contain us, in regard to whose arrangements the conception of probability could have no applicability.''  (Peirce 2.684, i.e., volume 2, section 684 of his collected papers)\cite{peirce_collected_1974}
\end{quote}
In traditional approaches, this \textit{probability fallacy} is made from the outset of the statistical inference (e.g. inverse probability or Bayesian methods\cite{bayes_essay_1763,stigler_laplaces_1986,laplace_memoir_1986}, upon the conclusion of the statistical inference (e.g. various \textit{probabilistic} interpretations of the probable error or of the confidence interval), or as justification of a particular statistic (frequentist coverage probability). A related fallacy (or unnecessary limitation) made in many traditional approaches is the restriction of consideration to a particular parametrized family of distributions (\textit{parameter fallacy}); Bayesian and most frequentist approaches absolutely require this \textit{ad hoc} restriction. Many statistical arguments are also only valid for a particular choice of coordinate system (\textit{coordinate fallacy}), as they rely on coordinate-dependent moments like the mean and standard deviation of the data (e.g. frequentist confidence interval for a Gaussian based on the data mean) or the mean and standard deviation of statistics associated with the model parameters of a specific distribution family (e.g., the Fisher information\cite{savage_rereading_1976} and the related Cram\'er-Rao efficiency\cite{cramer_mathematical_1946,rao_information_1945}). For example, a statistic that achieves the Cram\'er-Rao efficiency for estimation of $\sigma$ for a Gaussian distribution will not be efficient with respect to estimation of the variance $v=\sigma^2$. Equally troubling, for every statistic that asymptotically satisfies the Cram\'er-Rao efficiency, another statistic can always be found that is even more efficient (so-called superefficiency\cite{stoica_evil_1996}). And, again, what if the underlying data are \textit{not} actually drawn from a member of the distribution family assumed for these efficiency calculations (\textit{parameter fallacy})? Many statistical approaches are based on proofs that are only valid in the asymptotic limit of large numbers of data points (\textit{asymptotic fallacy}, e.g. the Cram\'er-Rao efficiency), which yields no insight into their reliability on the smaller data sets for which statistical analysis is most relevant. Another significant fallacy is the use of the likelihood function and likelihood ratio (\textit{likelihood fallacy}). The likelihood itself is a coordinate-dependent quantity. While the likelihood ratio is coordinate independent, its use mandates a statistics based merely on comparison of the \textit{relative} ``quality of fit'' of one model distribution to another from the same family, which provides no information about the absolute concordance of each individual model with the observed data (e.g. in the sense of Pearson's $\chi^2$-$p$ test\cite{pearson_criterion_1900}) and prevents meaningful cross-comparison of models from different distribution families. Parameter estimation via maximization of the likelihood can also fail for certain, otherwise well-defined distribution families\cite{pitman_basic_1979,cheng_estimating_1983,ranneby_maximum_1984,ekstrom_alternatives_2008}.

The specific aspects to which a proper approach to statistical inference should conform were outlined by Peirce in his extensive writings\cite{peirce_collected_1974}, particularly in his unorthodox yet well-motivated definition of \textit{induction}:
\begin{quote}
``Induction is the experimental testing of a theory. The justification of it is that, although the conclusion at any stage of the investigation may be more or less erroneous, yet the further application of the same method must correct the error. The only thing that induction accomplishes is to determine the value of a quantity. It sets out with a theory and it measures the degree of concordance of that theory with fact.'' (Peirce 5.145)\cite{peirce_collected_1974}
\end{quote}
Peirce's central notion of determining the ``degree of concordance'' of a theory will serve as a guiding light throughout this manuscript.
An inductive inference must avoid the above-mentioned \textit{probability fallacy} in order to be valid:
\begin{quote}
``Every argument or inference professes to conform to a general method or type of reasoning, which method, it is held, has one kind of virtue or another in producing truth. In order to be valid the argument or inference must really pursue the method it professes to pursue, and furthermore, that method must have the kind of truth-producing virtue which it is supposed to have. For example, an induction may conform to the formula of induction; but it may be conceived, and often is conceived, that induction lends a probability to its conclusion. Now that is not the way in which induction leads to the truth. It lends no definite probability to its conclusion. It is nonsense to talk of the probability of a law, as if we could pick universes out of a grab-bag and find in what proportion of them the law held good. Therefore, such an induction is not valid; for it does not do what it professes to do, namely, to make its conclusion probable. But yet if it had only professed to do what induction does (namely, to commence a proceeding which must in the long run approximate to the truth), which is infinitely more to the purpose than what it professes, it would have been valid.''  (Peirce 2.780)\cite{peirce_collected_1974}
\end{quote}
Peirce's arguments have been viewed as anticipating the work of the dominant frequentist schools of the twentieth century represented by Ronald Fisher, Jerzy Neyman, and Egon Pearson\cite{levi_induction_1980,hacking_theory_1980,mayo_error_1996}. While Peirce may have been sympathetic to some of what is contained in these approaches, his views are actually more closely aligned with the pragmatic viewpoint of Egon Pearson's father and Peirce's partial contemporary, Karl Pearson, and especially with the \textit{strict} and uncontroversial frequentist viewpoint espoused by Edwin Bidwell Wilson, a lifelong admirer of Peirce's work\cite{hacking_theory_1980}.

The main distinctions between Karl Pearson's viewpoint and that of the later frequentist schools arise, on the one hand, from Pearson's general emphasis on goodness-of-fit (his invention\cite{pearson_criterion_1900}) and, on the other hand, from the latter's embrace of the \textit{parameter} and \textit{likelihood fallacies} described above. Pearson clearly took a more pragmatic view on the choice of model distribution for fitting data (in the following, a ``normal curve'' is a Gaussian): ``I have never found a normal curve to fit anything if there are enough observations!''\cite{pearson_statistical_1935,inman_karl_1994}  Pearson elaborated his views (in the following, a ``graduation curve'' is simply an hypothesized model distribution):
\begin{quote}
``The reader will ask: `But if they do not represent laws of Nature, what is the value of graduation curves?' He might as well ask what is the value of scientific investigation! A good graduation curve --- that is, one with an acceptable probability --- is the only form of `natural law', which the scientific worker, be he astronomer, physicist or statistician can construct. Nothing prevents its being replaced by a better graduation; and ever better graduation is the history of science.''\cite{pearson_statistical_1935,inman_karl_1994}
\end{quote}
Pearson loosely uses the phrase ``acceptable probability'', for which he implied an acceptable $p$ value derived from his $\chi^2$ test\cite{pearson_criterion_1900}, which represented the first completely general measure of model \textit{concordance} comparable across arbitrary distributions. Pearson clearly would not have at all been surprised if a ``better graduation'' came in the form of an arbitrarily different model (e.g. from a distribution family differing from the one initially considered). To illustrate the \textit{parameter fallacy} from these different viewpoints, consider a collection of data points: $x_1$,$\ldots$,$x_n$. In the Fisher-Neyman-E.~Pearson approach, the data mean could be used to determine a fiducial interval (Fisher) or confidence interval (Neyman-E.~Pearson) on the central value $\mu$ (for fixed $\sigma$) for a Gaussian family of models. No matter the observed distribution of the data, a confidence interval or ``belt'' (e.g., encompassing the central 95\% range) can always be defined based on the value of the data mean (though this glosses over Neyman's paradox\cite{neyman_lectures_1952,redhead_neymans_1974}, which is discussed in \S\ref{sec:neyman} and resolved within the context of the fidelity). However, what if the data are not Gaussian, as Pearson indicates above? The lack of a goodness-of-fit measure means that this question lies outside the domain of the approaches developed by Fisher-Neyman-E.~Pearson; as pointed out above, these approaches in fact can only be defined upon complete restriction of consideration to a specific distribution family. But such ``parametrization of our ignorance'' yields only a very limited approach to statistical inference. Such circumscribed frequentist approaches therefore do not appear to adequately ``do what induction does (namely, to commence a proceeding which must in the long run approximate to the truth)'' (Peirce 2.780)\cite{peirce_collected_1974}. By contrast, goodness-of-fit measured by $\chi^2$ provides a much surer way of determining the absolute concordance of the full model distribution to the complete data. Minimization of $\chi^2$ also provides an optimal method for parameter estimation. However, $\chi^2$ can only be applied to binned data and works well only asymptotically ($>10$ data points in each bin). If only there were an optimal approach similar to $\chi^2$ but that was valid as well for data sets of \textit{arbitrary} size $n$ (something sought by Pearson as well\cite{pearson_method_1933}). In this manuscript, I will argue that maximum fidelity represents just such an approach.

The main distinction between Wilson's viewpoint and that of the Fisher-Neyman-E.~Pearson school can be found in their various notions of a confidence interval. Wilson was credited by Neyman\cite{neyman_lectures_1952} as the originator of the concept of confidence intervals based on Wilson's seminal paper ``On Probable Inference'' from 1927\cite{wilson_probable_1927}; however, Wilson demurred for reasons that will become clear below\cite{wilson_comparative_1964}. As this controversy is critical for understanding the viewpoint of the current manuscript, it is necessary to delve a bit deeper. In his 1927 paper, Wilson addresses inference of the true success rate, $p$, of a binary sequence of $n$ events with an observed success ratio of $p_0$. The conventional view was that:
\begin{quote}
``The probability that the true value of the rate $p$ lies outside its limits $p_0-\lambda\sigma_0$ and $p_0+\lambda\sigma_0$ is less than or equal to $P_{\lambda}$.''\cite{wilson_probable_1927}
\end{quote}
$P_{\lambda}$ is a probability that decreases for increasing $\lambda$.
Wilson criticized this view:
\begin{quote}
``Strictly speaking, the usual statement of probable inference as given above is elliptical. Really the chance that the true probability $p$ lies outside a specified range is either 0 or 1; for $p$ actually lies within that range or does not. It is the observed rate $p_0$ which has a greater or less chance of lying within a certain interval of the true rate $p$. If the observer has had the hard luck to have observed a relatively rare event and to have based his inference thereon, he may be fairly wide of the mark.''\cite{wilson_probable_1927}
\end{quote}
Wilson developed his own statistic for describing such binary processes (by ordering the outcomes according to their individual probabilities for an assumed $p$, akin to a \textit{likelihood} approach, albeit a coordinate-independent one due to the discrete nature of the data) and importantly discussed its interpretation as follows:
\begin{quote}
``The rule then may be stated as: If the true value of the probability $p$ lies outside the range \ldots [Wilson states the range of his interval in terms of $n$, $p_0$, and $\lambda$] \ldots the chance of having such hard luck as to have made an observation so bad as $p_0$ is less or equal to $P_{\lambda}$. And this form of statement is not elliptical. It is the proper form of probable inference.''\cite{wilson_probable_1927}
\end{quote}
In 1934, Clopper \& E.~Pearson\cite{clopper_use_1934} published their own interpretation of the confidence interval for such a binary process (Neyman cited their work as confirmation of his more general arguments\cite{neyman_two_1934}).  Clopper-E.~Pearson use the model cumulative distribution to define an exact, central interval over which a confidence level or, more accurately, 2D ``belt'' (e.g. corresponding to 95\%) could be assigned based on the assumption of an independent Bernoulli process for the events\cite{clopper_use_1934}. Wilson criticized E.~Pearson's probabilistic viewpoint in his 1942 paper ``On Confidence Intervals''\cite{wilson_confidence_1942}. Wilson first restated his view:
\begin{quote}
``I was trying to emphasize that we know nothing about the value of $p$, which must have whatever value it did have in the universe from which the sample was drawn, but that we could set limits based on probability calculations such that if $p$ lay between them the chance of getting the particular observation or any less probable one would exceed some preassigned value $P$ whereas if $p$ lay outside them the chance of getting the observation or any less probable one would be less than $P$.''\cite{wilson_confidence_1942}
\end{quote}
Wilson implied that Clopper and E.~Pearson's probabilistic interpretation went too far:
\begin{quote}
``Furthermore Clopper and Pearson (and Rietz) introduce the notion of confidence belts, which I did not have in mind, and state: `We cannot therefore say that for any specified value of $x$($= np_0$) the probability that the confidence interval will include $p$ is 0.95 or more. The probability must be associated with the whole belt, that is to say with the result of the continued application of a method of procedure to all values of $x$ met with in our statistical experience' --- a statement I should hesitate to make.''\cite{wilson_confidence_1942}
\end{quote}
Note that Wilson always speaks of a limit on the chance of obtaining the \textit{observed} $p_0$ based on the assumption of a particular value of $p$, whereas Clopper-E.~Pearson make a statement about the probability of $p$ lying within the confidence belt based on the observed value $p_0$. Based on the implications of both his and their confidence intervals, Wilson concluded:
\begin{quote}
``Thus although confidence intervals are based on probabilities, it is not certain that probabilities are based on them.''\cite{wilson_confidence_1942}
\end{quote}
Wilson's hesitation to conclude anything about the probability distribution of $p$ was prescient given the modern nuanced distinction between the confidence level and the coverage probability\cite{brown_interval_2001}. In 1964, Wilson deflected Neyman's praise for developing the notion of the confidence interval on the following basis:
\begin{quote}
``Long ago I suggested not making any guess as to the value of an unknown [binary] probability [parameter], but rather to consider what could be inferred from it or what value it would have to have to make some function of it take some assigned value.''\cite{wilson_comparative_1964}
\end{quote}
And in the footnote to this statement:
\begin{quote}
``Mr. Neyman, on page 222 of his \textit{Lectures and Conferences on Mathematical Statistics and Probability} (USDA, 1952), suggested that I then had the idea of confidence intervals. I would make no such claim. I was merely trying to say that for logical reasons statisticians should treat an unknown probability as unknown and in using a standard deviation should allow for any Lexian ratio that rates of the sort they were considering might be expected to show.''\cite{wilson_comparative_1964}
\end{quote}
Wilson is here implicitly criticizing Neyman's (and E.~Pearson's) viewpoint in two ways. First, by unequivocally stating that ``for logical reasons statisticians should treat an unknown probability as unknown'' (here, ``probability'' clearly refers to the Bernoulli \textit{parameter} $p$, not the probability distribution of $p$, but Wilson would have also agreed on the latter interpretation) and, second, by stressing the necessity of avoiding the \textit{parameter fallacy}, in this case the restriction of consideration to a Bernoulli process. Binary data obtained in the real world might be drawn from a Bernoulli process, but they could also have been drawn from a different process, e.g. one that deviates from a Bernoulli process as judged by the ``Lexian ratio''\cite{lexis_zur_1877}. Wilson considered these alternatives at the end of his 1927 paper:
\begin{quote}
``It is well known that some phenomena show less and some show more variation than that due to chance as determined by the Bernoulli expansion $(p+q)^n$. The value $L$ of the Lexian ratio is precisely the ratio of the observed dispersion to the value of $(npq)^{1/2}$ or $(pq/n)^{1/2}$ as the case may be. If we have general information which leads us to believe that the variation of a particular phenomenon be supernormal ($L> 1$), we naturally shall allow for some value of $L$ in drawing the inference. Thus if the Lexian ratio is presumed from previous analysis of similar phenomena to be in the neighborhood of 5, we may use $\lambda=10$ as properly as we should use $\lambda=2$ if the phenomenon were believed to be normal (Bernoullian).''\cite{wilson_probable_1927}
\end{quote}
Wilson here leaves open the possibility of an underlying distribution from an \textit{alternative} family accounting for the data as well, which precludes making any definite fundamental statement about the probability of the overall \textit{number} of successes to lie within certain bounds in the sense of Neyman-E.~Pearson. The inference of a binary success rate from a sequence of $n$ events is given a novel and more exact interpretation in \S\ref{sec:binary} within the context of the fidelity and the cumulative distribution. This fidelity-based approach is actually mathematically more related to Clopper-E.~Pearson's\cite{clopper_use_1934} through its use of the cumulative distribution, but is philosophically in line with Wilson's restrictive interpretation of the confidence interval.

Peirce similarly understood that the quantities obtained through induction should not be interpreted as normal probabilities:
\begin{quote}
``The theory here proposed does not assign any probability to the inductive or hypothetic conclusion, in the sense of undertaking to say how frequently that conclusion would be found true.  It does not propose to look through all the possible universes, and say in what proportion of them a certain uniformity occurs; such a proceeding, were it possible, would be quite idle.  The theory here presented only says how frequently, in this universe, the special form of induction or hypothesis would lead us right.  The probability given by this theory is in every way different--- in meaning, numerical value, and form --- from that of those who would apply to ampliative inference the doctrine of inverse chances.'' (Peirce 2.748)\cite{peirce_collected_1974}
\end{quote}
Peirce defines this ``probability'' quantity more exactly in the following passage:
\begin{quote}
``The third order of induction, which may be called Statistical Induction, differs entirely from the other two in that it assigns a definite value to a quantity. It draws a sample of a class, finds a numerical expression for a predesignate character of that sample and extends this evaluation, under proper qualification, to the entire class, by the aid of the doctrine of chances. The doctrine of chances is, in itself, purely deductive. It draws necessary conclusions only. The third order of induction takes advantage of the information thus deduced to render induction exact.'' (Peirce 7.120)\cite{peirce_collected_1974}
\end{quote}
Peirce based his views largely upon the similar consideration of inference on binary data sets (e.g. Peirce 2.687)\cite{peirce_collected_1974}, which certainly influenced Wilson's formulation of ``a proper form of probable inference''\cite{wilson_statistical_1926,wilson_probable_1927}.

The approach presented here, called maximum fidelity, avoids the \textit{probability fallacy} and the \textit{parameter fallacy}  --- the latter the mathematical \textit{crutch} of both Bayesian and most conventional frequentist approaches to inference --- as well as all of the other fallacies or \textit{ad hoc} limitations mentioned above. Maximum fidelity represents a literal interpretation of Peirce's and Wilson's guidelines, which have heretofore not been fully appreciated. Maximum fidelity ``sets out with a theory'' (a particular hypothesized probability distribution) and ``measures the degree of concordance of that theory with fact'' (Peirce 5.145)\cite{peirce_collected_1974}, obtaining the fidelity and associated concordance as universal ``predesignate characters'' of the data set of arbitrary size $n$. The fidelity is a seemingly optimal information-theoretic statistic defined by the model-based cumulative mapping of the observed data points and important symmetry considerations. Its basis in the cumulative mapping importantly makes it independent of the choice of coordinate system for the data (avoiding the \textit{coordinate fallacy}). Maximization of the fidelity yields the particular candidate model that is most concordant with the data, viz., the model that best represents the complete ``information'' contained in the data. Maximum fidelity unifies several typically distinct notions, including Peirce's foundational work on scientific and statistical inference, the cumulative distribution, information theory, general mathematical considerations (symmetry, boundaries, dimensionality), optimal parameter estimation, and the central notion of model concordance. Maximum fidelity also provides a clear setting for the application of Ockham's razor, which together with model concordance are the dual, competing considerations upon which all scientific inference is based.

\begin{figure}[]
   \begin{center}
\vspace{-4.5cm}
     \includegraphics*[width=6.in]{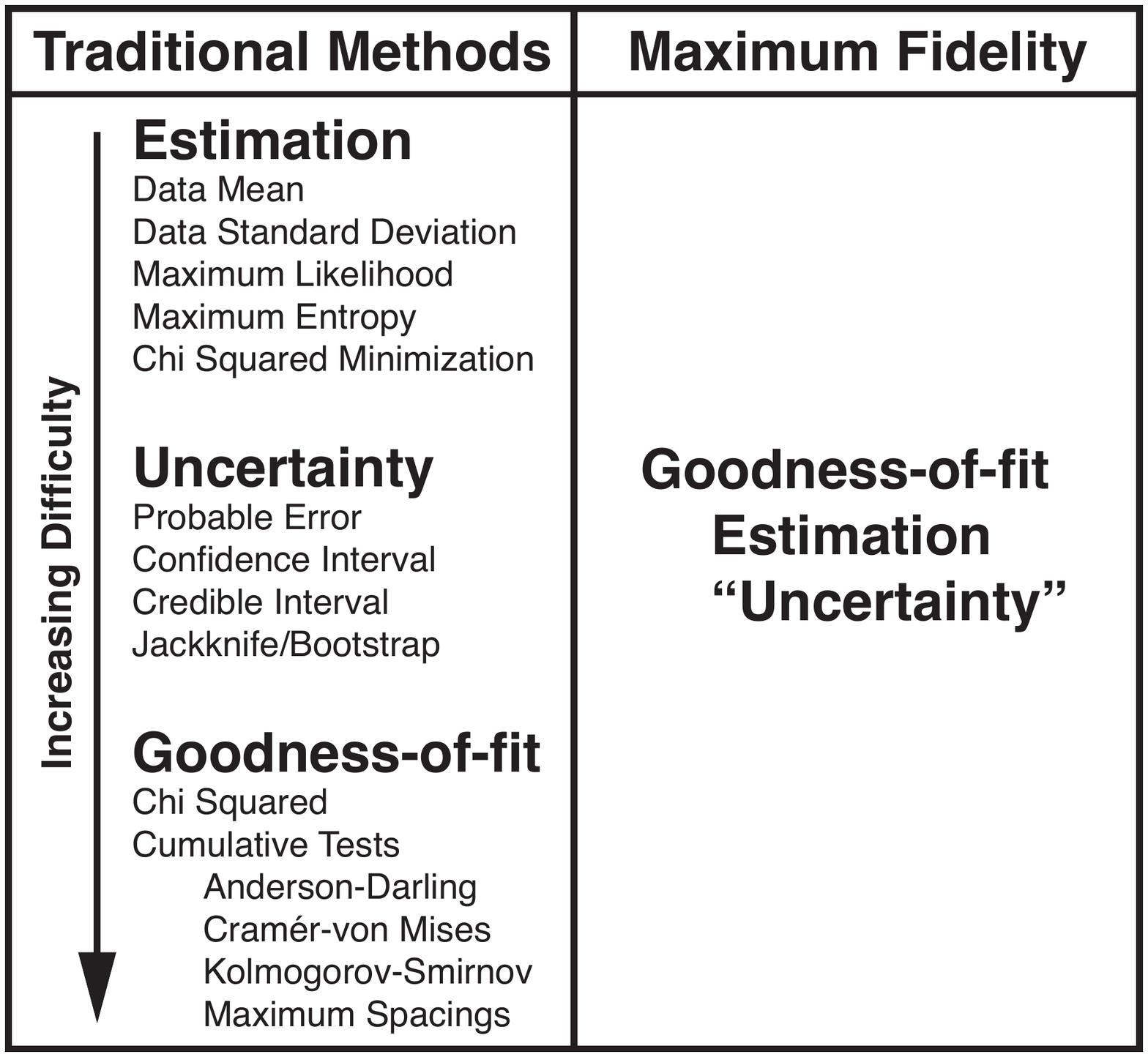}
\vspace{-4.5cm}
      \caption{Comparison of traditional statistical approaches with maximum fidelity.}
      \label{fig:statover}
   \end{center}
\end{figure}

The principal differences between maximum fidelity and traditional approaches are displayed in Fig.~\ref{fig:statover}. In traditional approaches, a heirarchy of ``difficulty'' (or ``ambiguity'' or ``ill-definedness'') is present, with parameter estimation the easiest, parameter uncertainty more difficult, and goodness-of-fit only possible for a very restricted class of problems. At each level, completely different frameworks, assumptions, and artificial limitations are typically required. Maximum fidelity subverts this heirarchy through its fundamental, unifying basis in model concordance (goodness-of-fit), with parameter estimation and a heavily qualified form of parameter ``uncertainty'' following naturally.

The simplest route to understanding maximum fidelity is to see how the fidelity statistic is derived, how the method is applied, and how well it performs in comparison with other approaches.
The outline of this manuscript is therefore as follows. In \S\ref{sec:fidelity}, a general derivation of the fidelity and the related fidelity statistic for univariate data on the circle and on the line is presented. In \S\ref{sec:estimation}, the superiority of maximum fidelity over all other methods --- including the ``gold standard'' method of maximum likelihood --- for parameter estimation on the circle and the line is demonstrated. In \S\ref{sec:concordance}, I show how the fidelity can be converted directly to an absolute concordance value ($p$-value, derived from the null hypothesis) through use of a highly accurate gamma-function approximation. In this section, I also argue for the \textit{general} superiority of this fidelity-based concordance value to other classical goodness-of-fit measures. In \S\ref{sec:joint}, joint analysis of independent data sets with maximum fidelity is explained. In \S\ref{sec:tests}, fidelity-based generalizations of Student's $t$ test and related tests are presented. In \S\ref{sec:neyman}, Neyman's paradox is resolved within the context of the fidelity. In \S\ref{sec:binned}, the straightforward extension of maximum fidelity to binned data is given. In \S\ref{sec:binary}, a solution to the classical problem of binary distributions within the context of the fidelity is given. In \S\ref{sec:multidim}, the application of maximum fidelity to higher dimensional data sets is demonstrated based on ``inverse Monte Carlo'' reasoning. In \S\ref{sec:nonparam}, an extension of maximum fidelity to the nonparametric (or, more accurately, model-independent) determination of whether two observed data sets were drawn from the same unknown distribution is proposed and shown to be \textit{generally} superior to other classical tests. In \S\ref{sec:discussion}, I conclude with a discussion of these results, in particular as they pertain to general scientific inference.

%%%%%%%%%%%%%%%%%%%%%%%%%%%%%%%%%%%%%%%%%%%%%%%%%%%%%%%%%%%%%%
\section{Derivation of the Fidelity}\label{sec:fidelity}
%%%%%%%%%%%%%%%%%%%%%%%%%%%%%%%%%%%%%%%%%%%%%%%%%%%%%%%%%%%%%%

The information measure I refer to as the \textit{fidelity} was first introduced by Shannon\cite{shannon_mathematical_1949} and Wiener\cite{wiener_cybernetics_1948} and later further examined by Kullback and Leibler\cite{kullback_information_1951}, after whom this measure is generally referred (e.g., the Kullback-Leibler divergence or Kullback-Leibler relative entropy). The fidelity, in fact, can be considered the most \textit{fundamental} information measure, from which other measures can be derived. (It is important to note that my use of the term ``fidelity'' should be distinguished from Shannon's use of the same term for a different, non-unique quantity\cite{shannon_mathematical_1949}.)

My choice of the word \textit{fidelity} to describe this information measure best reflects its role in comparing how well a specific distribution represents or \textit{translates} the information contained in a reference distribution. The fidelity can best be intuited from the example of text translation. The fidelity of translation from English to Greek will have one value determined by the degree of degeneracy of the translated words (the multiple Greek words that can be substituted for a given word in English, as well as the occasional need to use the same Greek word to substitute for different English words). The fidelity will have a different value when translating from Greek to English, which exhibits well the importance of the \textit{directionality} of the fidelity. In the simplest approach, the fidelity could be based on the substitution degeneracy of all the words in the translation dictionaries (for a particular translation direction), with equal ``weight'' given to each word. A more informative measure would come from assigning different weights for the words dependent on their frequency of usage in each language. The fidelity could also hypothetically be applied to a particular work (e.g., translating Shakespeare's \textit{Hamlet} to Greek, or translating Homer's \textit{Iliad} to English) to attempt to quantify the translation fidelity or, in the opposite sense, what might have been ``lost in translation''. The example of translation should not be interpreted too literally, though, as different words often have only slightly different nuances of meaning that cannot be expressed quantitatively, which severely complicates the task of assigning a single meaningful value for the fidelity. Translation nevertheless serves as a useful intuitive guide to the directional quantity represented by the fidelity.

\begin{figure}[]
   \begin{center}
\vspace{-2.2cm}
     \includegraphics*[width=6.1in]{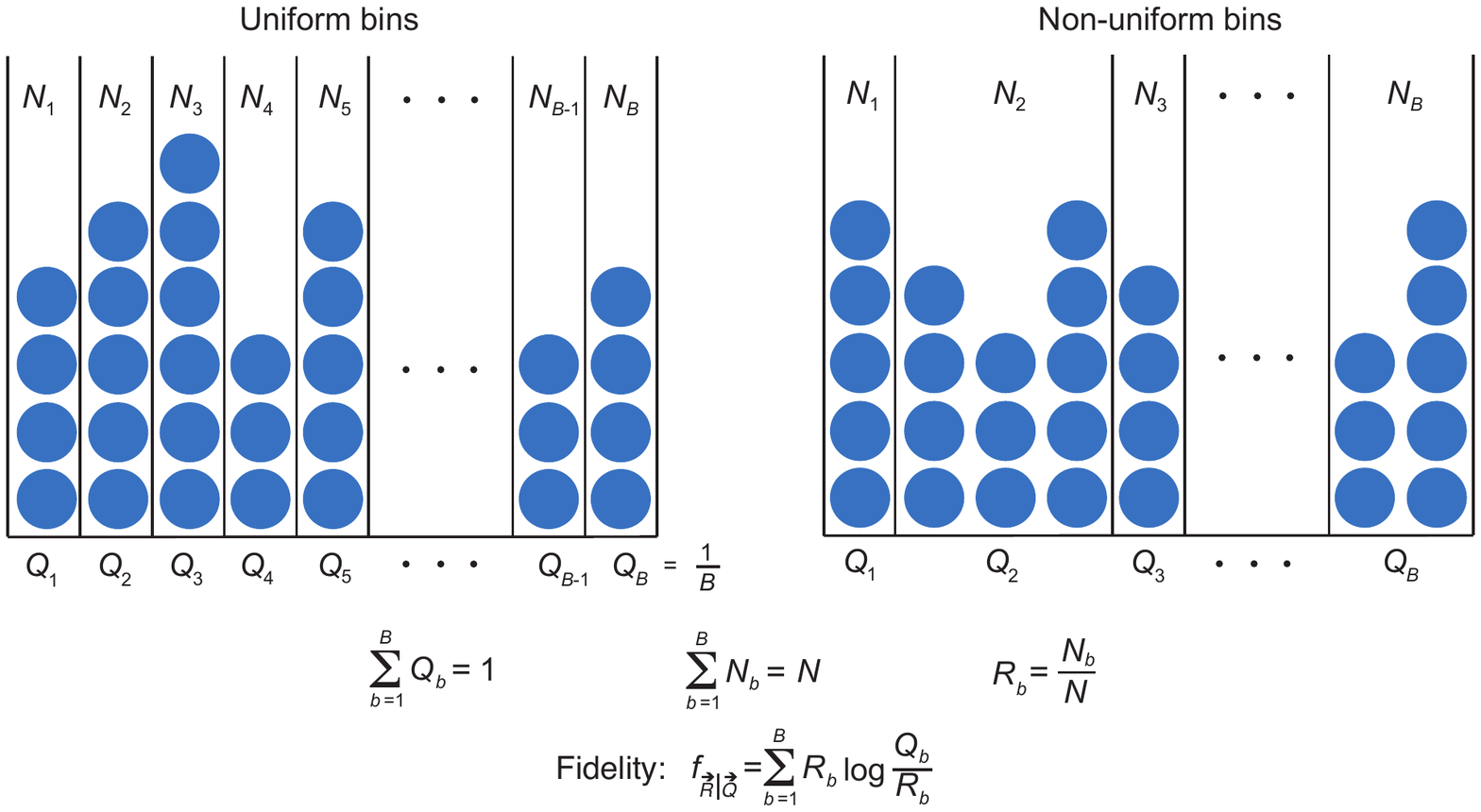}
\vspace{-11.8cm}
      \caption{The probability that a ball lands into a particular bin $b$ out of $B$ total bins is denoted by $Q_b$ with the number of balls in that bin given by $N_b$. The bins could be considered equivalent in width with a uniform probability $Q_b=1/B$ (left), or of different sizes and therefore probabilities $Q_b$ (right). For display purposes, I have chosen to give the balls a ``width'' similar to the width of the bins; the balls should in actuality be considered infinitesimal.}
      \label{fig:balls}
   \end{center}
\end{figure}

Formulas for the entropy or related information measures are often presented in a mysterious or purely mathematical manner. To avoid confusion, it is important to maintain a direct intuitive link to what these measures actually represent (which is accomplished by sticking to definitions based on discrete distributions). Most information measures can be derived from considering the distribution of balls into slots, as shown in Fig.~\ref{fig:balls}. In information theory, these balls can be considered as infinitesimal units of probability. In statistical mechanics, these balls can be considered as infinitesimal probability units over the phase space, as infinitesimal units of energy partitioned across all degrees of freedom of the system, or indeed as literal particles falling into different categories.

Consider $N$ identical balls distributed across $B$ bins with $N_b$ balls in each bin $b$.  We would like to calculate the frequency with which a particular distribution of the $N$ balls $N_1$, $N_2$, $\ldots$, $N_B$ will occur given the assumed known probabilities $Q_b>0$ ($\sum_{b=1}^B Q_b=1$) for a single dropped ball to land in each bin.  For bins of uniform width, $Q_b=1/B$ (Fig.~\ref{fig:balls}, left); however, we will assume that the bins can have arbitrary widths (Fig.~\ref{fig:balls}, right).  The frequency of occurence of a particular distribution $N_1$, $\ldots$, $N_B$ can be determined from the corresponding term in the multinomial expansion of the $N$th power of the sum of probabilities for a single ball:
\begin{eqnarray}
(Q_1+\ldots+Q_B)^N &=& 1^N  \\
\sum_{N_1,\ldots,N_B}\frac{N!}{N_1!\cdots N_B!}Q_1^{N_1}\cdots Q_B^{N_B} &=& 1,
\end{eqnarray}
as
\begin{equation}
\phi_{Q_1,\ldots,Q_B}^{N_1,\ldots,N_B}=\frac{N!}{N_1!\cdots N_B!}Q_1^{N_1}\cdots Q_B^{N_B}.
\end{equation}
In the limit of large $N_b$, we can take the logarithm and use the De Moivre-Stirling approximation:
\begin{equation}
\log N! %=\log \Gamma(N+1)
 \simeq N \log N - N+ \frac{\log (2\pi N)}{2}  +\mathcal{O}(1/N),
\end{equation}
to obtain:
\begin{eqnarray}
%\log \phi_{Q_1,\ldots,Q_B}^{n_1,\ldots,n_B} &=& N\log N - N + \frac{1}{2}\log N + \mathcal{O}(1)+\nonumber\\
%&&\sum_{i=1}^B \left[-n_b\log n_b + n_i - \frac{1}{2}\log n_b + \mathcal{O}(1) + n_b \log Q_b\right]\\
%&=& N\log N  + \frac{1}{2}\log N + \mathcal{O}(1)+\nonumber\\
%&&\sum_{i=1}^B \left[-n_b\log n_b  - \frac{1}{2}\log n_b + \mathcal{O}(1) + n_b \log Q_b\right]\\
\log \phi_{Q_1,\ldots,Q_B}^{N_1,\ldots,N_B} &=& N\log N - N + \frac{\log (2\pi N)}{2} + \mathcal{O}(1/N)+\nonumber\\
&&\sum_{b=1}^B\left[-N_b\log N_b + N_b  - \frac{\log (2\pi N_b)}{2}- \mathcal{O}(1/N_b) + N_b \log Q_b\right]\nonumber\\
&=& N\log N  +   \frac{\log (2\pi N)}{2} + \sum_{b=1}^B \left[ -N_b\log N_b  + N_b \log Q_b- \frac{\log (2\pi N_b)}{2} -  \mathcal{O}(1/N_b)\right] \nonumber\\
&=& N\sum_{b=1}^B \left[-\frac{N_b}{N}\log \frac{N_b}{N} + \frac{N_b}{N} \log Q_b\right] +   \frac{\log (2\pi N)}{2} - \frac{1}{2}\sum_{b=1}^B \left(\log (2\pi N_b) -     \mathcal{O}(1/N_b)\right)\nonumber\\
&\simeq& N\sum_{b=1}^B \left[\frac{N_b}{N}\log \frac{Q_b}{N_b/N}\right]    +   \frac{\log (2\pi N)}{2} - \frac{1}{2} \sum_{b=1}^B \log (2\pi N_b).
\end{eqnarray}
In the last line, we neglect the $\mathcal{O}(1/N_b)$ terms (limit of large $N_b$).  
Defining $R_b=N_b/N$, we are interested in the ratio, $\zeta$, of the probability of a given distribution to the distribution obtained for $R_b=Q_b$ for all $b$, which amounts to subtraction of the logarithm of this distribution in the below expression:
\begin{eqnarray}
\log \zeta_{Q_1,\ldots,Q_B}^{R_1,\ldots,R_B;N} &=& \log \phi_{Q_1,\ldots,Q_B}^{R_1,\ldots,R_B;N} - \log \phi_{Q_1,\ldots,Q_B}^{Q_1,\ldots,Q_B;N}\nonumber\\
 & \simeq & N \sum_{b=1}^B R_b\log{\frac{Q_b}{R_b}} - \frac{1}{2}\sum_{b=1}^B \log\left(2\pi\frac{R_b}{Q_b}\right).
\end{eqnarray}
Assuming $R_b$ and $Q_b$ are strictly greater than 0, then, in the $N\rightarrow\infty$ limit, we can neglect the last term, giving simply:
\begin{eqnarray}
\log \zeta_{Q_1,\ldots,Q_B}^{R_1,\ldots,R_B;N} & \simeq & N \sum_{b=1}^B R_b\log{\frac{Q_b}{R_b}}\nonumber\\
& \simeq & N \left(\sum_{b=1}^{B-1} R_b\log{\frac{Q_b}{R_b}} + R_B \log{\frac{Q_B}{R_B}}\right)\nonumber\\
 & \simeq & N \left(\sum_{b=1}^{B-1} R_b\log{\frac{Q_b}{R_b}} + \left(1-\sum_{b=1}^{B-1}R_b\right) \log{\frac{Q_B}{\left(1-\sum_{b=1}^{B-1} R_b\right)}}\right),
\end{eqnarray}
where I have removed the variable $R_B$ using the normalization condition, leaving $B-1$ independent variables.
The collection $R_1$, $\ldots$, $R_B$ that maximizes $\log \zeta_{Q_1,\ldots,Q_B}^{R_1,\ldots,R_B;N}$ (and therefore also $\phi_{Q_1,\ldots,Q_B}^{R_1,\ldots,R_B;N}$) can be found by setting the derivatives with respect to each $R_b$ (from $b=1$ to $B-1$) to 0:
\begin{eqnarray}
\frac{\partial}{\partial R_b} \log \zeta_{Q_1,\ldots,Q_B}^{R_1,\ldots,R_B;N} &=& N \left(\log \frac{Q_b}{R_b} - 1- \log \frac{Q_B}{\left(1-\sum_{b=1}^{B-1}R_b\right)}+ 1\right)\nonumber\\
&=& N \left(\log \frac{Q_b}{R_b} - \log \frac{Q_B}{R_B}\right),
\end{eqnarray}
which gives for each $b$ simply:
\begin{equation}
R_b = \frac{Q_b}{Q_B}R_B.
\label{eq:Rb}
\end{equation}
Plugging these into the normalization condition gives:
\begin{eqnarray}
\sum_{b=1}^B R_b &=& 1\\
\frac{R_B}{Q_B}\sum_{b=1}^B Q_b &=&1\\
R_B&=&Q_B,
\end{eqnarray}
which, according to Eq.~\ref{eq:Rb}, implies that $R_b=Q_b$ for all $b=1$,$\ldots$,$B$.  

That this unique extremum is in fact a global maximum can be proven by examining the Hessian matrix of second derivatives as follows. The diagonal terms of the $(B-1)\times(B-1)$ Hessian of $\log\zeta$ are:
\begin{eqnarray}
\frac{\partial^2}{\partial R_b^2} \log \zeta_{Q_1,\ldots,Q_B}^{R_1,\ldots,R_B;N}  &=& N\frac{\partial}{\partial R_b} \left(\log \frac{Q_b}{R_b}-\log\frac{Q_B}{R_B}\right)\nonumber\\
&=& -N \left(\frac{1}{R_b}+\frac{1}{R_B}\right).
\end{eqnarray}
All off-diagonal terms ($i\neq b$) are equal to:
\begin{eqnarray}
\frac{\partial}{\partial R_i}\frac{\partial}{\partial R_b} \log \zeta_{Q_1,\ldots,Q_B}^{R_1,\ldots,R_B;N}  &=& N\frac{\partial}{\partial R_i} \left(\log \frac{Q_b}{R_b}-\log\frac{Q_B}{R_B}\right)\nonumber\\
&=& -N\frac{1}{R_B}.
\end{eqnarray}
The Hessian can therefore be written as the sum of a diagonal matrix and a uniform matrix in the following way:
\begin{eqnarray}
H &=& -
\left(\begin{array}{ccc}
1 & \cdots & 0 \\
\vdots   &  \ddots & \vdots \\
0    &     \cdots     & 1
\end{array}
\right)
\left(\begin{array}{c}
\frac{N}{R_1}  \\
 \vdots \\
 \frac{N}{R_{B-1}}
\end{array}
\right)
-\frac{N}{R_B} \left(\begin{array}{ccc}
1 & \cdots & 1 \\
\vdots   &  \ddots & \vdots \\
1    &     \cdots     & 1
\end{array}
\right)
.
\end{eqnarray}
That $\log\zeta$ is strictly concave down (not a saddle point) at its unique extremum can be proven by verifying that all eigenvalues of $H$ are negative. Taking the natural basis for the first, diagonal matrix leads to the trivial eigenvalues of:
\begin{eqnarray}
\lambda_1&=&-\frac{N}{R_1}\nonumber\\
&\vdots&\nonumber\\
\lambda_{B-1} &=& -\frac{N}{R_{B-1}},
\end{eqnarray}
which are all negative (since we have assumed all $R_b$ are strictly greater than 0), implying that the first matrix is negative definite.
A natural set of eigenvectors for the second, uniform matrix can be assumed:
\begin{eqnarray}
v_1 &=& (1,\ldots,1)\nonumber\\
v_2 &=& (1,-1,0,\ldots,0)\nonumber\\
v_3 &=& (1,0,-1,0,\ldots,0)\nonumber\\
&\vdots&\nonumber\\
v_{B-1} &=& (1,0,\ldots,0,-1),
\end{eqnarray}
with corresponding eigenvalues:
\begin{eqnarray}
\lambda_1&=&-(B-1)\frac{N}{R_B}\nonumber\\
\lambda_2 &=& 0\nonumber\\
&\vdots&\nonumber\\
\lambda_{B-1 }&=& 0,
\end{eqnarray}
which implies that the second matrix is negative semi-definite. As the sum of a negative definite matrix and a negative semi-definite matrix, the full Hessian matrix $H$ is negative definite. This proves that $\log\zeta$ is strictly concave down at its unique extremum, therefore establishing this point as a global maximum.

At this maximum, $\log\zeta_{Q_1,\ldots,Q_B}^{R_1=Q_1,\ldots,R_B=Q_B;N}=0$ or $\zeta_{Q_1,\ldots,Q_B}^{R_1=Q_1,\ldots,R_B=Q_B;N}=1$. At all other points, $\zeta$ is exponentially suppressed as $\zeta=e^{Nf}$ with
\begin{equation}\label{eq:fidelity}
f=f_{\vec{R}|\vec{Q}}=\sum_{b=1}^B R_b\log{\frac{Q_b}{R_b}}
\end{equation}
in the discrete case.  For the continuous case (the $B\rightarrow\infty$ limit),
\begin{eqnarray}
\log \zeta_{Q_1,\ldots,Q_B}^{R_1,\ldots,R_B;N}&=&\lim_{B\rightarrow\infty}\left(N\sum_{b=1}^B R_b\log{\frac{Q_b}{R_b}}\right)\nonumber\\
&=&\lim_{B\rightarrow\infty}\left(N\sum_{b=1}^B r_bdb\log{\frac{q_bdb}{r_bdb}}\right)\nonumber\\
&=& N \int db\ r(b)\log{\frac{q(b)}{r(b)}},
\end{eqnarray}
implying
\begin{eqnarray}
f=f_{r(b)|q(b)}&=&\int db\ r(b)\log{\frac{q(b)}{r(b)}}.
\end{eqnarray}
In the above, we assume that $r(b)$ and $q(b)$ are sufficiently smooth (which represents an important assumption) to justify conversion of the sum to an integral. 

The fidelity provides a natural ordering of distributions specified by the $R_b$ (or continuous distributions $r(b)$) in terms of those that occur most often ($f$ near 0) to those that occur least often (increasingly negative $f$) given the basis distribution specified by the $Q_b$ (or the continuous distribution $q(b)$). The fidelity importantly groups together differently shaped distributions that have the same ``frequency of occurrence''. As one can observe from the above, it is best to define information measures over the discrete distributions $R_b$ and $Q_b$, as this preserves the meaning of the ``frequency of occurrence'' (the logarithmic terms in the fidelity ultimately derive from approximation of the factorials over \textit{discrete} distributions). This meaning is unfortunately obscured when information measures are defined from the start using continuous functions such as $r(b)$ and $q(b)$, which are perhaps only later discretized.

\begin{figure}[]
   \begin{center}
\vspace{-6cm}
     \includegraphics*[width=7in]{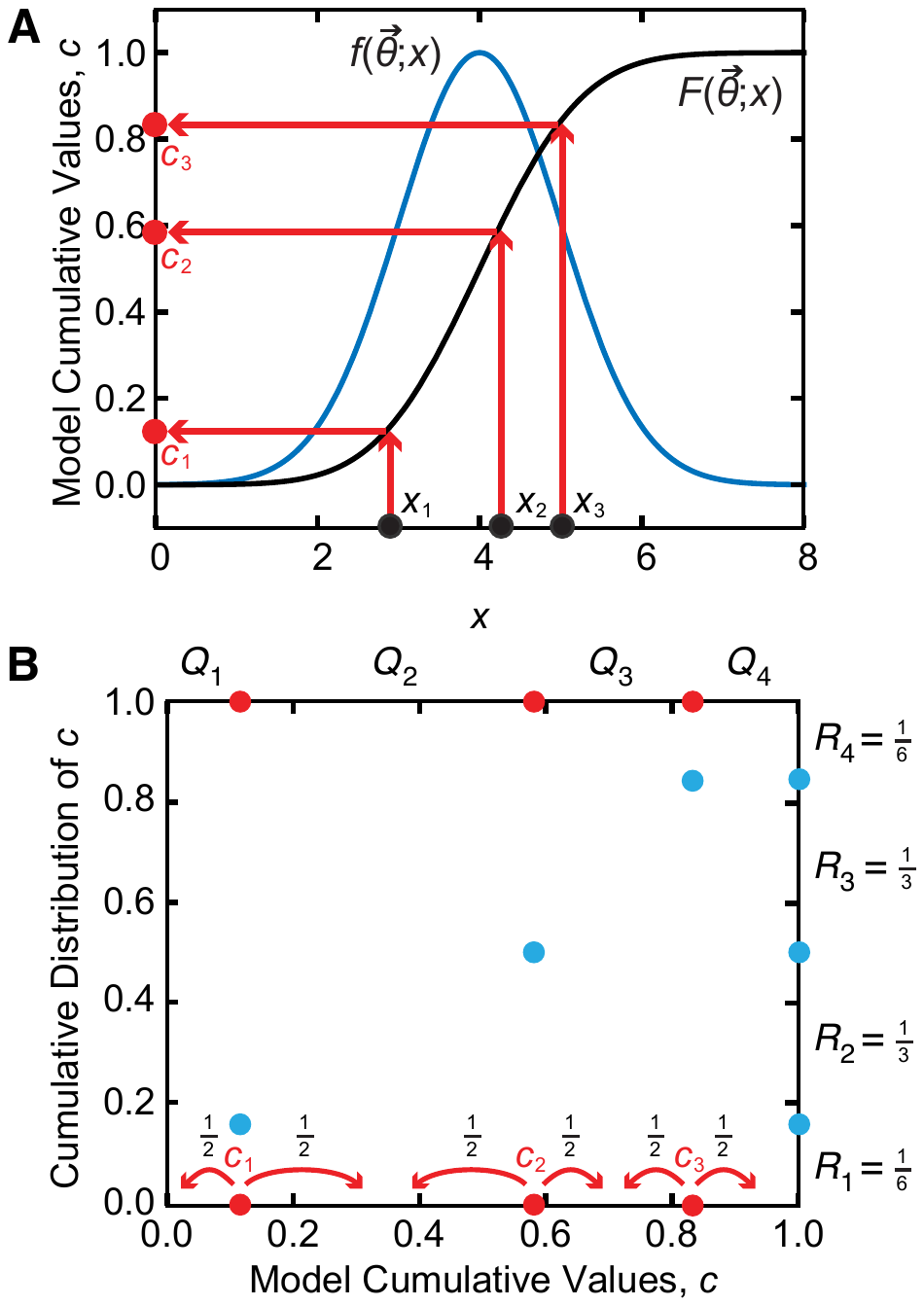}
\vspace{-7cm}
      \caption{(A) Mapping of observed data points to the unit interval by the cumulative distribution (black) of a Gaussian probability distribution with $\mu=4$ and $\sigma=1$ (blue). (B) Determination of $R_b$ from the splitting of the weight of the mapped data points to the adjacent cumulative intervals, with interval sizes, $Q_b$, determined by the model-based cumulative distribution mapping (the $x$-axis here is identical to the $y$-axis in A).}
      \label{fig:cumulative}
   \end{center}
\end{figure}

The fidelity, as derived above, is a completely general information measure. The associated fidelity \textit{statistic} for univariate data on the line can be derived from the following ``density'' argument.  Consider a sorted collection of $n$ points $x_i$, for example, drawn from a Gaussian distribution (Fig.~\ref{fig:cumulative}A). For maximum likelihood, the values of the model probability distribution, $f(x_i)$, corresponding to each observed point $x_i$ are all that are needed to construct the likelihood (which is just the product of all of these values).  Maximum fidelity, on the other hand, is based on the cumulative distribution. The observed points are mapped by the model cumulative distribution to the unit interval (Fig.~\ref{fig:cumulative}A). In order to calculate the fidelity statistic associated with these cumulative values, we assume that each point provides a local estimate of the distribution ``density'' on the cumulative interval.  The only unbiased way to do this (in a way that preserves the local information) is to split each point in half and distribute the resulting weight over the neighboring left and right intervals (Fig.~\ref{fig:cumulative}B). This allows us to determine how much weight should be distributed to each model-defined cumulative interval, based upon which we can calculate the discrete version of the fidelity defined in Eq.~\ref{eq:fidelity}.

The most important thing to note in Fig.~\ref{fig:cumulative}B is that the first and last interval contain only a half point, whereas all interior intervals contain a whole point (the sum of the contributions from both adjacent points). For a given hypothesized model distribution, a mapping of the individual points via the model cumulative distribution generates a series of cumulative values $c_1,\ldots,c_n$ ranging over the unit interval (Fig.~\ref{fig:cumulative}A).  The fidelity statistic is then determined from consideration of the relative weight ($R_b$) of each bin in each of the intervals ($Q_b$) created by the cumulative mapping of the observed data points (Fig.~\ref{fig:cumulative}B):
\begin{eqnarray}
\hspace{-1cm}f&=&\sum_{b=1}^B R_b\log{\frac{Q_b}{R_b}}\nonumber\\
&=& \frac{1}{2n} \log{\frac{c_{1}-0}{\frac{1}{2n}}}  + \frac{1}{2n} \log{\frac{1-c_{n}}{\frac{1}{2n}}}+ \sum_{i=1}^{n-1} \frac{1}{n} \log{\frac{c_{i+1}-c_i}{\frac{1}{n}}} \nonumber\\
&=& \frac{1}{2n} \log{[2n(c_1-0)]} + \frac{1}{2n} \log{[2n(1-c_n)]}+ \frac{1}{n} \sum_{i=1}^{n-1} \log{[n(c_{i+1}-c_i)]}.
\end{eqnarray}
The fidelity statistic can be interpreted in the following way. The $R_b$ values, denoting the relative weights in each bin (either a half point or a full point), should be considered as a particular outcome (e.g. giving the relative numbers of balls in each bin, as in Fig.~\ref{fig:balls}), the frequency of which we wish to ascertain based on the cumulative bin widths determined by the model-derived $Q_b$ values. Maximization of the fidelity therefore identifies the model for which the observed outcome has the highest ``frequency of occurrence'' (under our ``density'' assumption). The most optimal model would place the $n$ data points evenly over the cumulative interval at $c_k=(k-1/2)/n$ (for $k=1,\ldots,n$), yielding the maximum value for the fidelity of $f=0$. 

Note that the only information used by the fidelity derives from the model-based mapping of the discrete events to their cumulative values.  All other information about the function is ignored.  Any other functions having very different shape (large perturbations, highly oscillatory) yet preserving this mapping are indistinguishable in this approach. Candidate distributions that map the data points to the exact same cumulative values form a \textit{mapping equivalence class}. Candidate distributions that generate the same fidelity (and therefore ``frequency of occurrence'') form an even larger \textit{fidelity equivalence class}.

\begin{figure}[]
   \begin{center}
\vspace{-7.8cm}
     \includegraphics*[width=6.5in]{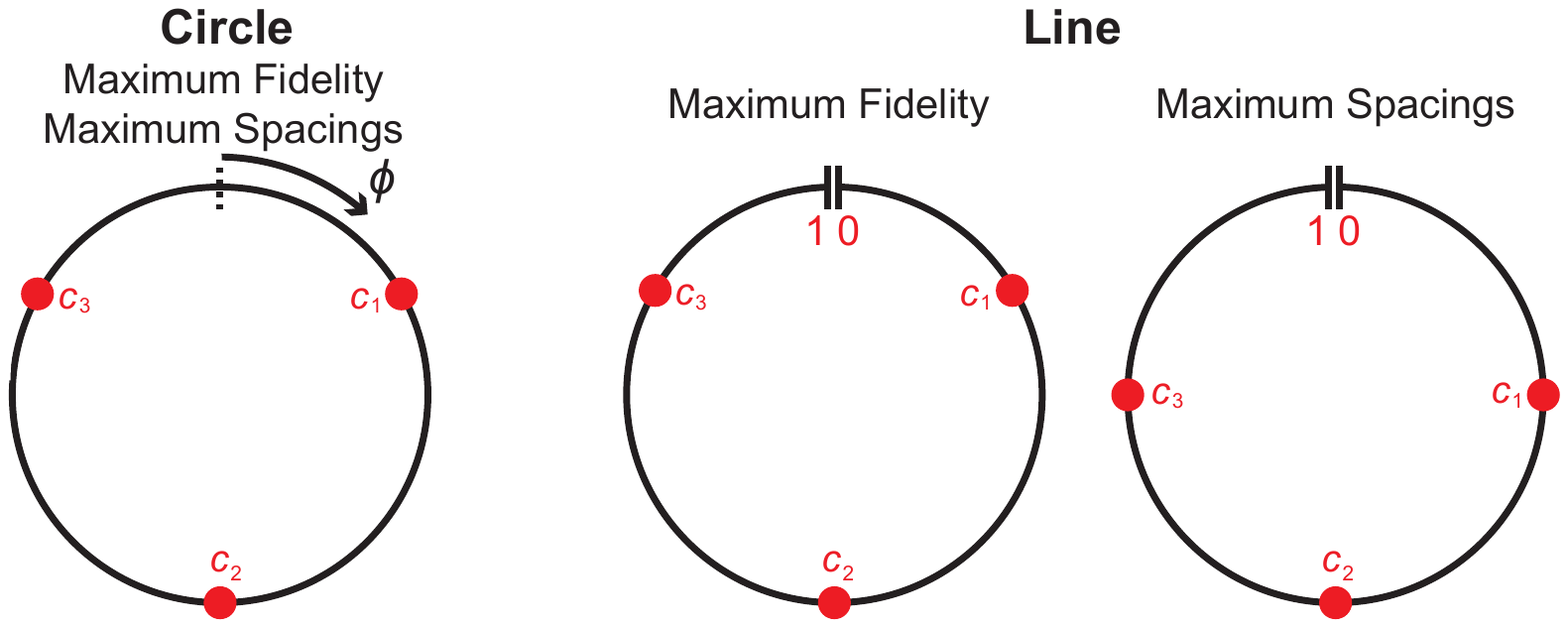}
\vspace{-10.6cm}
      \caption{Optimal cumulative mappings based on maximum fidelity and maximum spacings on the circle and on the line for $n=3$ data points. The $\phi$ parameter indicates the arbitrariness (phase degeneracy) of the location of the coordinate origin on the circle. The line interval has been bent into a circle (with boundaries 0 and 1 indicated at the top) to emphasize the symmetry of the solutions.}
      \label{fig:symmetry}
   \end{center}
\end{figure}

\begin{figure}[]
   \begin{center}
\vspace{-9cm}
     \includegraphics*[width=6.5in]{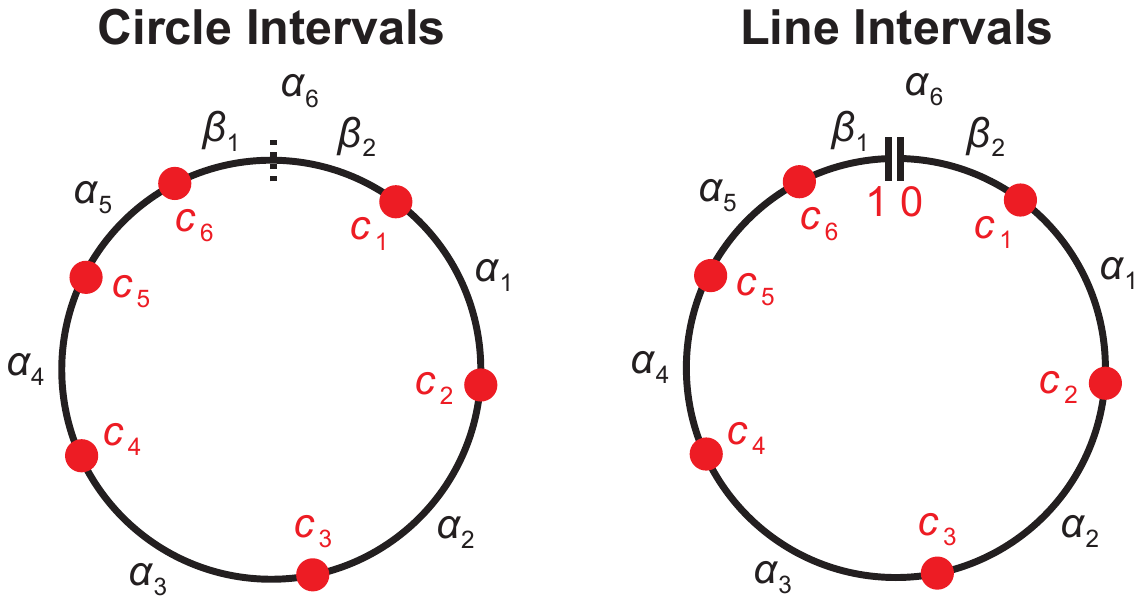}
\vspace{-9.7cm}
      \caption{Cumulative intervals on the circle and on the line for $n=6$ points.}
      \label{fig:intervals}
   \end{center}
\end{figure}

That the boundary intervals contain only a half point in Fig.~\ref{fig:cumulative}B constitutes the only difference between maximum fidelity and the related approach of maximum spacings, which attributes equal weight to all of the intervals. This difference is critical: For the spacings statistic, the intervals are treated as \textit{fundamental}, whereas for the fidelity, the data points are \textit{fundamental}. The maximum spacings statistic for $n$ points on the line (superscript $l$ for \textit{line}) is:
\begin{eqnarray}
s_n^l&=&   \frac{1}{n+1}\log{\left[(n+1) (1-c_{n})\right]} + \frac{1}{n+1}\log{\left[(n+1) (c_{1}-0)\right]} \nonumber\\
&& + \sum_{i=1}^{n-1}  \frac{1}{n+1}\log{\left[(n+1)(c_{i+1}-c_{i})\right]}.
\end{eqnarray}
The spacings statistic has its roots in the work of Moran\cite{moran_random_1951}, who used it to test the concordance of the hypothesis of a Poisson process with an observed data set (time series), but it has been explored by many other researchers since then\cite{darling_class_1953,lecam_theoreme_1958,pyke_spacings_1965,kale_test_1967,gebert_goodness_1969,kale_unified_1969,pyke_spacings_1972,cheng_estimating_1983,ranneby_maximum_1984,cheng_goodness--fit_1989,shao_strong_1999,ekstrom_consistency_2001,anatolyev_alternative_2005,torabi_new_2006,ekstrom_alternatives_2008}. Its connection with information theory has long been recognized\cite{gebert_goodness_1969}. Later investigations have examined the maximization of the spacings statistic as a method for both parameter estimation and general goodness-of-fit assessment\cite{cheng_estimating_1983,ranneby_maximum_1984,cheng_goodness--fit_1989,shao_strong_1999,ekstrom_consistency_2001,anatolyev_alternative_2005,torabi_new_2006,ekstrom_alternatives_2008}.
As derived above, for $n$ data points on the line, the fidelity (with a slight reordering of the boundary terms) is:
\begin{equation}\label{eq:fidline}
f_n^l=   \frac{1}{2n}\log{\left[2n (1-c_{n})\right]} + \frac{1}{2n}\log{\left[2n (c_{1}-0)\right]} + \sum_{i=1}^{n-1}  \frac{1}{n}\log{\left[n(c_{i+1}-c_{i})\right]}.
\end{equation}
For $n$ data points on the circle, the fidelity $f_n^c$, which is identical to the circular version of the spacings statistic $s_n^c$, is:
\begin{equation}\label{eq:fidcircle}
f_n^c= s_n^c= \frac{1}{n}\log{ \left[n\big((1-c_n)+(c_{1}-0)\big)\right]} + \sum_{i=1}^{n-1}  \frac{1}{n}\log{\left[n(c_{i+1}-c_{i})\right]}.
\end{equation}
The difference between the optimal solutions implied by these different statistics is graphically displayed at the top of Fig.~\ref{fig:symmetry} for $n=3$ points. In this figure, the cumulative intervals for the line are bent into a circle to emphasize the non-symmetric weighting of the boundary intervals for maximum spacings. The symmetric set of solutions on the circle (with phase freedom $\phi$) and the single solution shown for the fidelity on the line (which is symmetric with respect to the boundary) can be viewed as an instance of ``symmetry breaking''. A physical theory that breaks symmetry is one that selects a particular solution from a set of symmetric solutions. One can clearly see in Fig.~\ref{fig:symmetry} that the solution of the fidelity on the line belongs to the solution set on the circle, whereas the spacings solution does not correspond to any of the phase-degenerate solutions on the circle. In Fig.~\ref{fig:intervals}, $n=6$ data points are plotted with arbitrary values. The interval at the boundary, $\alpha_6$, created by the first and the last point is afforded the same weight for maximum fidelity as the other intervals both on the circle and on the line. On the circle, $\alpha_6$ is directly considered (as there is no real boundary on the circle); however, on the line, separate consideration is required of the two subintervals, $\beta_1$ and $\beta_2$, which only together contribute a weight equal to the non-boundary intervals.

%%%%%%%%%%%%%%%%%%%%%%%%%%%%%%%%%%%%%%%%%%%%%%%%%%%%%%%%%%%%%%
\section{Parameter Estimation}\label{sec:estimation}
%%%%%%%%%%%%%%%%%%%%%%%%%%%%%%%%%%%%%%%%%%%%%%%%%%%%%%%%%%%%%%

In the following sections, parameter estimation of probability distributions on the circle and on the line using maximum fidelity is compared with other  well-established approaches, including the ``gold standard'' method of maximum likelihood. I will focus in particular on parameter estimation on low $n$ data sets, for which statistical analysis is truly necessary and discrimination among the various estimation methods can be easily visualized. Invocation of the Cram\'er-Rao efficiency to ``prove'' that one estimator is better than another is avoided for the reasons listed in \S\ref{sec:intro}, and is anyway invalid for the low $n$ data sets (typically $n=5$) considered below. Other values for $n$ were tested (not shown) with similar results, proving their generality. Estimate bias with respect to the mean of the parameter estimate distribution is often discussed as a measure of optimality; however, this is completely dependent on the choice of parameter coordinate representation (e.g. $\sigma$ vs. $v=\sigma^2$). One might be tempted to define a good parameter estimate as one with no \textit{median} bias, which at least would be coordinate invariant. However, narrowness of the parameter estimate distribution around the true value is also an important consideration and there is no guarantee that an estimator that is median unbiased will also have the narrowest distribution. It is important to avoid taking too seriously (or quantitatively) the comparisons made below for parameter estimation within a fixed distribution family, as these comparisons are but a form of \textit{parameter fallacy} (see \S\ref{sec:intro}). 
The dual aspects of optimal estimators that nevertheless \textit{will} be focussed on below are low (but possibly non-zero) median bias and a narrow distribution around the true value. We will see that these qualitative considerations are already sufficient to determine the best overall estimation approach.

\subsection{Parameter Estimation on the Circle}

We first consider parameter estimation for circular distributions. Unlike distributions on the line, there are only a handful of circular distributions that are actually interesting to consider as models for real-world data\cite{mardia_directional_1999}. Here, we focus on the three flexible distribution families listed in Table~\ref{tab:circ_dist}. Here, $\beta$ generally denotes a location parameter and $\alpha$ generally represents a shape parameter for the distribution (though $\beta$ is \textit{not} a location parameter for the Wrapped Laplace distribution). In addition to the fidelity and the likelihood, we test the parameter estimation precision of several other circular statistics named after their originators: Ajne \cite{ajne_simple_1968}, Gini\cite{jammalamadaka_test_2004}, Kuiper\cite{kuiper_tests_1962}, Rao\cite{Rao_contributions_1969}, Rayleigh\cite{Rayleigh_resultant_1880,Rayleigh_problem_1905,Rayleigh_xxxi._1919} , and Watson\cite{watson_goodness--fit_1961}. All of these tests, aside from the Gini test, are well described by Mardia \& Jupp\cite{mardia_directional_1999}.

\begin{table}
\begin{center}
\begin{tabular}{|c|c|c|c|}
\hline
%\multicolumn{2}{c}{line} & \multicolumn{2}{c}{circle}  \\
distribution & expression     &   $\beta$  & $\alpha_0$  \\
\hline
%cardiod & $\frac{1}{2\pi}(1+2\alpha\cos{(\theta-\beta)})$ & $\beta=0$   \\
%wrapped Cauchy & $\frac{1}{2\pi}\frac{1-\alpha^2}{1+\alpha^2-2\alpha\cos{(\theta-\beta)}}$  & $\beta=0$     \\
Theta \T& $\frac{1}{2\pi}\frac{\sqrt{1+\alpha}}{1+\alpha\sin^2{(\theta-\beta)}}$ &  0 & 4\\
von Mises \T& $\frac{\exp{(\alpha \cos{(\theta-\beta)})}}{2\pi I_b(0,\alpha)}$  &   0   & 1/4  \\
%Wrapped Exponential \T  & $\frac{\alpha\exp{(-\alpha\theta)}}{1-\exp{(-2\pi\alpha)}}$   &  & 1  \\
Wrapped Laplace \T\B&  $\frac{\beta\alpha}{1+\alpha^2}\left(\frac{\exp{(-\beta\alpha\theta)}}{1-\exp{(-2\pi\beta\alpha)}}  + \frac{\exp{((\beta/\alpha)\theta)}}{\exp{(2\pi\beta/\alpha)}-1}     \right)$  & 1 &  1  \\
\hline
\end{tabular}
\end{center}
\caption{Probability distributions on the circle.}
\label{tab:circ_dist}
\end{table}

\begin{table}
\begin{center}
\begin{tabular}{|c|c|}
\hline
name     &  statistic        \\ 
\hline
Ajne  \T& $\frac{n}{4} - \frac{1}{2\pi n}\sum_{i=1}^n\sum_{j=1}^n \left( \pi-\left| \pi - 2\pi \left| c_i-c_j \right| \,\right| \right)  $ \\
Gini \T& $ 2\sum_{i=1}^{n+1}\sum_{j=i+1}^n \left|(c_j-c_{j-1}) - (c_i-c_{i-1})\right| $\\
Kuiper  \T&   $\max_i\left|c_i-\frac{i-1}{n}\right| + \max_i\left| c_i-\frac{i}{n} \right|    $ \\
Maximum Fidelity  \T& $ \sum_{i=1}^{n}\frac{1}{n}\log[n(c_{i+1}-c_i)]$\\
Maximum Likelihood \T& $\sum_{i=1}^n \log{p(x_i)}$  \\
Rao \T&  $\sum_{i=2}^{n+1}  \left|c_i-c_{i-1}-\frac{1}{n}\right| $   \\
Rayleigh \T& $2n    \left(\sum_{i=1}^n\cos(2\pi c_i)\right)^2      +   2n    \left(\sum_{i=1}^n\sin(2\pi c_i)\right)^2  $\\
Watson \T\B& $\frac{n}{12} + \left(\sum_{i=1}^n c_i^2  \right) - \frac{1}{n}\left( \sum_{i=1}^n c_i\right)^2-\frac{2}{n}\sum_{i=1}^n i c_i + \frac{n+1}{n}\sum _{i=1}^n c_i$    \\
\hline
\end{tabular}
\end{center}
\caption{Statistics on the circle ($c_{n+1} = 1+c_1$ is defined for convenience). For maximum fidelity and maximum likelihood, the statistic is maximized to find the optimal parameter estimate.  For all the other statistics, the minimum is sought.}
\label{tab:circ_stat}
\end{table}

\begin{figure}[!h]
   \begin{center}
\vspace{-1.5cm}
     \includegraphics*[width=5in]{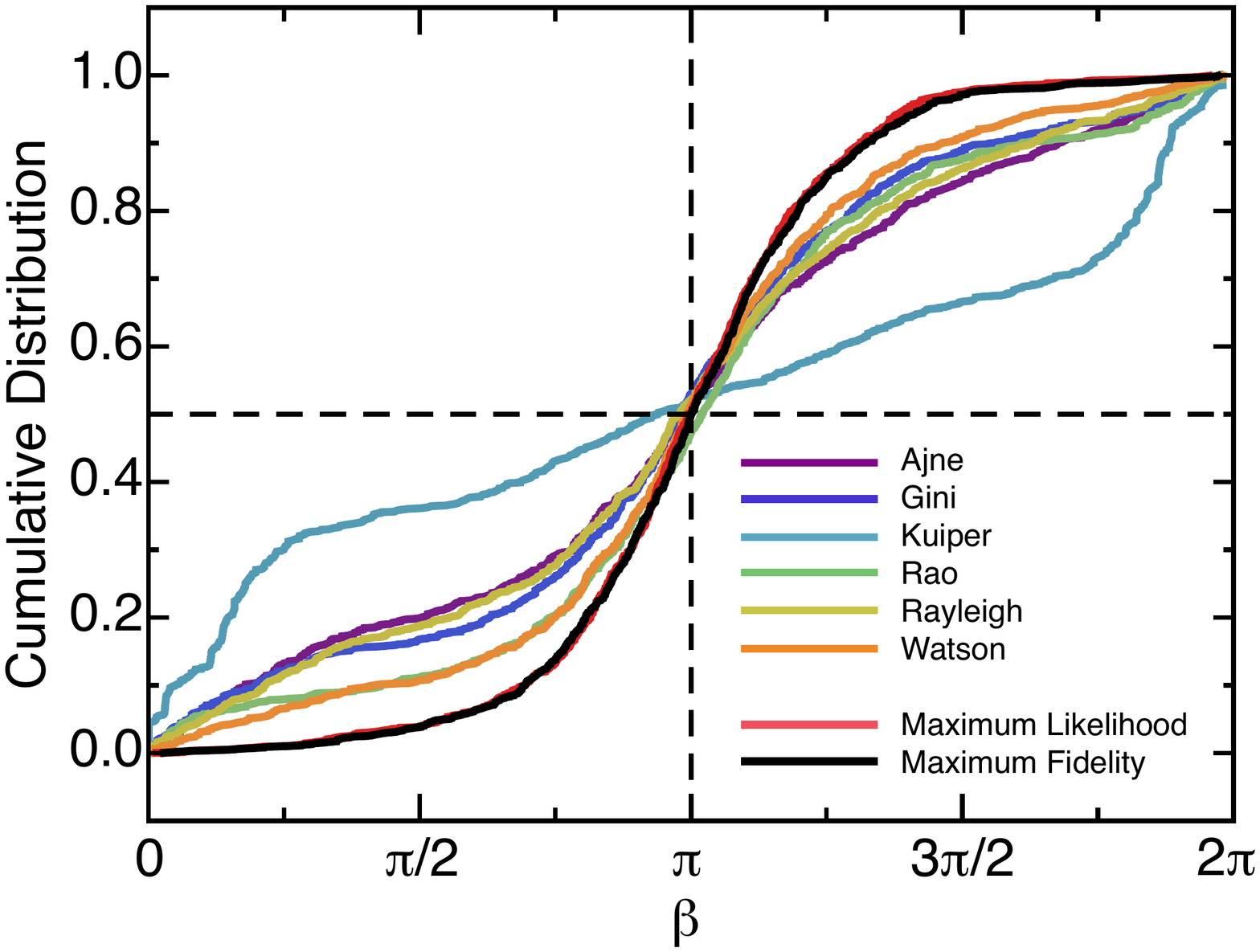}
\vspace{-1cm}
      \caption{Location parameter fitting on the circle for location parameters. For each distribution, $n=5$ data points were repeatedly drawn (1000 times) from a von Mises distribution with $\alpha=1$ and $\beta=\pi$ (see Table~\ref{tab:circ_dist}) by extremizing the various statistics listed in Table~\ref{tab:circ_stat}. To perform this fitting, the cumulative distribution and inverse cumulative functions must be calculated for the true model and the test parameter models. These functions were estimated by calculation on a grid of 101 points (across the full range of $\beta$) followed by Stineman interpolation\cite{wagon_mathematica_2010}. Cumulative distributions for these guesses using the different statistics are displayed, with the vertical dashed line indicating the true value and the horizontal dashed line indicating the median value of the cumulative distribution of guesses. A good estimator should be steep around the true value with low median bias.}
\label{fig:circ_par_loc}
   \end{center}
\end{figure}

We first test the precision of estimation of the location parameter $\beta$ for the von Mises distribution (Fig.~\ref{fig:circ_par_loc}). Here, data sets containing $n=5$ points were repeatedly generated (1000 times) from a von Mises distribution with $\alpha=1$ and $\beta_0=\pi$. Each data set was fit with the indicated statistics to generate the observed cumulative distributions of the best fit $\beta$ values. Due to symmetry considerations, the median value of the fit distribution should asymptotically approach the true value of $\beta_0=\pi$ for all of the statistics (they are all median unbiased), with only the steepness of each fit distribution revealing the precision of the estimator. As can clearly be seen, both maximum fidelity and maximum likelihood provide equivalently narrow (and overlapping) distributions about the true value. All other statistics perform significantly worse.

\begin{figure}[]
   \begin{center}
\vspace{-7cm}
     \includegraphics*[width=6.5in]{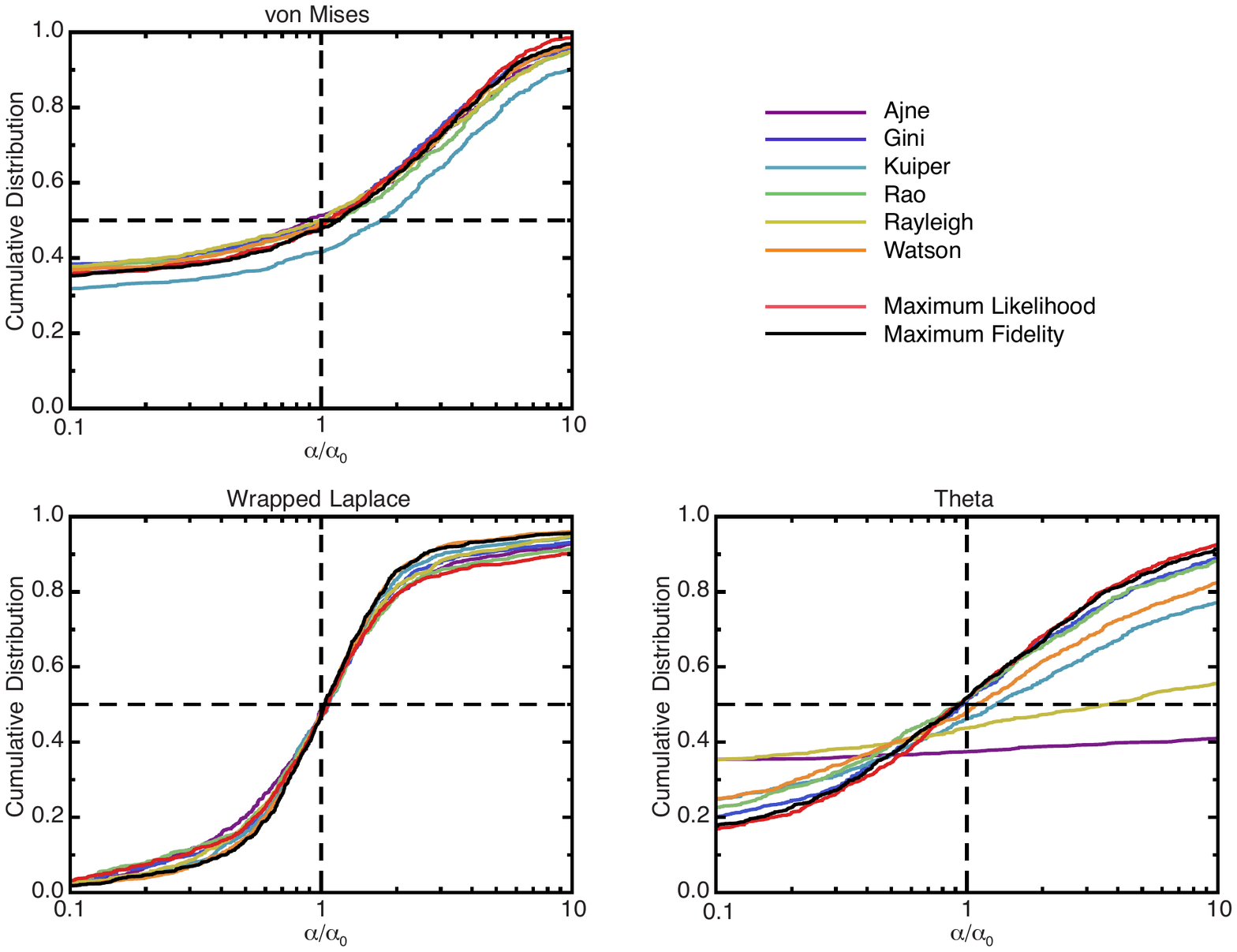}
\vspace{-6.5cm}
      \caption{Shape parameter fitting on the circle. For each distribution, $n=5$ data points were repeatedly drawn (1000 times) with the fit parameter $\alpha$ determined by extremizing the various statistics listed in Table~\ref{tab:circ_stat}. To perform this fitting, the cumulative distribution and inverse cumulative functions must be calculated. These functions were estimated by calculation on a grid of 101 points (evenly sampled over $\log_{10}(\alpha/\alpha_0)$ from -1.5 to 1.5) followed by Stineman interpolation. Cumulative distributions for these estimates using the different statistics are displayed, with the vertical dashed line indicating the true value and the horizontal dashed line indicating the median value of the cumulative distribution of estimates. A good estimator should be steep around the true value with low median bias. }
\label{fig:circ_par}
   \end{center}
\end{figure}

In Fig.~\ref{fig:circ_par}, the more revealing consideration of the estimation of the shape (or width) parameters $\alpha$ for each distribution is displayed (for these fits, $\beta$ is assumed known and the true $\alpha$ is $\alpha_0$, with the exact values taken for each distribution listed in Table~\ref{tab:circ_dist}). For the von Mises distribution, all statistics are equally steep around the true value of $\alpha_0$. All statistics, aside from the Kuiper statistic, also have acceptably low median bias. For the Theta distribution, maximum fidelity and maximum likelihood provide equivalently high precision estimates of the relevant parameter $\alpha$, with the median estimate lying near the true value of $\alpha_0$, flanked by a steep distribution above and below this value. All other distributions clearly have broader wings. For the Wrapped Laplace distribution, maximum fidelity and the Watson statistic provides the best distribution, with maximum likelihood interestingly displaying the worst estimate distribution (this could arise from the ``kink'' at the origin for the Wrapped Laplace, see the red curve in Fig.~\ref{fig:good_circle} for an example). Maximum fidelity, therefore, provides the most optimal shape parameter estimate for all of the tested distributions.

\subsection{Parameter Estimation on the Line}

We now examine the parameter estimation precision for multiple distributions on the line (Table~\ref{tab:line_dist}) of maximum fidelity and other statistics (Table~\ref{tab:line_stat}), including the classical statistics named after their (co-)originators: Anderson-Darling\cite{anderson_asymptotic_1952}, Cram\'er-von Mises\cite{cramer_mathematical_1946,Mises_probability_1928}, Gini\cite{jammalamadaka_test_2004}, and Kolmogorov-Smirnov\cite{kolmogorov_sulla_1933,smirnov_table_1948}. In addition to these, other statistics based on Equal Intervals (the sum of the deviation of each cumulative interval from $1/(n+1)$) and Order Statistics (sum of the logarithm of the order statistics\cite{whitworth_choice_1901,irwin_william_1967,wilks_order_1948}) are also considered.

\begin{table}[h!]
\begin{center}
\begin{tabular}{|c|c|c|c|c|}
\hline
%\multicolumn{2}{c}{line} & \multicolumn{2}{c}{circle}  \\
distribution  & expression    & domain  &  $\beta$  &  $\alpha_0$   \\ 
\hline
Beta \T & $\frac{1}{\mathrm{Beta}(\beta,\alpha)}(1-x)^{-1+\beta}x^{-1+\alpha}$  &  $0\leq x \leq1$ & 3  & 1 \\
Cauchy \T& $\frac{1}{\alpha\pi}\frac{1}{1+(x-\beta)^2/\alpha^2}  $ &   $-\infty<x<\infty$     & 0 & 1/2 \\
Exponential \T& $\alpha e^{-\alpha x}$&   $x\geq 0$   & & 1 \\
Extreme Value \T& $\frac{1}{\beta}e^{-e^{(\beta-x)/\alpha}}e^{(\beta-x)/\alpha}$ & $-\infty<x<\infty$  & 0 & 1 \\
$F$ Ratio \T& $\frac{1}{\mathrm{Beta}(\beta/2,\alpha/2)}\beta^{\beta/2}\alpha^{\alpha/2}x^{\alpha/2-1}(\beta+\alpha x)^{-(\beta+\alpha)/2}$ & $x\geq 0$ & 2 & 1 \\
Gamma \T& $\frac{1}{\Gamma(\beta)}\alpha^{-\beta}x^{\beta-1}e^{-x/\alpha}$  & $x\geq0$ & 3 &  1 \\
Gauss \T& $\frac{1}{\sqrt{2\pi}\alpha}\exp{\left(-\frac{(x-\beta)^2}{2\alpha^2}\right)}$  & $-\infty<x<\infty$ &   0  & 1 \\
Inverse Gamma \T& $\frac{1}{\Gamma(\beta)x}e^{-\alpha/x}(\alpha/x)^{\beta}  $ & $x>0$ & 3   &  1 \\
Laplace \T& $\frac{1}{2\beta}e^{-|x-\beta|/\alpha}$  & $-\infty<x<\infty$ & 0 &  1 \\
Levy \T& $\frac{1}{\sqrt{2\pi}\alpha}e^{-\frac{\alpha}{2(x-\beta)}}\left(\frac{\alpha}{x-\beta}\right)^{3/2}   $ & $x>\beta$  & 0 & 1 \\
Logistic \T& $\frac{1}{\alpha}\frac{e^{(x-\beta)/\alpha}}{\left(1+e^{-(x-\beta)/\alpha}\right)^2}$  & $-\infty<x<\infty$  & 0 & 1 \\
%Maxwell \T& $\frac{1}{\alpha^3}\sqrt{\frac{2}{\pi}}x^2e^{-\frac{x^2}{2\alpha^2}} $  & $x\geq 0$  & & 1 \\
Pareto \T& $\alpha^{\beta}\beta x^{-\beta-1}$  & $x>\beta$ & 1 &  1 \\
Rayleigh \T& $\frac{1}{\alpha^2}x e^{-\frac{x^2}{2\alpha^2}}$ & $x\geq 0$   & & 1 \\
Student  \T& $\frac{1}{\sqrt{\alpha}\,\,\mathrm{Beta}(\alpha/2,1/2)}\left(\frac{\alpha}{x^2+\alpha}\right)^\frac{1+\alpha}{2}$  & $-\infty<x<\infty$ &  & 3 \\
%Wald \T& $\frac{1}{\sqrt{2\pi}}\exp{\left(-\frac{\alpha (x-\beta)^2}{2x\beta^2}\right)}\sqrt{\frac{\alpha}{x^2}}$  & $x>0$  & 0 & 1 \\
Weibull \T\B & $\frac{\alpha}{\beta}(x/\beta)^{\alpha-1}e^{-(x/\beta)^{\alpha}}$  & $x>0$  & 1 & 2 \\
\hline
\end{tabular}
\end{center}
\caption{Probability distributions on the line.}
\label{tab:line_dist}
\end{table}

\begin{table}
\begin{center}
\begin{tabular}{|c|c|}
\hline
name     &  statistic       \\ 
\hline
Anderson-Darling
 \T& $2\sum_{i=1}^n\frac{i-1/2}{n} \left(\log c_i + \log(1-c_{n-i+1})    \right) $  \\
Cram\'er-von Mises \T& $\frac{1}{12n} + \sum_{i=1}^n \left(\frac{i-1/2}{n}-c_i   \right)^2$  \\
Equal Intervals \T&  $\sum_{i=1}^{n+1}  \left|c_i-c_{i-1}-\frac{1}{n+1}\right| $  \\
Gini \T& $ 2\sum_{i=1}^{n+1}\sum_{j=i+1}^n \left|(c_j-c_{j-1}) - (c_i-c_{i-1})\right| $      \\   
Kolmogorov-Smirnov \T& $\max\left(\max_i \left|  c_i-\frac{i-1}{n}\right|, \max_i \left| c_i-\frac{i}{n}  \right| \right)$  \\
Maximum Fidelity \T& $\frac{1}{2n}\log[2n(1-c_n)]  +  \frac{1}{2n}\log[2n(c_1-0)]  + \sum_{i=1}^{n-1} \frac{1}{n}\log[n(c_{i+1}-c_i)]$  \\
Maximum Likelihood \T& $\sum_{i=1}^n \log{p(x_i)}$    \\
Maximum Spacings \T& \hspace{-2cm}$\frac{1}{n+1}\log[(n+1)(1-c_n)]  +  \frac{1}{n+1}\log[(n+1)(c_1-0)]$  \\
& \hspace{+4.5cm}$+ \sum_{i=1}^{n-1} \frac{1}{n+1}\log[(n+1)(c_{i+1}-c_i)]$   \\
Order Statistics \T\B& $\sum_{i=1}^n \big( (i-1)\log(c_i-0)+(n-i)\log(1-c_i)    \big)$   \\
\hline
\end{tabular}
\end{center}
\caption{Statistics on the line ($c_0=0$ and $c_{n+1}=1$ are defined for convenience). For the Cram\'er-von Mises, Equal Intervals, Gini, and Kolmogorov-Smirnov statistics, the minimum is considered optimal.  For all the other statistics, the maximum is considered optimal.}
\label{tab:line_stat}
\end{table}

\begin{figure}[]
   \begin{center}
\vspace{-1cm}
     \includegraphics*[width=4.0in]{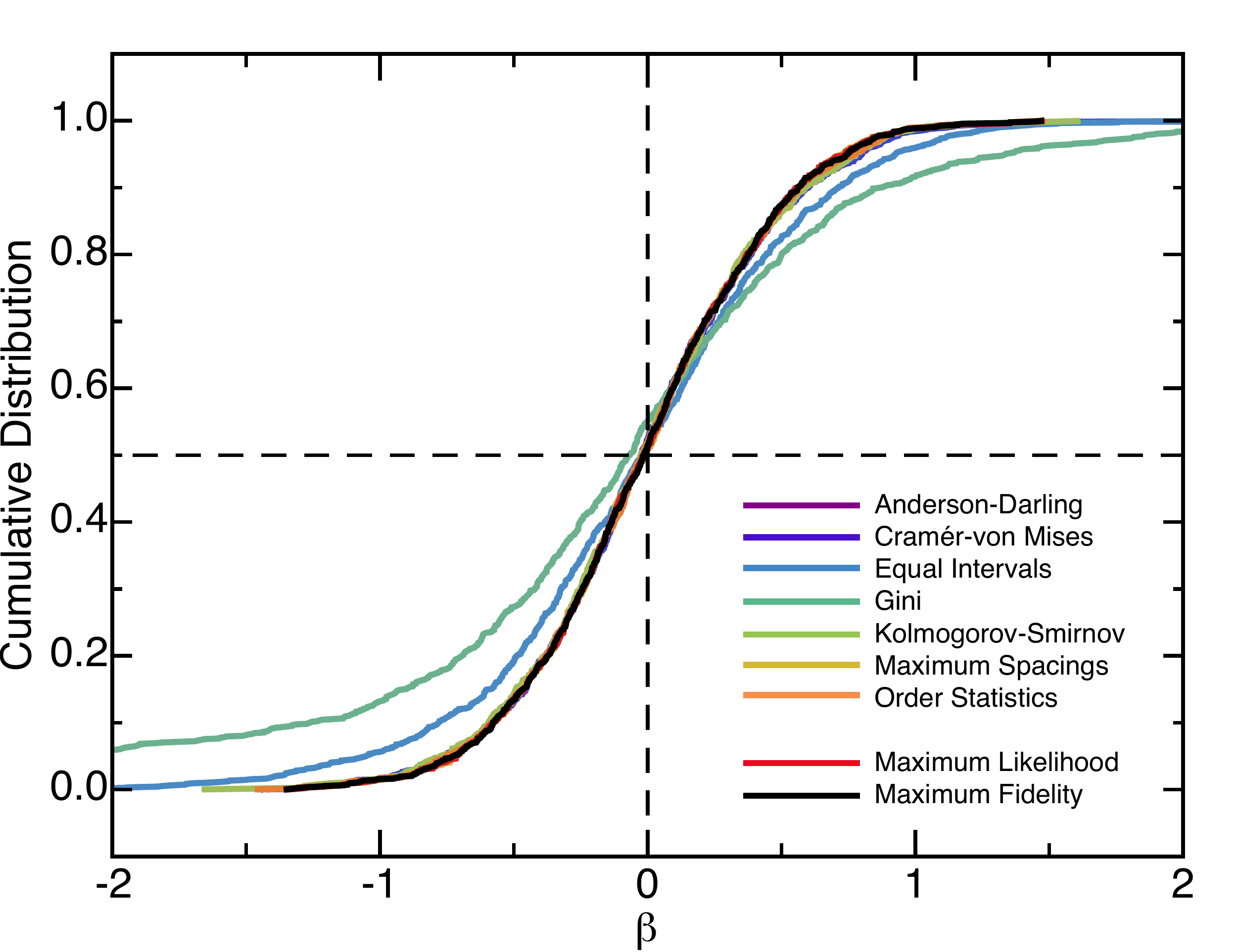}
\vspace{0cm}
      \caption{Fitting of a location parameter for a Gaussian on the line. For each model, $n=5$ data points were drawn repeatedly from a Gaussian with $\mu=0$ and $\sigma=1$. Fitting of the Gaussian mean $\mu$ (assuming a fixed $\sigma=1$)  was carried out for each realization (1000 total) using the listed statistics. All listed statistics lead to overlapping distributions except for the Equal Intervals and Gini statistics.}
      \label{fig:par_line_loc}
   \end{center}
\end{figure}

In Fig.~\ref{fig:par_line_loc}, a ``location'' test is performed: The estimation precision for the mean value of a Gaussian (with width $\sigma=1$ assumed known) is compared across the statistics listed in Table~\ref{tab:line_stat}. All statistics, excluding the Equal Intervals statistic and the Gini statistic, have equal precision for this ``location'' test.

\begin{figure}[]
   \begin{center}
\vspace{0cm}
     \includegraphics*[width=6.5in]{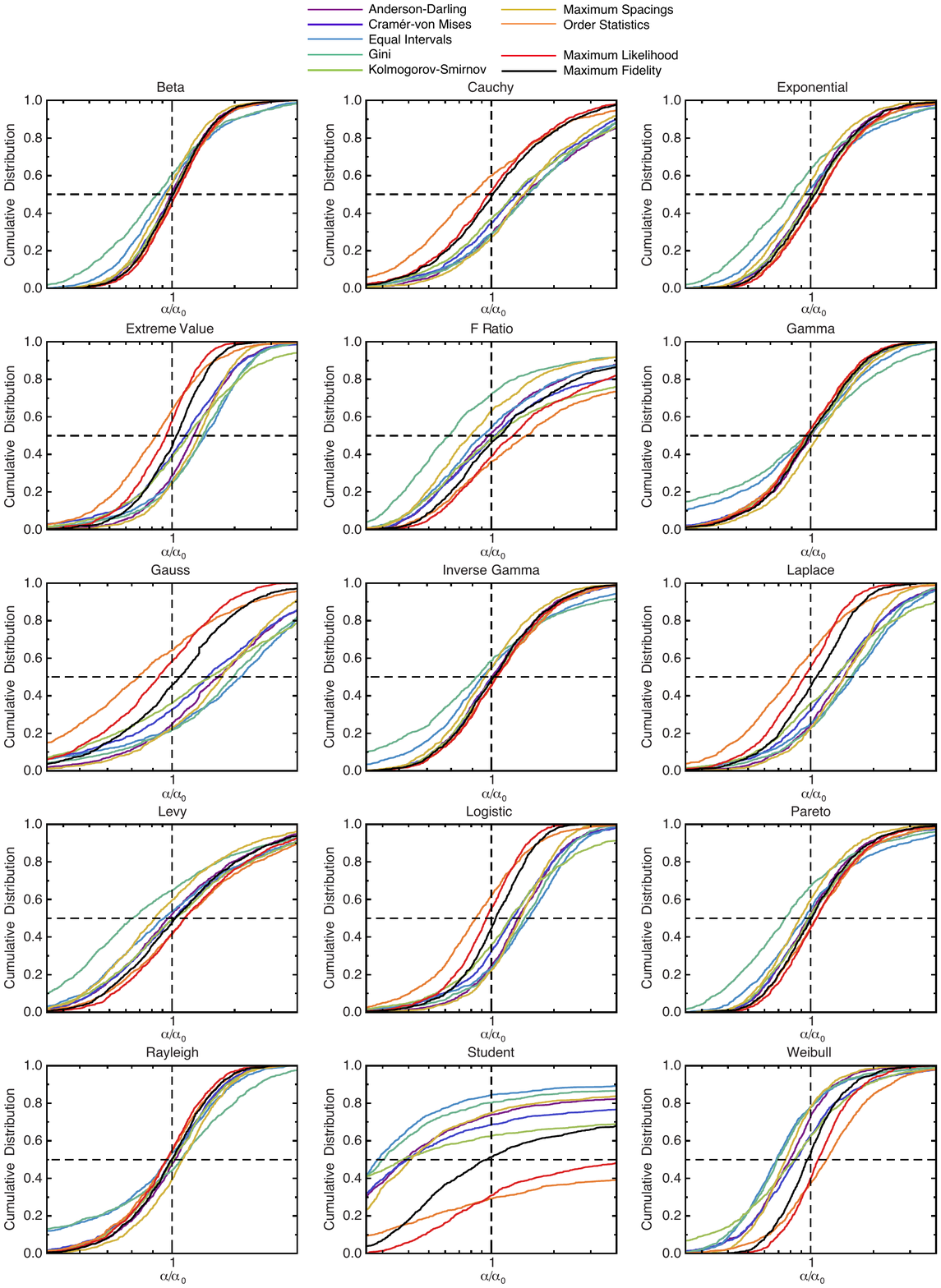}
\vspace{-2.5cm}
      \caption{Fitting of shape parameters on the line. $n=5$ data points were drawn repeatedly from the indicated distributions (values assumed for $\beta$ and $\alpha_0$ are listed in Table~\ref{tab:line_dist}). Fitting of the shape parameter $\alpha$, with $\beta$ fixed at its true value, was carried out for each realization using each indicated statistic (1000 total).}
      \label{fig:par_line}
   \end{center}
\end{figure}

In Fig.~\ref{fig:par_line}, the parameter that more generally dictates the shape or width of the distribution (rather than its location) is fit using the statistics listed in Table~\ref{tab:line_stat} for each indicated distribution family ($\beta$ and $\alpha_0$ values listed in Table~\ref{tab:line_dist}). For these statistics, the estimation of shape parameters appears to be a more stringent test than the estimation of ``location'' shown in Fig.~\ref{fig:par_line_loc}. Maximum fidelity provides the best estimate for all of the tested distributions, with maximum likelihood giving good but typically more median biased estimates. For the Student $t$ distribution, the estimate provided by maximum fidelity is far superior to that provided by maximum likelihood or any of the other statistics. In all cases, therefore, the fidelity is the clearly superior choice for parameter estimation of distribution width or shape.

\subsection{Simultaneous Fitting of Gaussian Mean and Width}

\begin{figure}[]
   \begin{center}
   \vspace{-0.2cm}
     \includegraphics*[width=4in]{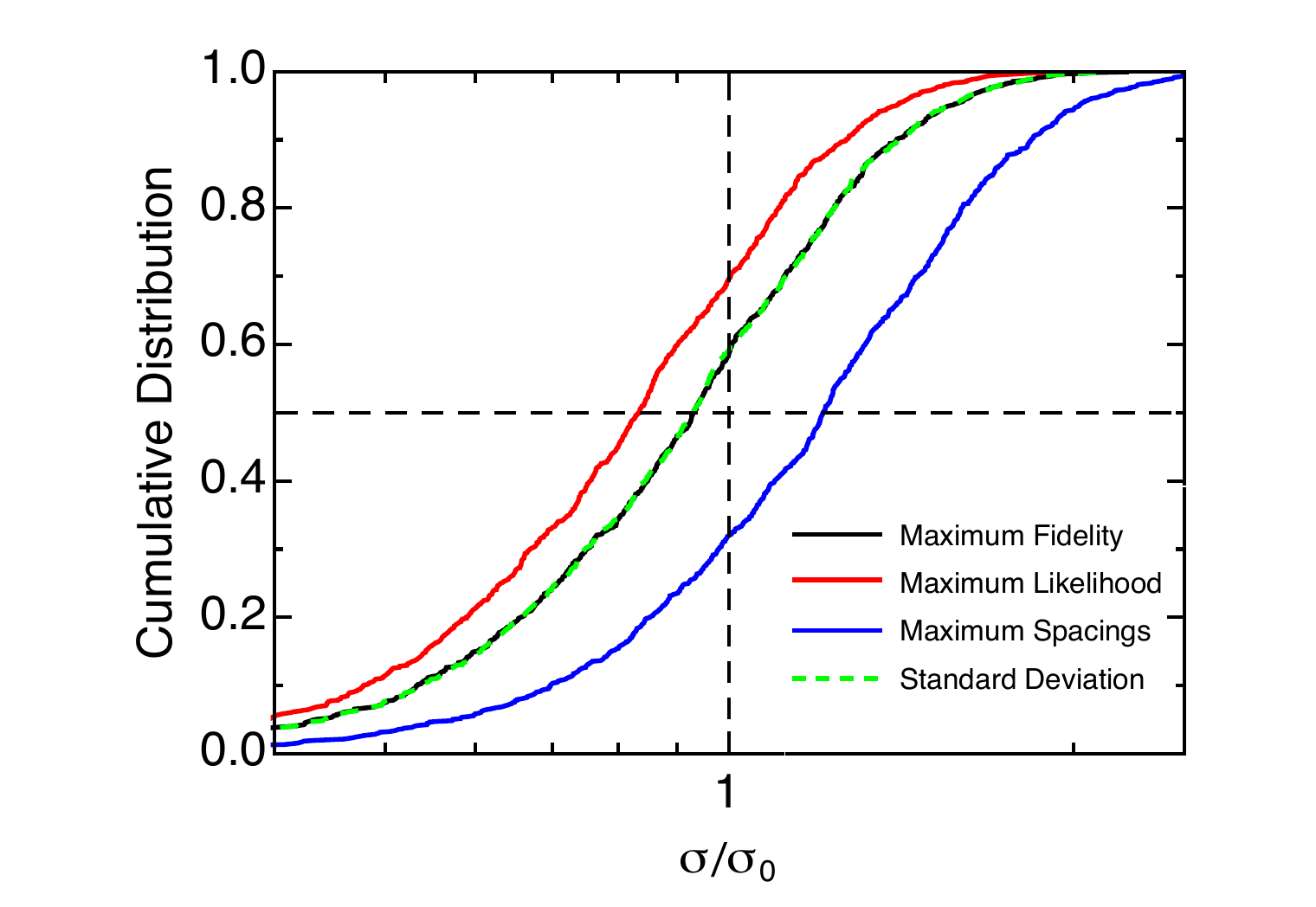}
     \vspace{-0cm}
      \caption{Estimate distributions for $\sigma$ obtained upon simultaneous fitting of Gaussian $\mu$ and $\sigma$ for $n=5$ data points (for 1000 Monte Carlo realizations). Estimate distributions obtained using maximum fidelity, maximum likelihood, maximum spacings, and the standard deviation, i.e. $\sigma_{\mathrm{SD}}=\sqrt{\sum_i (x_i-\langle x\rangle)^2/(n-1)}$.}
      \label{fig:gauss_mean_sigma}
   \end{center}
\end{figure}

In the previous sections, only a single model parameter was fit. However, the most important example of parameter fitting in statistics is the simultaneous fitting of the central value $\mu$ and width $\sigma$ of a hypothesized Gaussian. In Fig.~\ref{fig:gauss_mean_sigma}, a side-by-side comparison of maximum fidelity with maximum likelihood and maximum spacings is displayed for data sets containing $n=5$ points (1000 realizations) drawn from a Gaussian ($\mu=0$, $\sigma=1$). This trio of statistics is reminiscent of the choice of porridges presented to Goldilocks. Maximum likelihood and maximum spacings prefer Gaussian widths that are respectively narrower (``hotter'') or broader (``colder'') than the estimates obtained using the standard deviation, whereas the estimate distribution obtained from maximum fidelity perfectly overlaps that obtained from the standard deviation (``just right''). The standard deviation is generally considered to provide an ``optimal'' estimate (as the mean bias of the variance estimate $v=\sigma_{\mathrm{SD}}^2$ is zero). The bias of maximum likelihood towards distributions with higher peaks is well known (see also the below discussion of Fig.~\ref{fig:fidelityvslikelihood}), as is the need for a mean bias correction  of the likelihood estimates to match the true $\sigma^2$: $\sigma^2=\langle\sigma_{\mathrm{SD}}^2\rangle=(n/(n-1))\langle\sigma_{\mathrm{ML}}^2\rangle$. The $n/(n-1)$ correction of the maximum likelihood estimates not only removes the mean bias of $\sigma_{\mathrm{ML}}^2$, but also generates an estimate distribution that perfectly overlaps the standard deviation and maximum fidelity estimate distributions (not shown). That maximum spacings assigns too much weight to the boundary intervals is clear from Fig.~\ref{fig:symmetry}, which accounts for its preference for flatter distributions. Testing of different $n$ (not shown), leads to similar results, proving their generality. It is important to note that for any given data set, the separate $\sigma$ estimates obtained from maximum fidelity and from the standard deviation will generally differ, it is only the full parameter estimate distributions that are identical. While maximum fidelity coincides with the standard deviation in its power of estimation, as we will see in the next section, it also automatically provides a coordinate-independent measure of the degree of concordance of the hypothesized model with the data. Such valuable information is unobtainable from the standard deviation estimate or from the maximum likelihood estimate. Furthermore, parameter estimation by maximum fidelity is possible for any distribution family, unlike the standard deviation, and is also always reliable, unlike maximum likelihood (see \S\ref{sec:intro}). All of these considerations point to the fidelity as the most \textit{fundamental} statistic for parameter estimation.

\subsection{Comparison of Maximum Fidelity and Maximum Likelihood}\label{sec:fidlike}

\begin{figure}[]
   \begin{center}
\vspace{-2.3cm}
     \includegraphics*[width=6.5in]{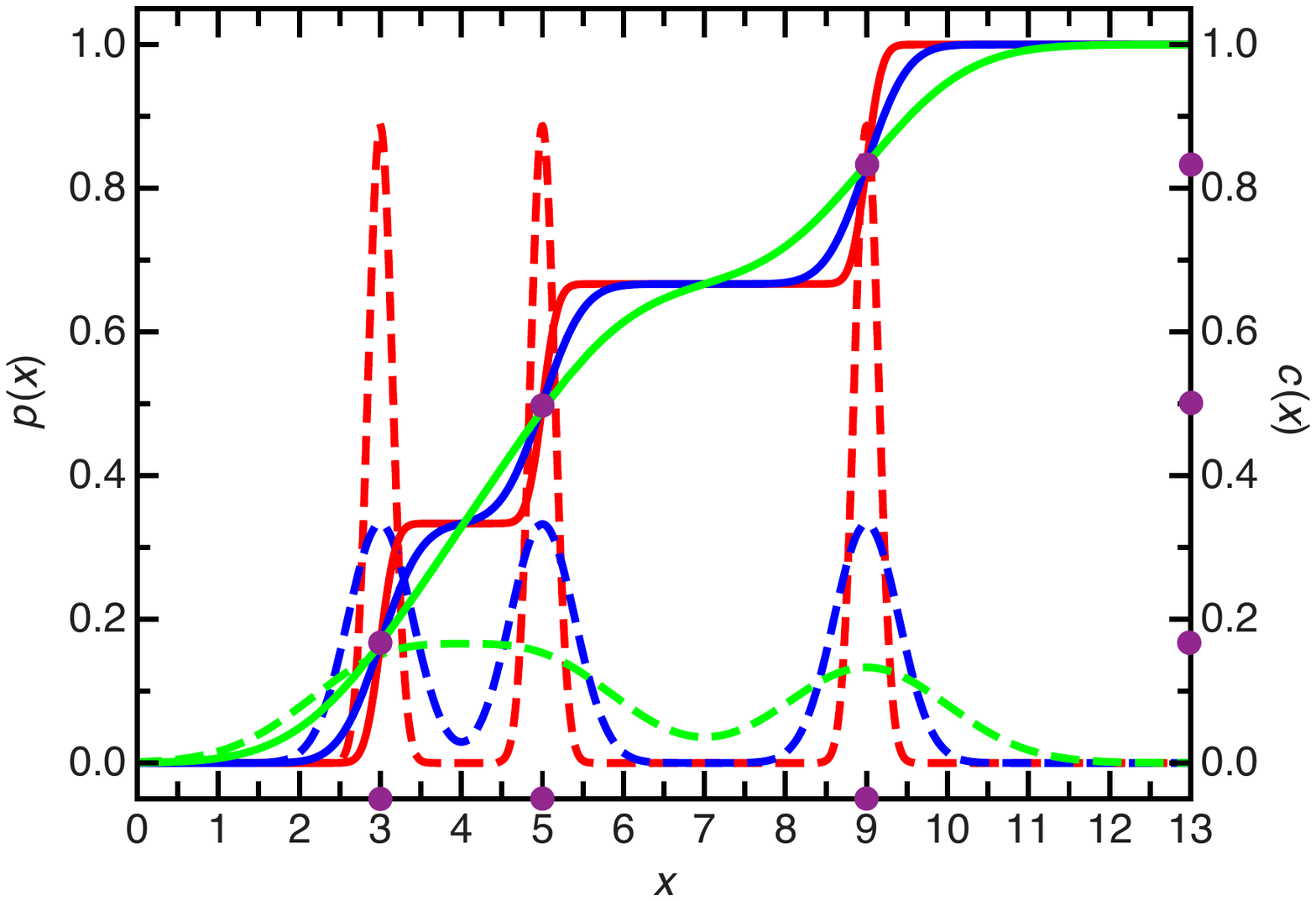}
\vspace{-3cm}
           \caption{Comparison of maximum fidelity and maximum likelihood. Three different normalized probability distributions are shown with increasingly higher peaks at the observed data values of $x=3$, 5, and 9. Maximum likelihood gives the highest likelihood to the red model, whereas maximum fidelity gives the same fidelity for all three models, as all models (by design) place the data points at the respective cumulative values of $c(x)=0.1667$, 0.5, and 0.8333.}
      \label{fig:fidelityvslikelihood}
   \end{center}
\end{figure}

A comparison of the fidelity and the likelihood is given in Fig.~\ref{fig:fidelityvslikelihood}. Model fitting in statistics is often restricted to a family of distributions (e.g. Gaussians) in order to enable the calculation of probable errors, confidence intervals, credible intervals, etc. Under these restrictions, the likelihood has traditionally played an important role. However, if we had complete freedom in the specification of the model, maximum likelihood would always lead to $\delta$-function spikes at each observed data point. By contrast, maximum fidelity is only concerned with the cumulative mapping of the data points.  In Fig.~\ref{fig:fidelityvslikelihood}, three data points are fit with increasingly narrower distributions, for which the likelihood takes on increasingly higher values. It is clear that the model with highest likelihood would indeed correspond to a sum of $\delta$-function at each observed data point. The fidelity, on the other hand, is identical for all of the models (including the $\delta$ function model), as they all map the data points to the same cumulative values.

It is clear that maximum spacings and maximum fidelity are asymptotically equivalent. But what is the asymptotic connection of these statistics with the likelihood? Ranneby has already shown that the likelihood and the spacings statistic represent different ways of approximating the \textit{continuous} version of the fidelity (or Kullback-Leibler divergence)\cite{ranneby_maximum_1984}. However, as argued above in \S\ref{sec:fidelity}, it is better to start with the more intuitive \textit{discrete} version of the fidelity and derive the likelihood as an approximation. Starting from Eq.~\ref{eq:fidline}, we can express $f_n^l$ in terms of the probability distribution as follows:
\begin{eqnarray}
f_n^l &=&   \frac{1}{2n}\log{\left[2n (1-c_{n})\right]} + \frac{1}{2n}\log{\left[2n (c_{1}-0)\right]} + \sum_{i=1}^{n-1}  \frac{1}{n}\log{\left[n(c_{i+1}-c_{i})\right]},\nonumber\\
     &=&   \frac{1}{2n}\log{\left[2n \int_{x_n}^{\infty} p(x)dx\right]} + \frac{1}{2n}\log{\left[2n \int_{-\infty}^{x_1} p(x)dx\right]} \nonumber\\
     &&\hspace{2.2in}+ \sum_{i=1}^{n-1}  \frac{1}{n}\log{\left[n\int_{x_{i}}^{x_{i+1}}p(x)dx\right]}.
\end{eqnarray}
For large $n$, we can approximate these integrals using forward finite differences (final value of $p(x)$ on each interval times the $x$ interval width):
\begin{eqnarray}
\lim_{n\rightarrow\infty} f_n^l &\approx&   \frac{1}{2n}\log{\left[2n(\infty-x_{n})p(\infty) \right]} + \frac{1}{2n}\log{\left[2n (x_1-(-\infty)) p(x_1) \right]} \nonumber\\
&&\hspace{2.2in}+ \sum_{i=1}^{n-1}  \frac{1}{n}\log{\left[n (x_{i+1}-x_i)p(x_{i+1})\right]},
\end{eqnarray}
or backward finite differences:
\begin{eqnarray}
\lim_{n\rightarrow\infty} f_n^l  &\approx&   \frac{1}{2n}\log{\left[2n(\infty-x_{n})p(x_n) \right]} + \frac{1}{2n}\log{\left[2n (x_1-(-\infty)) p(-\infty) \right]} \nonumber\\
&&\hspace{2.2in} + \sum_{i=1}^{n-1}  \frac{1}{n}\log{\left[n (x_{i+1}-x_i)p(x_{i})\right]}.\nonumber\\
\end{eqnarray}
Taking the average of these two approximations yields a better, more symmetric approximation of:
\begin{eqnarray}
\lim_{n\rightarrow\infty} f_n^l &\approx& \frac{1}{n}\log{\left[n^{n-1}(x_2-x_1)\ldots(x_n-x_{n-1})\right]}\nonumber\\
&&+\frac{1}{4n}\log{\left[(2n)^4p(-\infty)p(x_1)(x_1-(-\infty))^2p(x_n)p(\infty)(\infty-x_n)^2\right]}\nonumber\\
&&\hspace{0in}+\frac{1}{n}\log{\left[p(x_1)^{1/2} p(x_2)\ldots p(x_{n-1})p(x_n)^{1/2} \right]}.
\end{eqnarray}
The fidelity on the line is therefore asymptotically equivalent to the sum of a model-independent normalization term (interior intervals), a model-dependent term boundary term, and $1/n$ times the logarithm of the likelihood (aside from the $1/2$ exponents on the first and last terms). One can likewise show (with the result being equivalent for forward or backward finite differences) that the circular fidelity, $f_n^c$, is asymptotically:
\begin{eqnarray}
\lim_{n\rightarrow\infty} f_n^c &\approx& \frac{1}{n}\log{\left[n^n (x_2-x_1)\ldots (1+x_1-x_n)\right]}\nonumber\\
&&\hspace{1.8in}+ \frac{1}{n}\log{\left[p(x_1) p(x_2)\ldots p(x_{n-1})p(x_n)\right]},
\end{eqnarray}
which is equivalent to the sum of a model-independent normalization term and $1/n$ times the logarithm of the likelihood function.

%%%%%%%%%%%%%%%%%%%%%%%%%%%%%%%%%%%%%%%%%%%%%%%%%%%%%%%%%%%%%%
\section{Conversion of Fidelity to Concordance}\label{sec:concordance}
%%%%%%%%%%%%%%%%%%%%%%%%%%%%%%%%%%%%%%%%%%%%%%%%%%%%%%%%%%%%%%

In this section, I show how the fidelity can be efficiently converted to an absolute measure of model concordance ($p$ value) under assumption of the null hypothesis (i.e. the hypothesis that the trial distribution used to calculate the fidelity is the true distribution from which the data were drawn).

Arguably the most important and widely used tool in statistics is Karl Pearson's $\chi^2$-$p$ test, which was the first quantitative and universal approach for measuring model goodness-of-fit\cite{pearson_criterion_1900}. Unfortunately, this test requires binned data with a sufficiently high level ($>10$) of counts in each bin. It is therefore only reliably defined asymptotically. Since this time, several other goodness-of-fit measures have been proposed based on various statistics of the cumulative distribution of the data points\cite{dagostino_goodness--fit-techniques_1986}. While these alternatives are valid for any $n$, none of these can compete with the generality and power of Pearson's $\chi^2$ (in its range of validity). In this section, I will derive a highly accurate gamma-based approximation for the $p$ value of the fidelity, which is valid for any number of data points $n$. As the fidelity is not limited to asymptotic situations (high $n$), it provides an even more fundamental basis for goodness-of-fit than $\chi^2$. As maximum fidelity is asymptotically consistent with maximum likelihood (\S\ref{sec:fidlike}) and it is well-known that maximum likelihood is asymptotically equivalent to minimum $\chi^2$, maximum fidelity is also asymptotically equivalent to minimum $\chi^2$. As all three estimation approaches are \textit{consistent} estimators in the asymptotic limit (see Ranneby (1984) for a proof of the consistency of the spacings statistic\cite{ranneby_maximum_1984} that is also valid for the fidelity), minimum $\chi^2$ can be thought of as simply the asymptotic version of maximum fidelity. This glosses over the ``binned'' nature of the data required for $\chi^2$ analysis. Often, data are artificially binned (resulting in degraded resolution) in order to simply apply the power of $\chi^2$ analysis (with the number of data points in each bin ensured to be greater than 10). Use of the fidelity, however, obviates the need for such artificial binning. For data that is \textit{collected} in a binned fashion, however, maximum fidelity can also be applied in a consistent fashion (see \S\ref{sec:binned}), which further bolsters the above statement of minimum $\chi^2$ as the asymptotic version of maximum fidelity.

Calculation of the expected distribution of a cumulative-based statistic under the null hypothesis simply requires repeatedly drawing $n$ data points from a uniform distribution on the unit interval and calculating the resulting statistic. The statistics of random points distributed on the line, as well as the intervals that they create, was first  explored by Whitworth\cite{whitworth_choice_1901,irwin_william_1967}. In our particular case, we are interested in the statistics of the product of intervals on the line created by the $n$ randomly drawn data points. Knowledge of this null distriution of the fidelity, will allow us to convert the fidelity to a concordance measure ($p$ value). The corresponding null distribution for the fidelity or spacings statistics cannot be expressed explicitly.  However, through use of the Laplace transform, Darling\cite{darling_class_1953} was able to calculate the exact solution for the characteristic function corresponding to the distribution of the spacings statistic (sum of the log of the intervals for $n$ points on the line):
\begin{eqnarray}
\phi_n^s(\xi)=\left\langle\mathrm{e}^{i\xi s_n^l}\right\rangle&=&\frac{n!}{2\pi i}\int_{c-i\infty}^{c+i\infty}e^z\left(\int_0^{\infty}e^{-rz+i\xi\log r}dr\right)^{n+1}dz\nonumber\\
&=&\Gamma(n+1) \frac{\left((n+1)^{\frac{i\xi}{n+1}}\Gamma(1+\frac{i\xi}{n+1})\right)^{n+1}}{\Gamma\left((n+1)(1+\frac{i\xi}{n+1})\right)}.
\end{eqnarray}
By the same method, the corresponding result for the distribution of the fidelity for $n$ points on the circle is:
\begin{equation}
\phi_n^c(\xi)=\left\langle\mathrm{e}^{i\xi f_n^c}\right\rangle=\Gamma(n) \frac{\left(n^{\frac{i\xi}{n}}\Gamma(1+\frac{i\xi}{n})\right)^{n}}{\Gamma\left(n\left(1+\frac{i\xi}{n}\right)\right)},
\end{equation}
and the distribution of the fidelity for $n$ points on the line is:
\begin{equation}
\phi_n^l(\xi)=\left\langle\mathrm{e}^{i\xi f_n^l}\right\rangle=\Gamma(n+1) \frac{\left(n^{\frac{i\xi}{n}}\Gamma(1+\frac{i\xi}{n})\right)^{n-1}\left((2n)^{\frac{i\xi}{2n}}\Gamma(1+\frac{i\xi}{2n})\right)^{2}}{\Gamma\left((n-1)(1+\frac{i\xi}{n})+2\left(1+\frac{i\xi}{2n}\right)\right)}.
\end{equation}
The first two cumulants (corresponding to the mean and the variance) of the null distribution for each statistic follow from differentiation of the log of the above characteristic functions\cite{james_statistical_2006}.  The result for the spacings distribution on the line (already obtained by Darling\cite{darling_class_1953}) is:
\begin{eqnarray}
\mu_n^s&=&-\gamma+\big(\log{(n+1)}-\psi^{(0)}(n+1)\big)\\
(\sigma_n^s)^2&=&\frac{1}{n+1}\bigg[\frac{\pi^2}{6}-1-\big((n+1)\,\psi^{(1)}(n+1)-1\big)\bigg].
\end{eqnarray}
The corresponding cumulants for the fidelity distribution on the circle are:
\begin{eqnarray}
\mu_n^c&=&-\gamma+\big(\log{n}-\psi^{(0)}(n)\big)\\
(\sigma_n^c)^2&=&\frac{1}{n}\bigg[\frac{\pi^2}{6}-1-\big(n\,\psi^{(1)}(n)-1\big)\bigg],
\end{eqnarray}
and for the fidelity distribution on the line are:
\begin{eqnarray}
\mu_n^l&=&-\gamma+\bigg(\log{n}-\psi^{(0)}(n+1)+\frac{\log{2} }{n}\bigg)\\
(\sigma_n^l)^2&=&\frac{1}{n}\bigg[\frac{\pi^2}{6}-1-\bigg(n\,\psi^{(1)}(n+1)-1+\frac{\pi^2}{12n}\bigg)\bigg].
\end{eqnarray}
In the above expressions, $\gamma= 0.5772156649\ldots$ is the Euler constant and
\begin{equation}
\psi^{(m)}(z)=\bigg(\frac{d}{dz}\bigg)^{m+1} \log{\Gamma(z)}
\end{equation}
are the polygamma functions.
Also, in the above, all expressions in parentheses approach zero asymptotically ($n\rightarrow\infty$). 

\begin{figure}[]
   \begin{center}
\vspace{-6.5cm}
     \includegraphics*[width=6in]{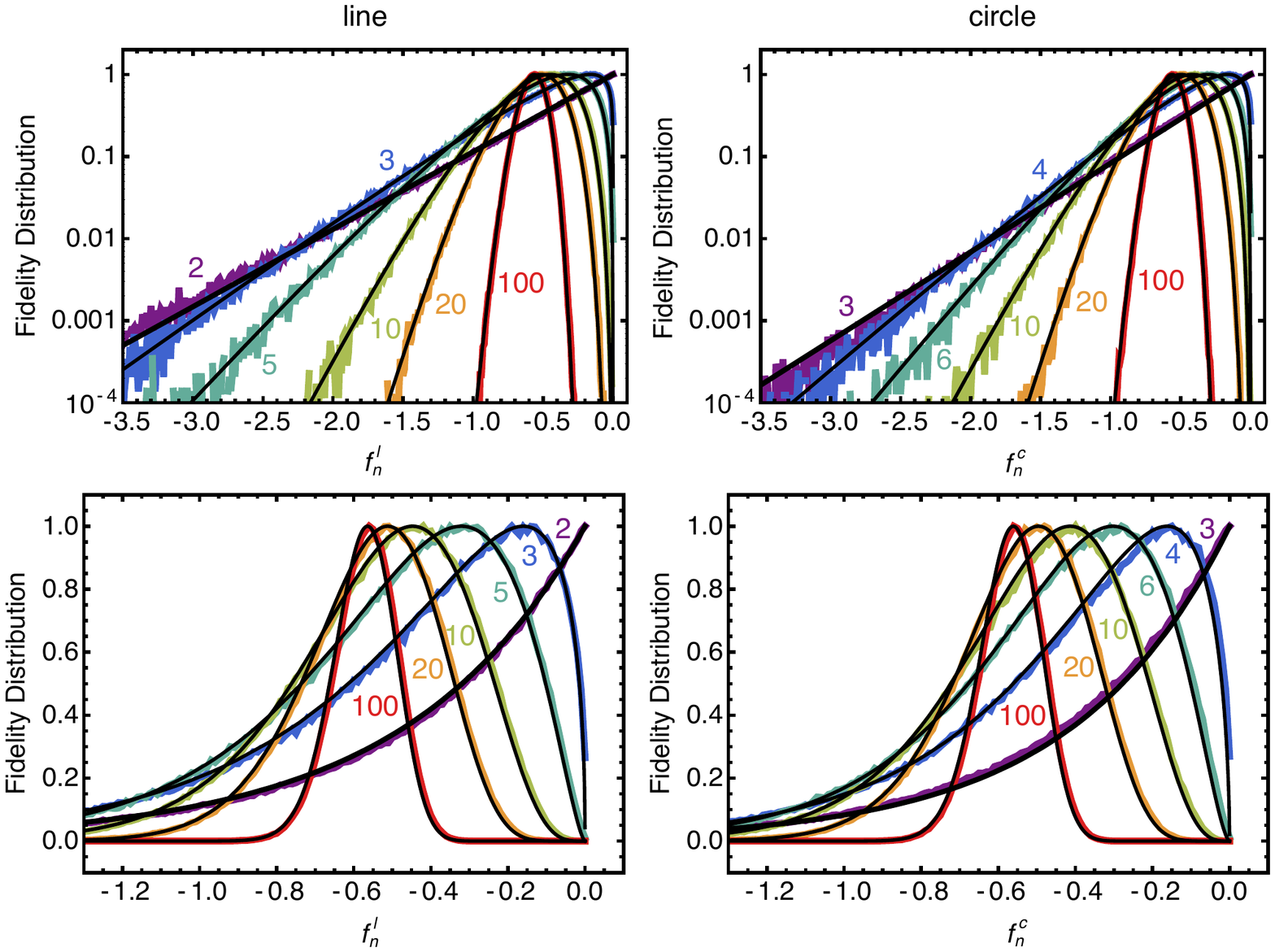}
\vspace{-4.cm}
      \caption{Peak-normalized probability distributions of the fidelity on the line and circle under the null hypothesis (1 million realizations at each $n$). Overlaid in black are the approximations for the fidelity summarized in Table~\ref{tab:approx}.}
      \label{fig:fid_dist}
   \end{center}
\end{figure}

A simple gamma-function approximation to the null distribution of the fidelity based on the first two cumulants can be defined as follows:
\begin{equation}
P^g_n(f^g_n)\simeq\frac{\left(\beta^g_n\right)^{\alpha^g_n}}{\Gamma(\alpha^g_n)}  \left(-f^g_n\right)^{\alpha^g_n-1}  \mathrm{e}^{-\beta^g_n f^g_n},
\end{equation}
where the geometry $g$ is either $c$ (fidelity on the circle) or $l$ (fidelity on the line).
The corresponding $p$ value is then given by:
\begin{equation}
p^g_n(f^g_n) = \int_{-\infty}^{f^g_n}P^g_n(x)dx\simeq\frac{\Gamma(\alpha^g_n,-\beta^g_n f^g_n)}{\Gamma(\alpha^g_n)}=Q(\alpha^g_n,-\beta^g_n f^g_n),
%=\frac{\Gamma(\mu_n^2/\sigma_n^2,\mu_n f_n/\sigma_n^2)}{\Gamma(\mu_n^2/\sigma_n^2)}.
\label{eq:concordance}
\end{equation}
where $Q(s,x)$ is the regularized gamma function.
Setting the first two cumulants of the gamma distribution (which completely specify it) to the exact solutions given above yields:
\begin{eqnarray}
\alpha^g_n&=&\frac{(\mu^g_n)^2}{(\sigma^g_n)^2}\\
\beta^g_n&=&-\frac{\mu^g_n}{(\sigma^g_n)^2}.
\end{eqnarray}
With this definition, remarkably good agreement with Monte Carlo calculations is obtained for all $n\geq4$ on the circle and $n\geq3$ on the line (see Fig.~\ref{fig:fid_dist}).  For $n=1$ on the circle, the distribution is trivial.  For $n=2$ on the circle and $n=1$ on the line, the distributions can be solved for exactly:
\begin{eqnarray}
P^g_n(f^g_n)&=&\frac{\mathrm{e}^{2f_n^g}}{\sqrt{1-\mathrm{e}^{2f_n^g}}}\\
p^g_n(f^g_n)&=&1-\sqrt{1-\mathrm{e}^{2f_n^g}},
\end{eqnarray}
For $n=3$ on the circle and $n=2$ on the line, the distributions are approximated very well by a simple exponential (Fig.~\ref{fig:fid_dist}):
\begin{eqnarray}
P_n^g(f_n^g)&\simeq&-\frac{1}{\mu_n^g}\mathrm{e}^{-f_n^g/\mu_n^g}\\
p_n^g(f_n^g)&\simeq&\mathrm{e}^{-f_n^g/\mu_n^g}.
\end{eqnarray}
The complete system of approximations for the $p$ value --- summarized in Table~\ref{tab:approx} with computed coefficients listed in Tables~\ref{tab:circ_alpha_beta} and \ref{tab:line_alpha_beta} for a range of $n$ --- serves as a solid and extremely efficient foundation for assessment of model concordance using maximum fidelity.  For the spacings statistic, Cheng \& Stephens\cite{cheng_goodness--fit_1989} presented a similar, though slightly less accurate approach obtained by fitting the null distribution to a $\chi^2$ distribution with $n$ degrees of freedom plus a displacement term correction.

\begin{table}[]
\begin{center}
\begin{tabular}{|c|c||c|c|}
%\hline
%\multicolumn{2}{c}{line} & \multicolumn{2}{c}{circle}  \\
\hline
$n$  &  $p^l_n(f^l_n)$    &   $p^c_n(f^c_n)$ & n \\ 
\hline
1 \T&  $1-\sqrt{1-\mathrm{exp}(2f_1^l)} $   &   1 & 1\\
2 \T&  $\mathrm{exp}\left(-f_2^l/\mu_2^l\right)$ &  $1-\sqrt{1-\mathrm{exp}\left(2f_2^c\right)} $  & 2   \\
$\geq3$ \T& $Q\left(\alpha^l_n,-\beta^l_n f^l_n\right)$   &  $\mathrm{exp}\left(-f_3^c/\mu_3^c\right)$  & 3\\
   \T\B&    &   $Q\left(\alpha^c_n,-\beta^c_n f^c_n\right)$  & $\geq4$\\
\hline
\end{tabular}
\end{center}
\caption{Complete set of explicit solutions and power-law or gamma-function approximations for the $p$ value of the fidelity for each $n$ on the circle and the line.}
\label{tab:approx}
\end{table}

\begin{table}[]
\begin{center}
\begin{tabular}{|c|c|c|c|c|}
\hline
$n$  & $\mu^c_n$ & $\sigma^c_n$ & $\alpha^c_n$    & $\beta^c_n$   \\ 
\hline
1 & 0 & 0 & --- & --- \\
2 & -0.30685281944 & 0.42134661097 & 0.5303727787 & 1.7284272626     \\
3 &  -0.40138771133 & 0.39163412615 & 1.0504299307 & 2.6169957401   \\
4 &   -0.44703897221 & 0.35694615977 & 1.5685029651 & 3.5086492736    \\
5 &   -0.47389542090 & 0.32812171161 & 2.0859077028 & 4.4016202960   \\
6 &   -0.49157386411 & 0.30468462723 & 2.6030138779 & 5.2952650007  \\
7 &   -0.50408985094 & 0.28538641001 & 3.1199621882 & 6.1892977659   \\
8 &   -0.51341560118 & 0.26922062265 & 3.6368172156 & 7.0835736336 \\
9 &   -0.52063256552 & 0.25545730989 & 4.1536126209 & 7.9780115497   \\
10 &   -0.52638316097 & 0.24357149054 & 4.6703676327 & 8.8725627622   \\
11 &   -0.53107298117 & 0.23318045574 & 5.1870940218 & 9.7671962343    \\
12 &   -0.53497069509 & 0.22399992418 & 5.7037993931 & 10.661891288    \\
13 &   -0.53826132075 & 0.21581464577 & 6.2204888775 & 11.556633623  \\
14 &   -0.54107642552 & 0.20845868134 & 6.7371660617 & 12.451413043   \\
15 &   -0.54351212546 & 0.20180199027 & 7.2538335284 & 13.346222078   \\
16 &   -0.54564027099 & 0.19574114482 & 7.7704931842 & 14.241055137   \\
17 &   -0.54751564917 & 0.19019278130 & 8.2871464679 & 15.135907952   \\
18 &   -0.54918076474 & 0.18508889795 & 8.8037944850 & 16.030777205   \\
19 &   -0.55066909903 & 0.18037342145 & 9.3204380992 & 16.925660284    \\
20 &   -0.55200738359 & 0.17599966025 & 9.8370779953 & 17.820555101   \\
21 &   -0.55321721942 & 0.17192838836 & 10.353714723 & 18.715459966    \\
22 &   -0.55431625140 & 0.16812638422 & 10.870348728 & 19.610373502   \\
23 &   -0.55531903429 & 0.16456530295 & 11.386980378 & 20.505294568    \\
24 &   -0.55623768074 & 0.16122079606 & 11.903609974 & 21.400222219 \\
25 &   -0.55708235289 & 0.15807181728 & 12.420237772 & 22.295155658  \\
26 &  -0.55786163973 & 0.15510006970 & 12.936863984 & 23.190094215    \\
27 &   -0.55858285021 & 0.15228956183 & 13.453488791 & 24.085037315    \\
28 &   -0.55925224308 & 0.14962624783 & 13.970112348 & 24.979984471   \\
29 &   -0.55987520898 & 0.14709773409 & 14.486734786 & 25.874935260    \\
30 &   -0.56045641592 & 0.14469303786 & 15.003356221 & 26.769889316   \\
31 &   -0.56099992644 & 0.14240238764 & 15.519976752 & 27.664846322   \\
32 &   -0.56150929264 & 0.14021705698 & 16.036596465 & 28.559806000   \\
33 &   -0.56198763397 & 0.13812922520 & 16.553215436 & 29.454768105    \\
34 &   -0.56243770112 & 0.13613186018 & 17.069833733 & 30.349732422  \\
35 &   -0.56286192896 & 0.13421861914 & 17.586451414 & 31.244698760    \\
36 &   -0.56326248056 & 0.13238376421 & 18.103068531 & 32.139666950   \\
37 &   -0.56364128415 & 0.13062209032 & 18.619685131 & 33.034636842    \\
38 &   -0.56400006410 & 0.12892886329 & 19.136301256 & 33.929608300    \\
39 &   -0.56434036717 & 0.12729976649 & 19.652916943 & 34.824581204   \\
40 &   -0.56466358482 & 0.12573085473 & 20.169532225 & 35.719555443  \\
50 &   -0.56718233290 & 0.11268251022 & 25.335668201 & 44.669353630   \\
100 &   -0.57220733165 & 0.07999483737 & 51.166170842 & 89.418936130    \\
200 &   -0.57471358157 & 0.05667582819 & 102.82700128 & 178.91869024    \\
1000 &   -0.57671558157 & 0.02538570267 & 516.11323043 & 894.91813110   \\
\hline
\end{tabular}
\end{center}
\caption{Gamma-function approximation coefficients for determining the $p$ value from the fidelity on the circle.}
\label{tab:circ_alpha_beta}
\end{table}

\begin{table}[]
\begin{center}
\begin{tabular}{|c|c|c|c|c|}
\hline
$n$  & $\mu^l_n$ & $\sigma^l_n$ & $\alpha^l_n$    & $\beta^l_n$   \\ 
\hline
1 & -0.30685281944 & 0.42134661097 & 0.5303727787 & 1.7284272626 \\
2 &  -0.46027922916 & 0.47107983211 & 0.9546710017 & 2.0741127151   \\
3 &   -0.50367198448 & 0.41605669559 & 1.4655159820 & 2.9096634857 \\
4 &  -0.52375217707 & 0.37216444132 & 1.9805323051 & 3.7814302103   \\
5 &  -0.53526598479 & 0.33877009357 & 2.4964861617 & 4.6640104783  \\
6 &  -0.54271600068 & 0.31267266168 & 3.0127662005 & 5.5512757995  \\
7 &   -0.54792596801 & 0.29166509010 & 3.5291885826 & 6.4409952963   \\
8 &  -0.55177220361 & 0.27432407161 & 4.0456825350 & 7.3321608239 \\
9 &   -0.55472732324 & 0.25971176754 & 4.5622161702 & 8.2242499676  \\
10 & -0.55706844292 & 0.24718899787 & 5.0787734475 & 9.1169649117   \\
11 &  -0.55896869203 & 0.23630560624 & 5.5953456092 & 10.010123445   \\
12 &   -0.56054176338 & 0.22673516237 & 6.1119275604 & 10.903607830  \\
13 & -0.56186538378 & 0.21823485573 & 6.6285161824 & 11.797338604  \\
14 &  -0.56299448405 & 0.21062004241 & 7.1451094857 & 12.691260196 \\
15 &  -0.56396898009 & 0.20374758571 & 7.6617061562 & 13.585332574    \\
16 &  -0.56481857220 & 0.19750464279 & 8.1783052999 & 14.479526174   \\
17 &   -0.56556581502 & 0.19180092539 & 8.6949062935 & 15.373818683   \\
18 &   -0.56622814360 & 0.18656323700 & 9.2115086925 & 16.268192948  \\
19 &   -0.56681924742 & 0.18173153891 & 9.7281121740 & 17.162635564    \\
20 &  -0.56735002456 & 0.17725606569 & 10.244716499 & 18.057135905    \\
21 &   -0.56782925844 & 0.17309517558 & 10.761321490 & 18.951685440    \\
22 &  -0.56826410683 & 0.16921372498 & 11.277927009 & 19.846277238   \\
23 &  -0.56866046122 & 0.16558182269 & 11.794532952 & 20.740905613     \\
24 & -0.56902321488 & 0.16217386361 & 12.311139238 & 21.635565854   \\
25 &  -0.56935646566 & 0.15896777083 & 12.827745803 & 22.530254026 \\
26 &   -0.56966367125 & 0.15594439507 & 13.344352596 & 23.424966817   \\
27 &  -0.56994776945 & 0.15308703435 & 13.860959576 & 24.319701417     \\
28 &  -0.57021127234 & 0.15038104652 & 14.377566712 & 25.214455430    \\
29 & -0.57045634068 & 0.14781353401 & 14.894173977 & 26.109226797      \\
30 & -0.57068484324 & 0.14537308573 & 15.410781349 & 27.004013741   \\
31 &  -0.57089840448 & 0.14304956395 & 15.927388810 & 27.898814720    \\
32 &   -0.57109844324 & 0.14083392740 & 16.443996347 & 28.793628387   \\
33 &   -0.57128620426 & 0.13871808353 & 16.960603947 & 29.688453564    \\
34 &  -0.57146278405 & 0.13669476430 & 17.477211599 & 30.583289214   \\
35 & -0.57162915237 & 0.13475742130 & 17.993819297 & 31.478134420    \\
36 &  -0.57178616999 & 0.13290013672 & 18.510427032 & 32.372988371    \\
37 &   -0.57193460359 & 0.13111754743 & 19.027034798 & 33.267850343    \\
38 &   -0.57207513829 & 0.12940477988 & 19.543642592 & 34.162719691    \\
39 &  -0.57220838818 & 0.12775739416 & 20.060250408 & 35.057595838   \\
40 &  -0.57233490531 & 0.12617133563 & 20.576858242 & 35.952478263   \\
50 &  -0.57331938929 & 0.11299717384 & 25.742937131 & 44.901563792  \\
100 &   -0.57527585985 & 0.08010572577 & 51.573330873 & 89.649739321    \\
200 & -0.57624784567 & 0.05671497003 & 103.23410946 & 179.14879897      \\
1000 &   -0.57702243439 & 0.02538919914 & 516.52029838 & 895.14768854   \\
\hline
\end{tabular}
\end{center}
\caption{Gamma-function approximation coefficients for determining the $p$ value from the fidelity on the line.}
\label{tab:line_alpha_beta}
\end{table}

\subsection{Concordance Landscape}

\begin{figure}[]
   \begin{center}
     \vspace{-8cm}
     \includegraphics*[width=6.5in]{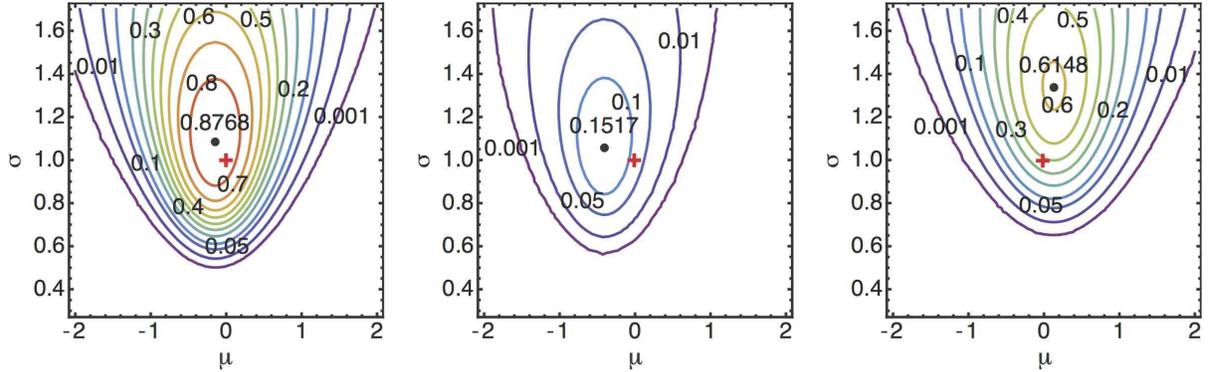}
     \vspace{-9.6cm}
      \caption{Concordance ($p$ value) landscapes for three Monte Carlo realizations of $n=20$ points on the line drawn from a Gaussian distribution with $\mu=0$ and $\sigma=1$. }
      \label{fig:gauss_land}
   \end{center}
\end{figure}

\begin{figure}[]
   \begin{center}
     \vspace{-9cm}
     \includegraphics*[width=6.5in]{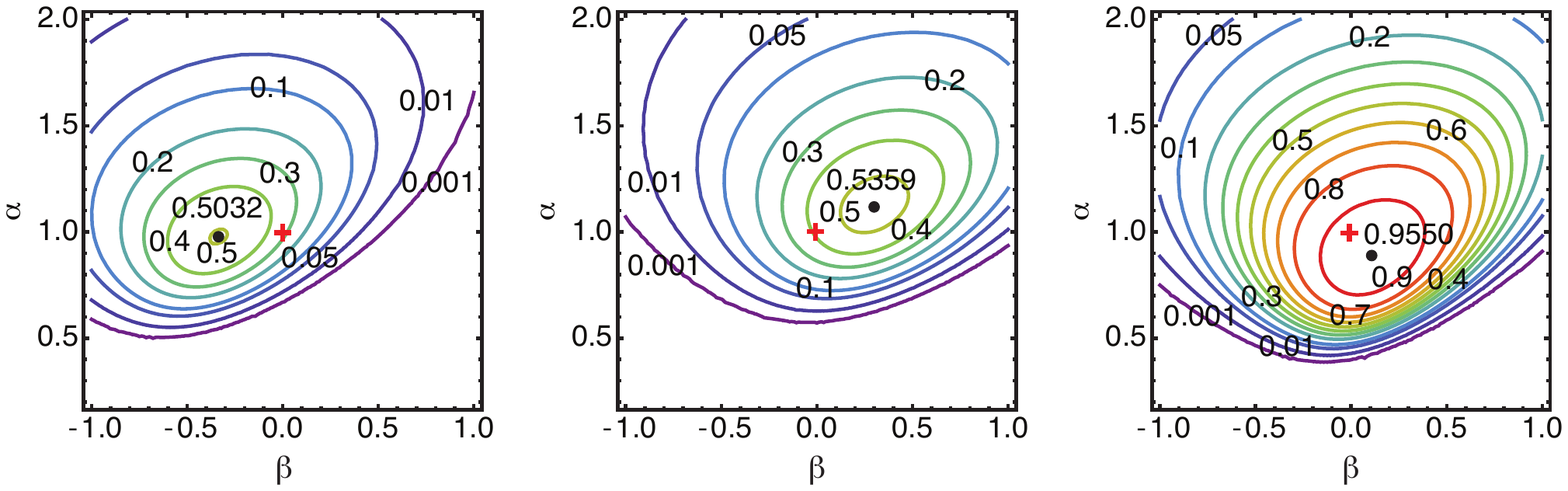}
     \vspace{-8.8cm}
      \caption{Concordance ($p$ value) landscapes for three Monte Carlo realizations of $n=40$ points on the line drawn from an Extreme Value distribution with $\beta=0$ and $\alpha=1$. }
      \label{fig:extreme_land}
   \end{center}
\end{figure}

Armed with the powerful approximations derived in the previous section, we can directly convert the fidelity to concordance ($p$ value) for a specific model fit to an observed data set (as was previewed in Fig.~\ref{fig:fidover}). We can also immediately plot entire concordance landscapes over the fitting parameters for a particular observed data set. This is demonstrated for a Gaussian fit (mean $\mu$ and standard deviation $\sigma$) to three different observed data sets (all with $n=20$) in Fig.~\ref{fig:gauss_land}. Note that the maximum $p$ value in each case depends on the data set. Concordance determined through the fidelity is, of course, completely general for any model distribution, not just Gaussian distributions. In Fig.~\ref{fig:extreme_land}, similar concordance landscapes are shown for the fitting of Extreme Value distributions to three different data sets (all with $n=40$).

In the concordance landscapes above, every point in the space represents an individual model for the data and gives an \textit{absolute} $p$-value measure of goodness-of-fit, echoing Peirce's dictum that induction ``sets out with a theory and it measures the degree of concordance of that theory with fact.'' (Peirce 5.145)\cite{peirce_collected_1974} If we decided we were interested in maximizing the fidelity over the model parameters and only considering the ``$p$ value obtained at the maximum'', then we should be wary of a few things, including the number of independent model degrees of freedom. We should first consider the \textit{absolute} $p$ value at the maximum as defined above. If this value is very low, we can of course reject the entire family of distributions. If the $p$ value is acceptable, there are two choices: We could retain the absolute $p$ value and look at how large a domain of acceptability exists (e.g., examine the zone over which the absolute $p$ value is $p>0.05$ in Figs.~\ref{fig:gauss_land} and \ref{fig:extreme_land}). Alternatively, we could attempt to ``correct'' the absolute $p$ value for the independent parameter degrees of freedom over which we have maximized our model, as Fisher proposed for Pearson's $\chi^2$ test\cite{fisher_interpretation_1922,baird_fisher/pearson_1983}. This ``correction'' provides information at the distribution family level, asking whether, for example, the data can be fit by a Gaussian (parameters undefined). It theoretically allows rejection of otherwise suitable fits determined by the absolute $p$ values. However, such a correction requires further assumptions about the shape of the fidelity or concordance landscape (sufficiently Gaussian) that are valid only asymptotically (large $n$). It also prohibits comparison of the particular maximized model with other models on the concordance landscape (these degrees of freedom are already ``accounted for'' in the maximization process). Such bias correction would also lead to a different corrected $p$ value for the exact same model arrived at in a different manner (different $p$ values would obtain for the exact same model arrived at by fitting $\mu$ and $\sigma$ simulatenously, or by fixing $\mu$ and fitting $\sigma$; these corrected $p$ values would nevertheless both be smaller than the absolute $p$ value), which prevents assigning an \textit{absolute} concordance value to each model. While such a correction allows comparisons at the distribution family level (e.g. ``the data set is better fit with a Gaussian than with a Lorentzian''), it prevents cross-comparison (e.g. of a particular Gaussian model with a particular Lorentzian).

While providing a certain type of information, ``correction'' of the $p$ value for the number of degrees of freedom can actually discourage examination of other information about the model fit. Consider a single data point $x_1$ fit by a Gaussian with fixed $\sigma$ and unknown location parameter $\mu$. Maximization of the fidelity would lead to $\mu=x_1$, with a cumulative value of $c_1=1/2$ at $x_1$, corresponding to a maximal fidelity of $f=0$ and a $p$ value of 1. By a literal application of Fisher's correction, we have zero degrees of freedom after accounting for our single data point and our one free model parameter, leading to an indefinite $p$ value (division by zero), and therefore an indefinite statement on the goodness-of-fit of this class of models. It is clear in this case, however, that the concordance landscape of absolute $p$ values contains important information both about the range of models with good concordance $p>0.05$ as well as those that provide unacceptable fits (and exactly how unacceptable). This is true for data sets containing multiple points as well. Consider $n$ data points fit with a model with $n$ linearly independent parameters. The corrected $p$ value is undefined (division by 0); however, the absolute $p$ value evaluated across the parameter landscape is still highly informative.

All of the preceding glosses over one important point about scientific inference. When we attempt to fit data with a model, often we have already looked at the data (\textit{inspection bias}) or, if we have not looked, we have a physical model in mind (\textit{theory bias}) or we are simply attempting to fit the data using a well-known distribution family (\textit{distribution bias}). All of these biases are, unfortunately, unquantifiable. Therefore, when we attempt to fit our data with a Gaussian and then notice the wings are not fit well, so we attempt a Lorentzian, this bias cannot be quantified, even if we ``correct'' our absolute $p$ value for the degrees of freedom of the model. Instead of choosing a different model with the same number of free parameters, we could have alternatively added an extra degree of freedom to the model (e.g. a skewness or a second Gaussian). Fisher's correction supposedly gives us guidance on how to determine whether the introduction of that extra degree of freedom actually generated an informative improvement in the ``corrected'' $p$ value of the maximized model fits. However, we may have introduced exactly that degree of freedom in order to explain some feature in the data. This type of bias appears unavoidable. Therefore, it is likely best to retain the absolute $p$ values and instead judiciously apply Ockham's razor to isolate the simplest models that are sufficiently concordant with the data (see \S\ref{sec:discussion}).

For those who are nevertheless interested in ``correcting'' the $p$ value, Cheng \& Stephens give some useful guidance for ``correcting'' the spacings statistic based on its asymptotic equivalence with the likelihood function\cite{cheng_goodness--fit_1989}; their analysis should be directly applicable to the fidelity.

\subsection{Goodness-of-fit on the Circle}

\begin{figure}[]
   \begin{center}
     \vspace{1cm}
     \includegraphics*[width=6.5in]{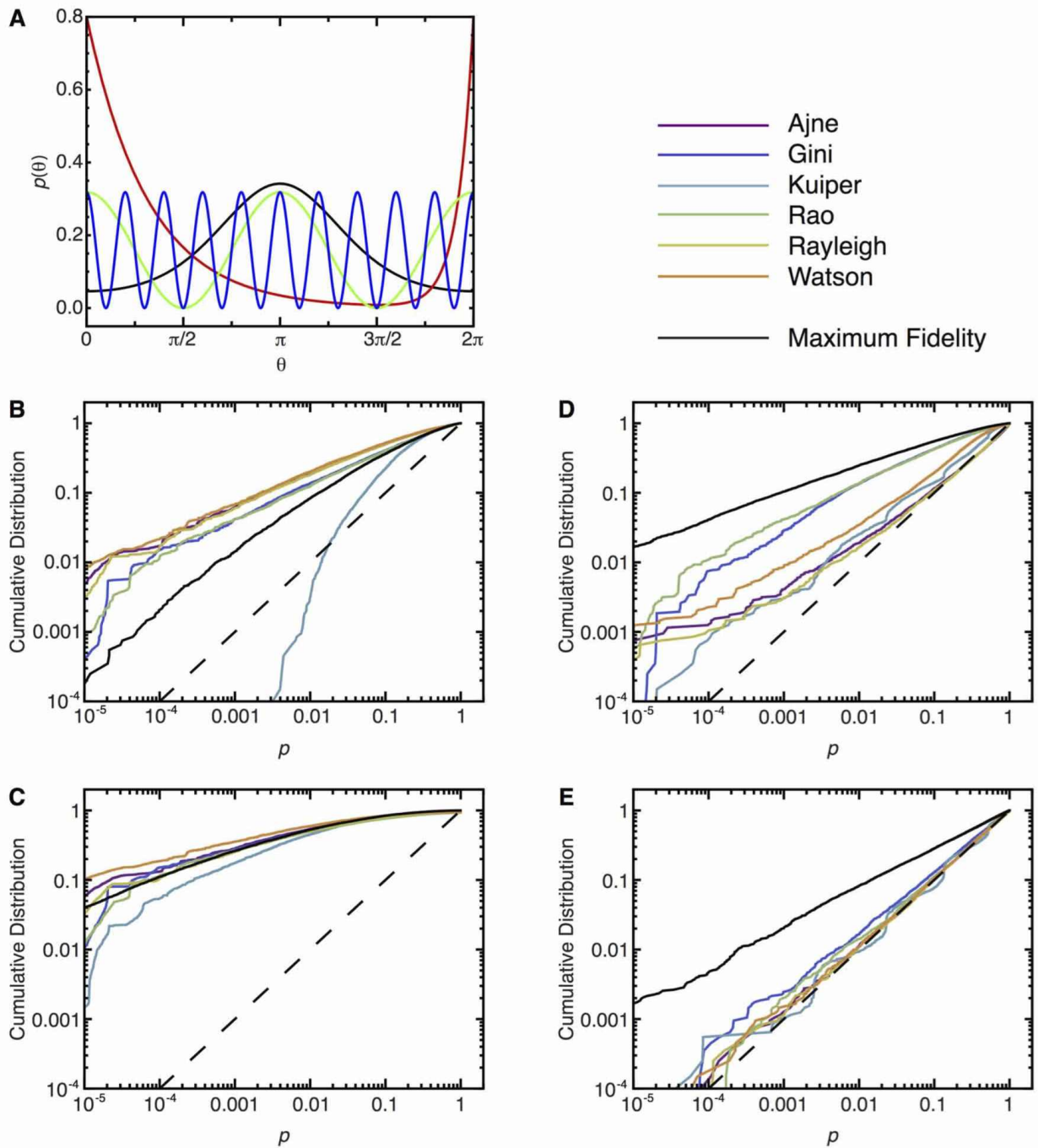}
     \vspace{-5.5cm}
      \caption{Concordance values for hypothesis testing on the circle. (A) Probability distributions used to generate data sets with $n=10$ points for testing against the assumption of a uniform distribution. (B) Cumulative distributions of $p$ values obtained using the the indicated goodness-of-fit tests for 1000 realizations of $n=10$ points drawn from a von Mises distribution with $\alpha=1$, $\beta=0$ (black curve in A) and tested against the assumption that the data were drawn from a uniform distribution. In the further panels, cumulative distributions of $p$ values similarly obtained are shown for data drawn from (C) a Wrapped Laplace model with $\alpha=0.5$ and $\beta=2$ (red curve in A),  (D) a cosine distribution: $\frac{1}{2\pi}\left(1+\cos{2\theta}\right)$ (green curve in A), and (E) another cosine distribution: $\frac{1}{2\pi}\left(1+\cos{10\theta}\right)$ (green curve in A).}
      \label{fig:good_circle}
   \end{center}
\end{figure}

A comparison of different goodness-of-fit measures applied to multiple circular distributions is presented in Fig.~\ref{fig:good_circle}. 
For distributions on the circle, a good reference distribution is simply the uniform distribution. For the cumulative curves shown in Figs.~\ref{fig:good_circle}B--E, $n=10$ points were repeatedly drawn by Monte Carlo (10000 realizations) from the displayed distributions in Fig.~\ref{fig:good_circle}A and compared against the assumption of a uniform distribution. How well a particular distribution can distinguish these deviations from uniformity is assessed by comparing the cumulative distribution of the $p$ values obtained by each approach. A distribution that can discriminate well should lie towards the upper left of these panels. Note that most of the tests have more power than the fidelity for discrimination of a ``location'' parameter (Fig.~\ref{fig:good_circle}B). However, the fidelity has equivalent or higher power for discrimination of shape parameters (Fig.~\ref{fig:good_circle}C--E). For the extreme case of a highly oscillating distrbution shown (Fig.~\ref{fig:good_circle}E), the fidelity performs much better, with the other statistics having essentially no discriminatory power. That the fidelity has the power to discriminate in the manner shown in Fig.~\ref{fig:good_circle}E is due to its consideration of the \textit{local} information of the distribution. This sensitivity to local information and thus more general discrepancies, however, comes at a slight cost for the specific discrimination of differences in ``location'' shown in Fig.~\ref{fig:good_circle}B (see further discussion in \S\ref{sec:discussion}).

\subsection{Goodness-of-fit on the Line}

\begin{figure}[]
   \begin{center}
     \vspace{-2cm}
     \includegraphics*[width=6.5in]{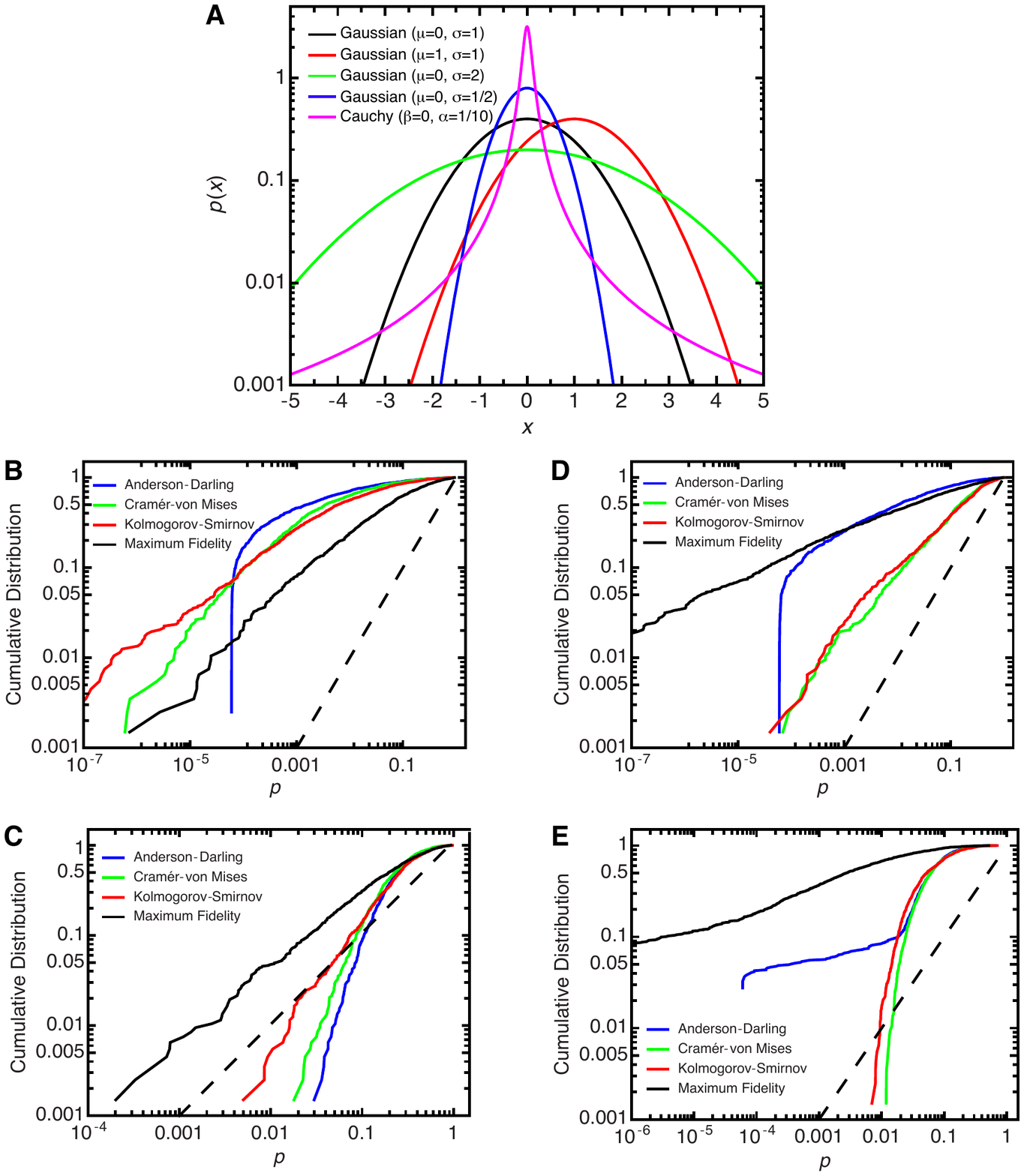}
     \vspace{-3.8cm}
      \caption{(A) Data points ($n=10$) were repeatedly drawn (1000 realizations) from the displayed distributions (red, green, blue, magenta) and tested against the assumption of a Gaussian model with $\mu=0$ and $\sigma=1$ (black). Cumulative distributions of the obtained $p$ values are shown using the indicated goodness-of-fit tests for data points drawn from (B) a Gaussian with $\mu=1$ and $\sigma=1$ (red model in A), (C) a Gaussian with $\mu=0$ and $\sigma=2$ (green model in A), (D) a Gaussian with $\mu=0$ and $\sigma=1/2$ (blue model in A), and (E) a Cauchy model with $\beta=0$ and $\alpha=1/10$ (magenta model in A).}
      \label{fig:good_Gauss}
   \end{center}
\end{figure}

In Fig.~\ref{fig:good_Gauss}, multiple methods for determining goodness-of-fit are compared for different distributions on the line. To perform this comparison, data points were drawn from the different model distributions shown in Fig.~\ref{fig:good_Gauss}A and compared against the assumption of a specific Gaussian model with $\mu=0$ and $\sigma=1$ (black curve in Fig.~\ref{fig:good_Gauss}A). In Fig.~\ref{fig:good_Gauss}B, data were drawn from a shifted Gaussian, with the fidelity exhibiting less discriminatory power than the other tests. In Fig.~\ref{fig:good_Gauss}C, data were drawn from a broader Gaussian with the same mean. In Fig.~\ref{fig:good_Gauss}D, data were drawn from a narrower Gaussian. In Fig.~\ref{fig:good_Gauss}E, data were drawn from a Cauchy distribution. In Figs.~\ref{fig:good_Gauss}C--E, the fidelity provides the best discrimination. For the slight drop in power for ``location'' testing (Fig.~\ref{fig:good_Gauss}B), one gains much greater power in discriminating more general, shape-dependent differences. In the opposite sense, the sensitivity of the other tests to ``location'' differences, comes at the cost of identifying other types of discrepancies between the distributions. On the whole, the fidelity offers a more balanced test for general differences between distributions (see further discussion in \S\ref{sec:discussion}).

%%%%%%%%%%%%%%%%%%%%%%%%%%%%%%%%%%%%%%%%%%%%%%%%%%%%%%%%%%%%%%
\section{Joint Analysis}\label{sec:joint}
%%%%%%%%%%%%%%%%%%%%%%%%%%%%%%%%%%%%%%%%%%%%%%%%%%%%%%%%%%%%%%

Highly accurate joint analysis of separate data sets on the line or the circle can be performed based on the following definitions (see Ranneby (1984) for a similar formalism for the spacings statistic\cite{ranneby_maximum_1984}):
\begin{eqnarray}
%F_{\vec{n}}&=&\sum_{i=1}^k F_{n_i}\\
%M_{\vec{n}}&=&\sum_{i=1}^k  M_{n_i}\\
%\Sigma^2_{\vec{n}}&=&\sum_{i=1}^k \Sigma^2_{n_i}\\
f_{\vec{n}}&=&\frac{1}{n}\sum_{i=1}^k n_i f_{n_i}\\
\mu_{\vec{n}}&=&\frac{1}{n}\sum_{i=1}^k n_i \mu_{n_i}\\
\sigma^2_{\vec{n}}&=&\frac{1}{n^2}\sum_{i=1}^k n_i^2 \sigma^2_{n_i},
\end{eqnarray}
where $n=\sum_{i=1}^k n_i$.  The $p$ value in this generalized case is:
\begin{equation}\label{eq:joint}
p_{\vec{n}}(f_{\vec{n}}) = \int_{-\infty}^{f_{\vec{n}}}P_{\vec{n}}(x)dx\simeq Q(\mu_{\vec{n}}^2/\sigma_{\vec{n}}^2,\mu_{\vec{n}} f_{\vec{n}}/\sigma_{\vec{n}}^2).
%=\frac{\Gamma(\mu_n^2/\sigma_n^2,\mu_n f_n/\sigma_n^2)}{\Gamma(\mu_n^2/\sigma_n^2)}.
\end{equation}
This approximation works extremely well for all combinations of circular and linear data (including combinations of circular data with linear data) with only the following three exceptions:  $n_1=1$ (line) with $n_2=1$ (line), $n_1=1$ (line) with $n_2=2$ (circle), and $n_1=2$ (circle) with $n_2=2$ (circle).  For these cases, a power law should be used based on the mean value obtained above to determine the $p$ value:
\begin{equation}
p_{\vec{n}}(f_{\vec{n}}) \simeq \exp(-f_{\vec{n}}/\mu_{\vec{n}}).
\end{equation}
This is similar to the need for a power-law approximation for single data sets having either $n=2$ points on the line or $n=3$ points on the circle. Useful examples of joint analysis can be found in \S\ref{sec:tests} in the context of generalized ``$t$ tests'' for comparing two data sets and in \S\ref{sec:multidim} for the case of multidimensional data.

\begin{figure}[]
   \begin{center}
     \vspace{0cm}
     \includegraphics*[width=5in]{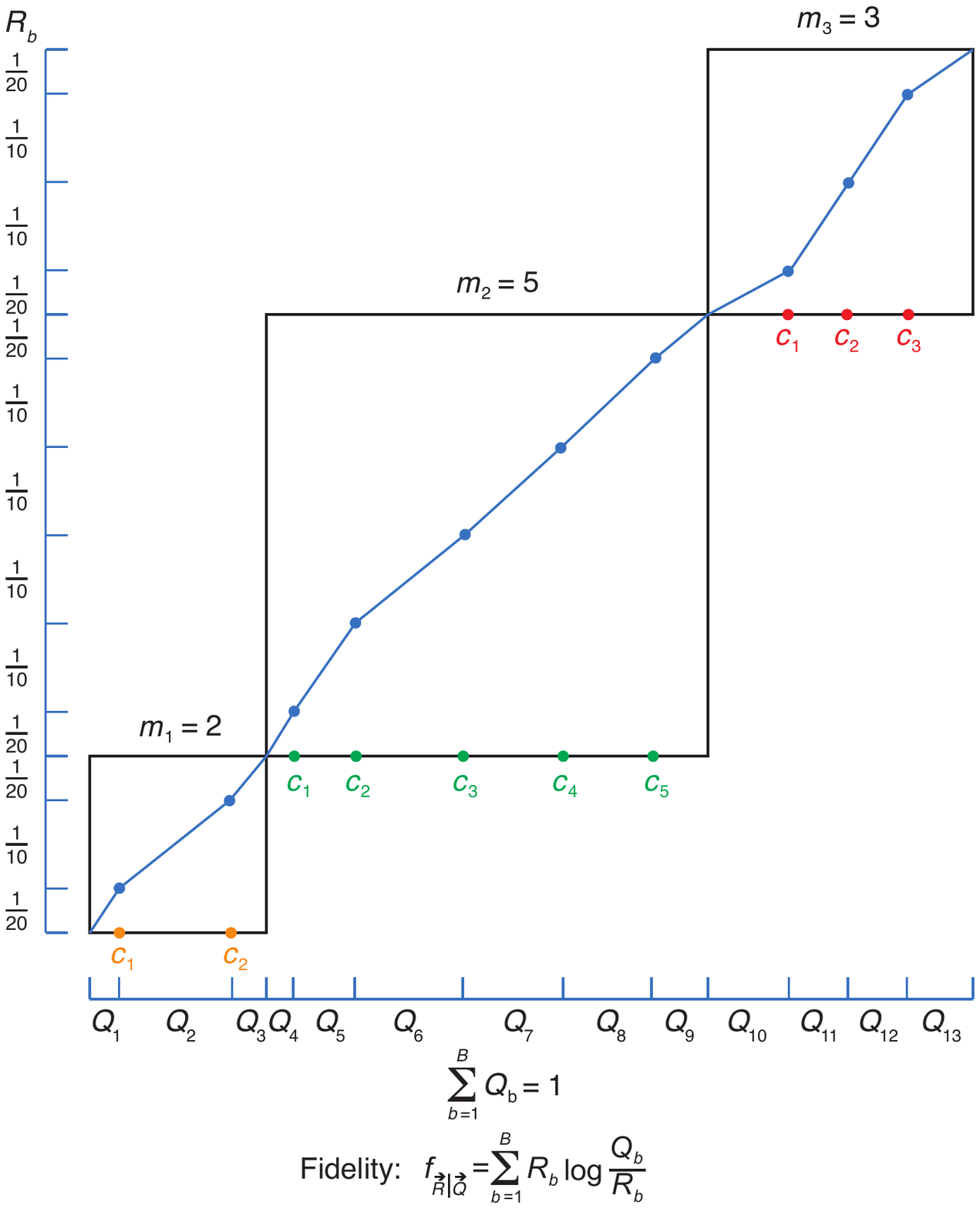}
     \vspace{-7cm}
      \caption{Illustration of joint analysis. Three independent data sets are shown with $m_1=2$, $m_2=5$, and $m_3=3$ points. An hypothesized model is simultaneously fit to all three data sets leading to the model-derived cumulative values, $c_i$. Here we assume that the data cannot simply be combined together into a single data set containing all of the points due to the different circumstances under which each individual data set was taken. Each point contributes equal weight to the final fidelity, with boundary intervals contributing the usual half weight. The fidelity is computed according to the displayed equation (same as Eq.~\ref{eq:fidelity}), with the $p$ value obtained according to Eq.~\ref{eq:joint}.}
      \label{fig:joint}
   \end{center}
\end{figure}

%%%%%%%%%%%%%%%%%%%%%%%%%%%%%%%%%%%%%%%%%%%%%%%%%%%%%%%%%%%%%%
\section{Generalization of Student's $t$ Test}\label{sec:tests}
%%%%%%%%%%%%%%%%%%%%%%%%%%%%%%%%%%%%%%%%%%%%%%%%%%%%%%%%%%%%%%

One of the most commonly used tests in statistics is Student's $t$ test \cite{student_probable_1908}, which was developed by William Gosset to determine the likelihood that two separate data sets (drawn from Gaussian distributions with equal variance) have the same mean. Student's $t$ test was subsequently generalized in many different ways, e.g. to treat the case of unequal variances by Welch (Welch's $t$ test)\cite{welch_generalization_1947} or to treat the comparison of more than two data sets (Fisher's $F$ test, which underlies ANOVA). The important assumption in all of these tests is that all of the true, underlying distributions are Gaussian. In the language of maximum fidelity, these tests are all consistently treated in the same manner, with, importantly, no assumption of Gaussianity required. This is possible because the fidelity and associated $p$ value can be applied to arbitrary distributions to, in this case, allow direct assessment of the concordance of the complete model (comprised of the individual models for each data set) to the complete data (all data sets at once).

\begin{figure}[]
   \begin{center}
     \vspace{-1cm}
     \includegraphics*[width=6.5in]{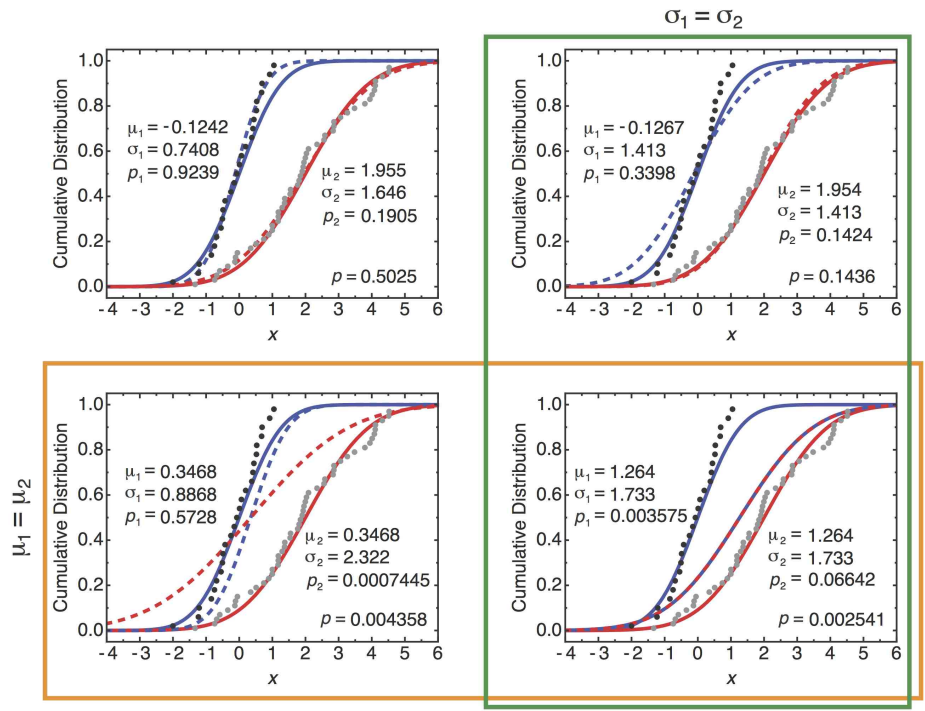}
     \vspace{-1.5cm}
      \caption{Maximum fidelity generalization of $t$ tests applied to Gaussian distributions. Two data sets of $n_1=25$ and $n_2=50$ data points were drawn from two different Gaussian distributions having $\mu_1^0=0$ and $\sigma_1^0=1$ (solid blue curve) or $\mu_2^0=2$ and $\sigma_2^0=1.5$ (solid red curve). In each panel, different assumptions are made to fit the data. In the upper left, different $\mu$ and $\sigma$ were independently fit to each data set to find the complete model that maximizes the overall fidelity. In the upper right, different $\mu$ were assumed but the same $\sigma$ to find the complete model that maximizes the overall fidelity. In the lower left, different $\sigma$ were assumed but the same $\mu$. In the lower right, both $\mu$ and $\sigma$ were assumed to be the same for each data set. In each panel, $p_1$ and $p_2$ correspond to the concordance of the individual model fits and $p$ corresponds to the overall concordance of both models to the data.}
      \label{fig:gaussian_t_test}
   \end{center}
\end{figure}

Consider two independent data sets with $n_1$ and $n_2$ data points. We would like to compare the concordance of different hypothesized complete models to the complete data. We obtain three informative $p$ values, $p_1$ for the model fit to the first data set alone, $p_2$ for the model fit to the second data set alone, and $p$ for the complete model fit to both data sets together (using the joint approach explained in \S\ref{sec:joint}). First, we hypothesize that the points in each data set are drawn from the same underlying Gaussian distribution specified by $\mu$ and $\sigma$.  A parameter space search determines the solution that maximizes the fidelity.  The concordance ($p$ value) for this $\mu_0$ and $\sigma_0$ is recorded (Fig.~\ref{fig:gaussian_t_test}, lower right).  A second hypothesis could be that both distributions have the same $\sigma$ but different means $\mu_1$ and $\mu_2$ (Fig.~\ref{fig:gaussian_t_test}, upper right). A third hypothesis could test how much improvement occurs upon assuming the same mean $\mu$ but different $\sigma_1$ and $\sigma_2$ (Fig.~\ref{fig:gaussian_t_test}, lower left).  Finally, we can test if distinct $\mu_1$ and $\mu_2$ and $\sigma_1$ and $\sigma_2$ provide a significant improvement in $p$ value (Fig.~\ref{fig:gaussian_t_test}, upper left). These plots are only shown for the maximum value of $p$ (as a \textit{representative} value of the full concordance landscape) with no correction for the ``degrees of freedom'' implied by the assumption, see \S\ref{sec:concordance} for further discussion of this issue.

\begin{figure}[]
   \begin{center}
     \vspace{-1cm}
     \includegraphics*[width=6.5in]{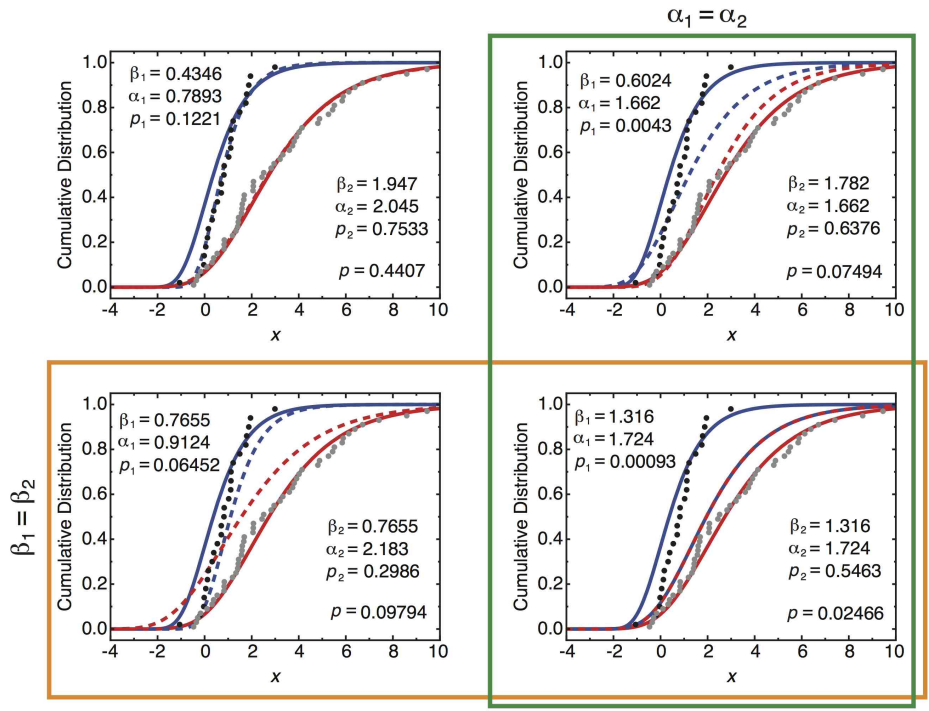}
     \vspace{-1.5cm}
      \caption{Maximum fidelity generalization of $t$ tests applied to Extreme Value distributions.  Two data sets of $n_1=25$ and $n_2=50$ data points were drawn from two different Extreme Value distributions having $\beta_1^0=0$, $\alpha_1^0=1$ (solid blue curve) or $\beta_2^0=2$, $\alpha_2^0=2$ (solid red curve). In each panel, different assumptions are made to fit the data. In the upper left, different $\beta$ and $\alpha$ were independently fit to each data set to find the complete model that maximizes the overall fidelity. In the upper right, different $\beta$ were assumed but the same $\alpha$ to find the complete model that maximizes the overall fidelity. In the lower left, different $\alpha$ were assumed but the same $\beta$. In the lower right, both $\beta$ and $\alpha$ were assumed to be the same for each data set. In each panel, $p_1$ and $p_2$ correspond to the concordance of the individual model fits and $p$ corresponds to the overall concordance of both models to the data.}
      \label{fig:extreme_value_t_test}
   \end{center}
\end{figure}

To demonstrate that maximum fidelity is not limited to merely the comparison of Gaussian distributions, a similar example is shown in Fig.~\ref{fig:extreme_value_t_test} based on data drawn from two different Extreme Value distributions. It is important to note that mixed cases are also treated in the exact same manner: the first data set could be assumed to be drawn from a Gaussian distribution and the second from an Extreme Value distribution. The associated $p$ value will always be calculated for the combined fit in exactly the same way and with the same meaning.

%%%%%%%%%%%%%%%%%%%%%%%%%%%%%%%%%%%%%%%%%%%%%%%%%%%%%%%%%%%%%%
\section{Neyman paradox}\label{sec:neyman}
%%%%%%%%%%%%%%%%%%%%%%%%%%%%%%%%%%%%%%%%%%%%%%%%%%%%%%%%%%%%%%

\begin{figure}[]
   \begin{center}
     \vspace{-4cm}
     \includegraphics*[width=5.5in]{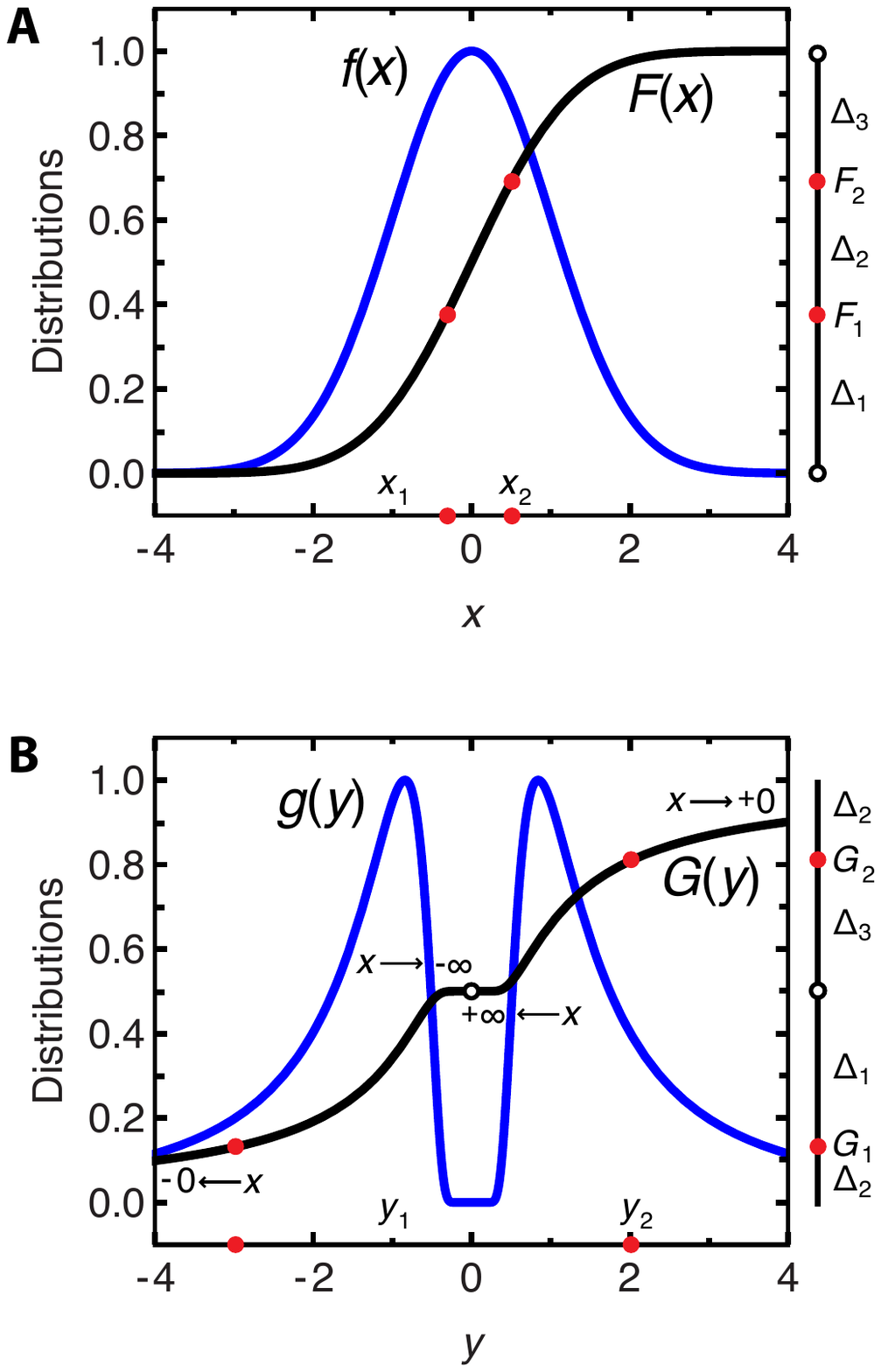}
     \vspace{-4.5cm}
      \caption{Resolution of the Neyman paradox in the context of maximum fidelity. (A) A Gaussian distribution $f(x)$ with $\mu=0$ and $\sigma=1$ and associated cumulative distribution $F(x)$ are displayed, along with two points $x_1=-1/3$ and $x_2=1/2$ and their associated cumulative values $F_1$ and $F_2$ shown on the right. The boundaries at $x\rightarrow-\infty$ and $x\rightarrow+\infty$ are displayed on the cumulative interval as open circles. (B) The transformation of $y=1/x$ leads to the transformed probability distribution $g(y)$ and cumulative distribution $G(y)$. The original $x_1$ and $x_2$ points are transformed to $y_1=-3$ and $y_2=2$. The boundaries at $x\rightarrow\pm\infty$ are mapped to $y=0$ (open circle). Proper tracking of the transformed boundaries preserves the value of the fidelity, defined by the cumulative intervals $\Delta_1$, $\Delta_2$, and $\Delta_3$.}
      \label{fig:neyman_i}
   \end{center}
\end{figure}

It is important to address a philosophical and mathematical argument made by Neyman, generally referred to as ``Neyman's paradox''\cite{neyman_lectures_1952,redhead_neymans_1974}. This paradox involves the supposed nonuniqueness of the \textit{central} interval for a particular quantity estimated from the data, such as the mean or median (though this argument can be generalized to other parameters as well). Neyman argued that the transformation $y=1/x$ (as shown in Fig.~\ref{fig:neyman_i}) moves points originally near the center of the distribution to the edges of the transformed distribution (and vice versa), leading to completely different values for the mean and median of the data, which, if used for parameter estimation or confidence intervals, would lead to a radically different model preference. 

The fidelity affords a simple resolution of Neyman's paradox. Neyman's mathematical \textit{sleight of hand} neglects one important aspect: the topology of the line. This is pictorially demonstrated in Fig.~\ref{fig:neyman_i}. Upon the transformation of $y=1/x$, the Gaussian probability distribution $f(x)$ and its associated cumulative distribution $F(x)$ shown in Fig.~\ref{fig:neyman_i}A are transformed to $g(y)$ and $G(y)$ in Fig.~\ref{fig:neyman_i}B. The two points near the center of the distribution $f(x)$ are moved toward the extremities of the transformed distribution $g(y)$. Imagine if we had only observed the single point $x_1$ in Fig.~\ref{fig:neyman_i}A. The closeness of $x_1$ to the median value of $F(x)$ implies a good fit according to the fidelity; however, the transformed value $y_1$ is moved to the extreme left of the transformed distribution, which appears at first glance to indicate a poor fit. This neglects the topological change that the transformation $y=1/x$ generates. The boundaries in Fig.~\ref{fig:neyman_i} at $x\rightarrow-\infty$ and $x\rightarrow+\infty$ are moved to $y=0$ in the transformed system. If the topology is preserved upon the transformation, $y=0$ does not exist as it represents the new boundary. The white circle at $y=0$ denotes the non-existent, boundary value nature of this point, with the values as $y\rightarrow+\infty$ wrapping around to $y\rightarrow-\infty$. Proper tracking of the boundaries and the cumulative intervals defined in the original and the transformed coordinates leads to preservation of the exact value of the fidelity defined by the intervals $\Delta_1$, $\Delta_2$, and $\Delta_3$ displayed on the right sides of Figs.~\ref{fig:neyman_i}A and B. If we are not sure where the boundaries fall on the line, then we should assume no boundaries, i.e. a circular coordinate system. In this case as well, the fidelity would be trivially invariant to such a coordinate inversion.

\begin{figure}[]
   \begin{center}
     \vspace{-3cm}
     \includegraphics*[width=4in]{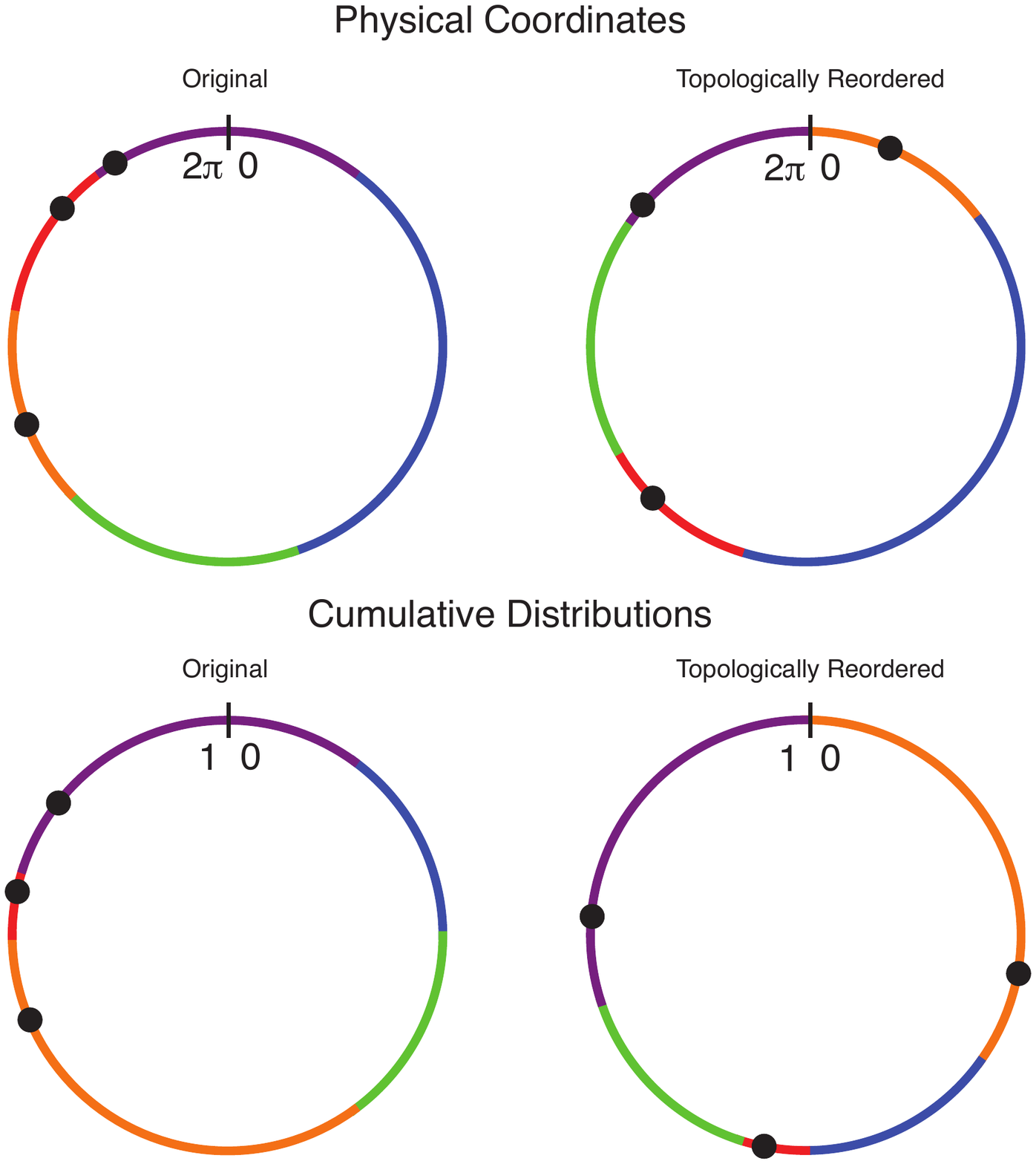}
     \vspace{-1.5cm}
      \caption{Arbitrary topological reordering of a physical coordinate system could be used to make an observed data set fit any desired model.  In the top panel, data on the circle are expressed in their original coordinate system and in a topologically reordered coordinate system. In the bottom panel, a particular model-based cumulative mapping of these data points leads to a poor fit for the original topology but a much improved fit on the reordered topology.}
      \label{fig:neyman_ii}
   \end{center}
\end{figure}

\begin{figure}[]
   \begin{center}
     \vspace{-6cm}
     \includegraphics*[width=5in]{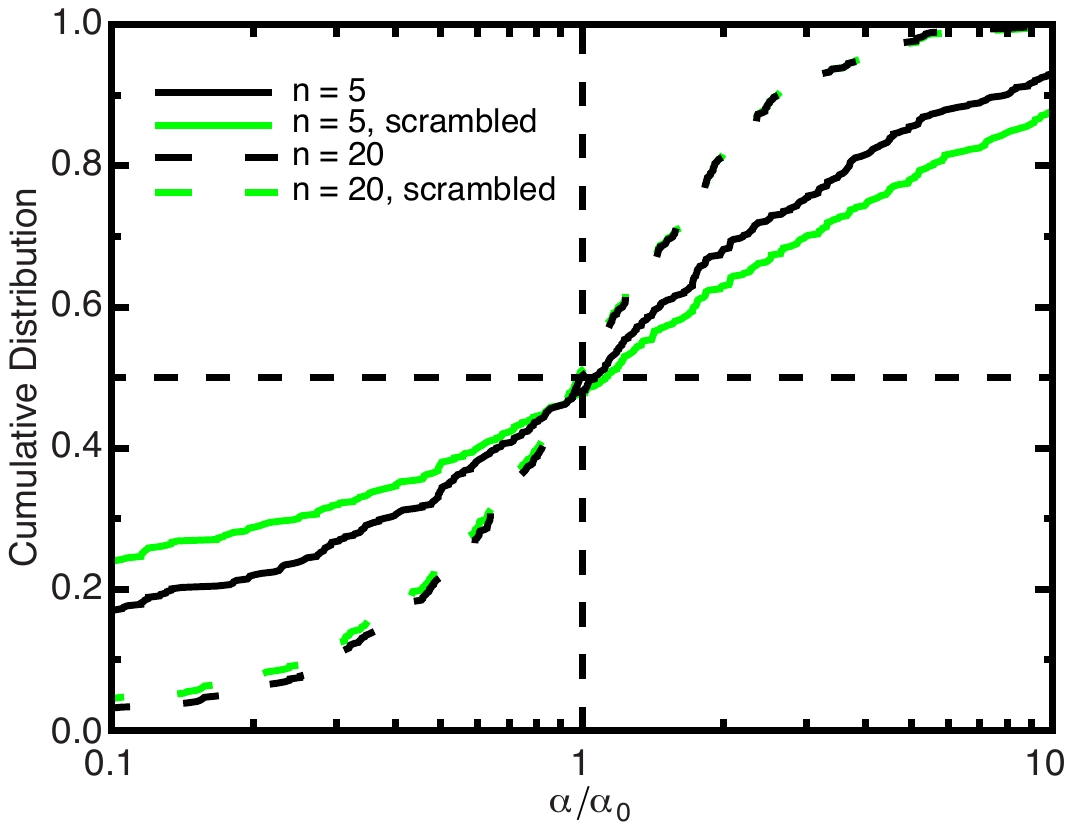}
     \vspace{-6cm}
      \caption{Parameter estimation of a model (Theta distribution with $\alpha_0=4$ and fixed $\beta=0$, see Table~\ref{tab:circ_dist}) carried out on the original topology of the physical coordinate system and on the scrambled topology shown in the top panel of Fig.~\ref{fig:neyman_ii}. If the topology is scrambled before the data are observed, parameter estimation and $p$ value goodness-of-fit can still be carried out with only a slight reduction of efficiency for the toploogically reordered system (with this reduction likely due to the model discontinuities introduced by the topological reordering).}
      \label{fig:neyman_iii}
   \end{center}
\end{figure}

To further illustrate why the topology should be preserved upon a change in the data coordinate system, consider Fig.~\ref{fig:neyman_ii}. In the upper panels, we have a circular coordinate system (for simplicity) with three data points in its natural coordinate system and in a topologically reordered system obtained by ordering the colored regions of the circle differently. With regard to the cumulative distribution, a fit that looks poor in the original topology can be made to look much better in a topologically reordered system. Arbitrary topological reordering can make any model fit any data set perfectly. Such topological reordering as carried out in Fig.~\ref{fig:neyman_ii} or topological redefinition of the boundaries (as Neyman carries out, thus making the points at $x\rightarrow+\infty$ directly neighbor the points at $x\rightarrow-\infty$) should not be permitted in a proper theory of inductive inference. If the topology were scrambled before the data were examined, however, the correct model parameters could \textit{still} be estimated with almost as efficiently as in the unaltered coordinate system (Fig.~\ref{fig:neyman_iii}). As I am not aware of any justification for such topological reordering (or boundary reassignment) of real-world data, it seems safe to conclude that such rearrangements should be restricted for the proper application of statistical inference. In other words, we should at least be able to agree on the overall topology of the data coordinate system, even if we cannot agree on the best choice of topology-preserving coordinate system. If we cannot agree on this restriction, then we must indeed embrace Neyman's argument as an argument against a fundamental basis for induction, leaving us with only Neyman's unsatisfactory notion of good \textit{inductive behavior}\cite{neyman_lectures_1952}.

%%%%%%%%%%%%%%%%%%%%%%%%%%%%%%%%%%%%%%%%%%%%%%%%%%%%%%%%%%%%%%
\section{Binned Data}\label{sec:binned}
%%%%%%%%%%%%%%%%%%%%%%%%%%%%%%%%%%%%%%%%%%%%%%%%%%%%%%%%%%%%%%

\begin{figure}[]
   \begin{center}
     \vspace{-8.5cm}
     \includegraphics*[width=6.5in]{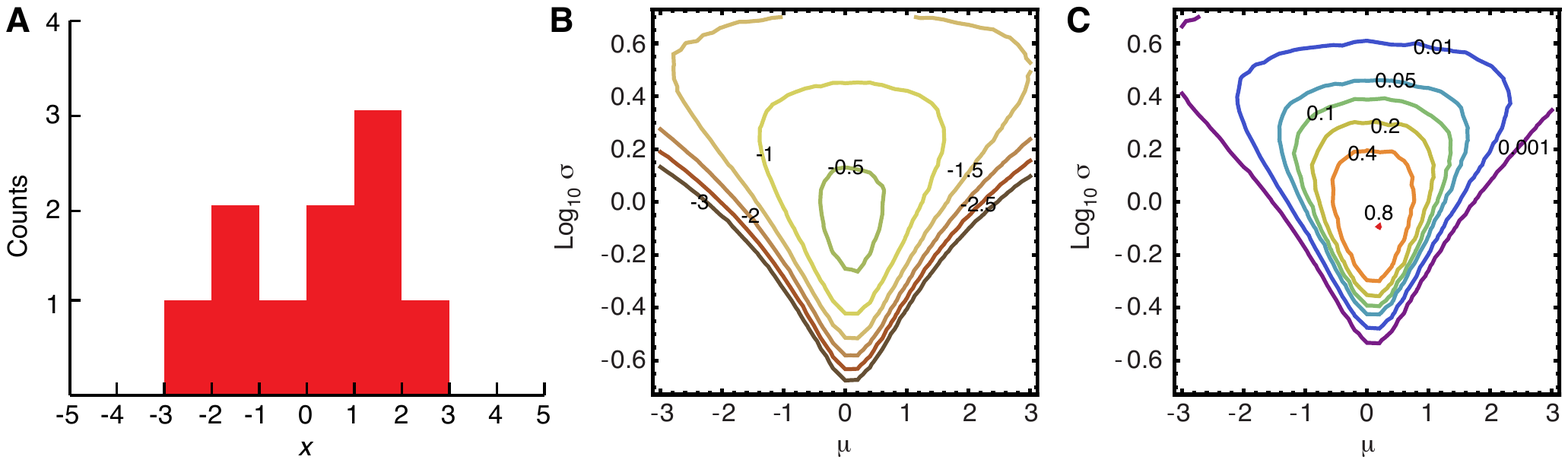}
     \vspace{-10cm}
      \caption{Fitting of a Gaussian model to binned data. (A) Binned data. (B) Contour map of the fidelity. (C) Contour map of the associated $p$ values. See text for further details.}
      \label{fig:binned}
   \end{center}
\end{figure}

The generalization of maximum fidelity to univariate binned data is straightforward but requires more computation. Binned data hypothesized to come from a Gaussian are given in Fig.~\ref{fig:binned}A. The fidelity landscape and corresponding concordance landscape shown in Figs.~\ref{fig:binned}B and C were generated in the following way. The landscapes were sampled at 30 values along the displayed abscissa and ordinate. At each value of $\mu$ and $\sigma$, the cumulative distribution was tabulated at the boundaries of all of the $x$ bins in Fig.~\ref{fig:binned}A that contain points. Within each of these bins, the points were randomly assigned cumulative values drawn from a uniform distribution spanning the cumulative boundaries. The fidelity was then calculated. This was repeated 999 times, with the \textit{representative} median fidelity value (the 500th ordered value) then displayed in Fig.~\ref{fig:binned}B. The fidelity was then converted to concordance ($p$ value) according to Eq.~\ref{eq:concordance}. With modern desktop computers, the generation of Figs.~\ref{fig:binned}B and C take only a fraction of a second. This is what I will refer to as the ``exact method'' for parameter estimation and concordance determination for binned data.

A more economical method valid only for parameter estimation would be to maximize a fidelity estimate obtained by calculating the cumulative values at the bin boundaries and then redistributing the data points in a way that maximizes the fidelity within that bin. One can then efficiently maximize this fidelity estimate without the tedious 999 Monte Carlo simulations for each trial choice of model parameter values required above. Specifically, for a bin $k$ containing $p$ data points, the data points are spread evenly over the cumulative range of the bin according to the following formula:
\begin{equation}\label{eq:fidspace}
c_{j+i} = \left(\frac{i-1/2}{p}\right)\left(c^u_k - c^l_k\right) + c^l_k, 
\end{equation}
where $i$ ranges over the list of data points (from $1$ to $p$), $j+1$ denotes the first data point in bin $k$ (out of the list of total data points), and $c^l_k$ and $c^u_k$ denote the lower and upper values of the model-specific cumulative distribution over the bin $k$. The fidelity is then calculated as usual (Eq.~\ref{eq:fidline}). This method for estimating the fidelity for binned data has already been used in the real-world context of determining fluorophore lifetimes and lifetime populations in Walther et al. (2011)\cite{walther_precise_2011} (where it was also shown that maximum fidelity allowed for better estimation of lifetime populations than maximum likelihood).

%%%%%%%%%%%%%%%%%%%%%%%%%%%%%%%%%%%%%%%%%%%%%%%%%%%%%%%%%%%%%%
\section{Binary Distributions}\label{sec:binary}
%%%%%%%%%%%%%%%%%%%%%%%%%%%%%%%%%%%%%%%%%%%%%%%%%%%%%%%%%%%%%%

Consider $k_0$ successes out of $n$ trials of a particular binary process. One could hypothesize an independent binomial process (Bernoulli process) to describe the observed set of events, with the probability of the particular outcome of $k_0$ successes given by the corresponding term in the binomial expansion:
\begin{equation}
\left(q+(1-q)\right)^n=\sum_{k=0}^{n}{n\choose k} q^k(1-q)^{n-k}.
\end{equation}
Here, $q$ denotes the probability of success for each individual event. Obtaining a measure of how adequate a particular $q$ value is (or a range of $q$ values) to account for the observed $k_0$ out of $n$ successes is a problem with a suprisingly long and convoluted history already partially explored in \S\ref{sec:intro}. The cumulative distribution values for a binomial distribution with success rate $q$ defined at the observed value of $k_0$ are:
\begin{eqnarray}\label{eq:binary}
c_l(n,q;k_0)&=&\sum_{k=0}^{k_0-1}{n\choose k} q^k(1-q)^{n-k},\\
c_m(n,q;k_0)&=&\frac{1}{2}{n\choose {k_0}} q^{k_0} (1-q)^{n-k_0}+\sum_{k=0}^{k_0-1}{n\choose k} q^k(1-q)^{n-k},\\
c_h(n,q;k_0)&=&\sum_{k=0}^{k_0}{n\choose k} q^k(1-q)^{n-k},
\end{eqnarray}
which have the following meanings: $c_l(n,q;k_0)$ (``low'') corresponds to the cumulative value obtained from summing the contributions from $k=0$ to one less than the observed value $k=k_0-1$, $c_h(n,q;k_0)$ (``high'') corresponds to the cumulative value obtained from summing the contributions from $k=0$ to $k=k_0$, and $c_m(n,q;k_0)$ (``middle'') corresponds to the midpoint value of the cumulative interval defined for the particular outcome $k_0$ ($c_m$ is the average of $c_l$ and $c_h$). Under the assumption of a binomial process, the particular binomial distribution (parametrized by $q$) that maximizes the fidelity is the one that centralizes the cumulative distribution at the observed value of $k_0$. This is equivalent to maximizing the fidelity for a data set containing a single point, in this case $k_0$ on the discretized interval from $k=0$ to $n$. This requires solving:
\begin{equation}
c_m(n,q;k_0)=\frac{1}{2}
\end{equation}
for $q$. For example, for $n=10$ and $k_0=3$, the value of $q$ that maximizes the fidelity is $q\approx0.306089$ (which differs from the value $q=0.3$ obtained from maximum likelihood). Due to the fact that $k_0$ can be located at one of the boundaries ($k=0$ or $k=n_0$), it is best to retain the cumulative value (or interval) at the point $k_0$ as an indicator of goodness-of-fit rather than convert to a \textit{symmetrized} $p$ value (see further discussion below). The cumulative function value(s) determined over the observed bin $k_0$ already contains all of the information we need (in fact, it contains even more information than the $p$ value, as the particular left-right asymmetry of the cumulative value from its optimal value of 1/2 is preserved).

\begin{figure}[]
   \begin{center}
\vspace{0cm}
     \includegraphics*[width=3.5in]{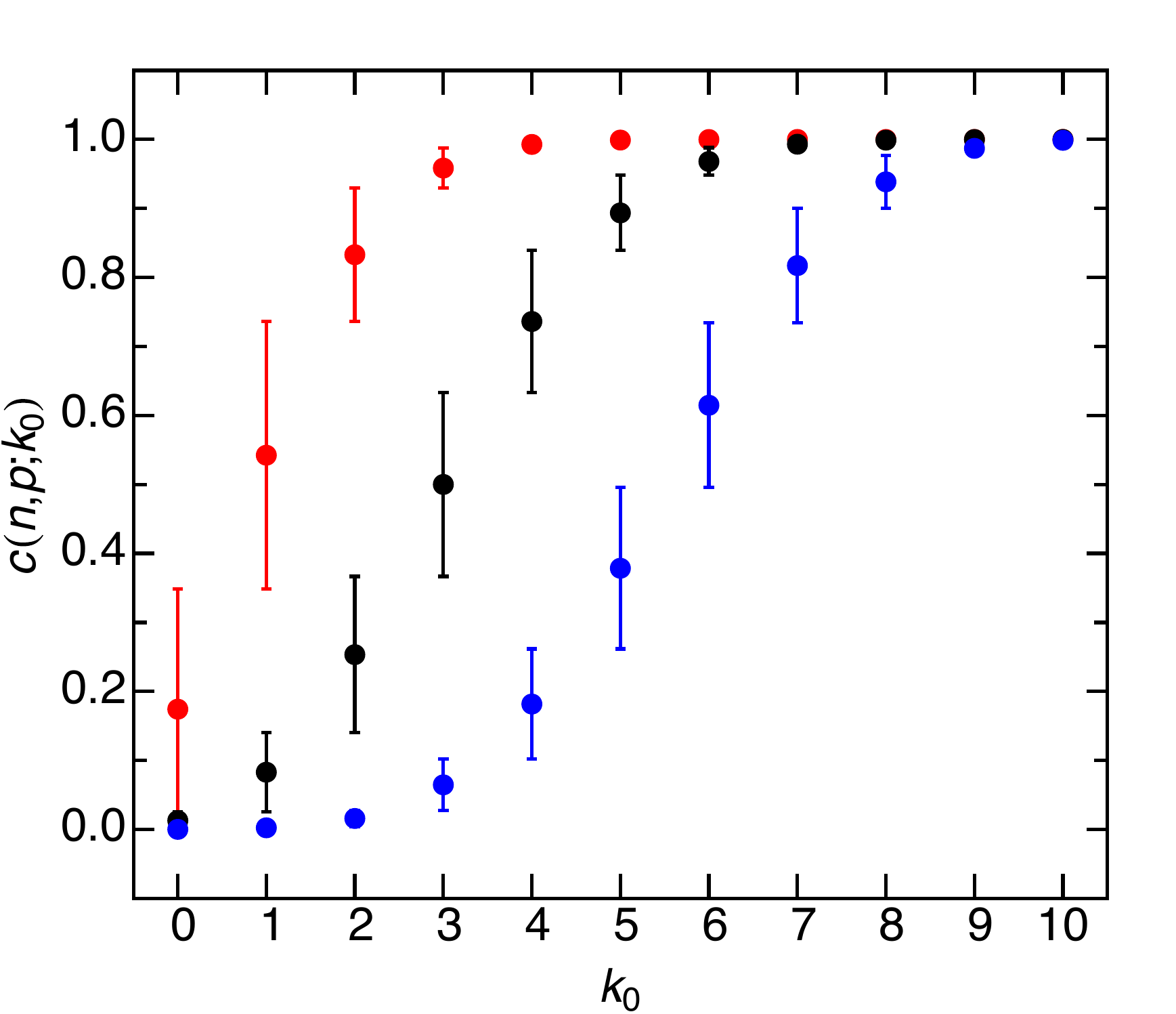}
\vspace{0cm}
      \caption{Binomial cumulative distributions ($n=10$) for $q=0.1$ (red), 0.306089 (black), and 0.55 (blue). The cumulative interval for a given number of successes $k_0$ is defined by its midpoint value given by $c_m(10,q;3)$ and its symmetric lower  and upper extents respectively defined by $c_l(10,q;3)$ and $c_h(10,q;3)$.}
      \label{fig:binary}
   \end{center}
\end{figure}

The cumulative distributions $c_m(n,q;k_0)$ (along with the lower and upper extent of the cumulative intervals respectively given by $c_l(n,q;k_0)$ and $c_h(n,q;k_0)$) are shown as a function of $k_0$ for $n=10$ and for different values of the probability $q$ in Fig.~\ref{fig:binary}. Consider the outcome $k_0=3$. For the $q=0.306089$ distribution (black intervals), the cumulative interval at $k_0$ is centered at 0.5, symmetrically spanning the values $c=0.367$ to $0.633$. Determination of goodness-of-fit (measured by the cumulative interval) for other values of $q$ proceeds as follows. For the assumption of $q=0.1$ (red intervals), the outcome of $k_0=3$ would occur over the cumulative interval spanning $c=0.930$ to $0.987$. For the assumption of $q=0.55$ (blue intervals), the cumulative interval at $k_0=3$ spans $c=0.027$ to $0.102$. Because both distributions map the value $k_0=3$ to cumulative intervals at the extremities of the cumulative distribution, they can be considered poor fits. The central 90\% ``confidence interval'' over which the full cumulative bin at $k_0=3$ is never greater than $c=0.95$ ($c_h(10,q;3)=0.95$) or less than $c=0.05$ ($c_l(10,q;3)=0.05$) corresponds to the range $q=0.150$ to $0.507$ (the Clopper-E.~Pearson ``exact'' interval\cite{clopper_use_1934}). One could alternatively calculate the values at which $k_0=3$ is centered at $c=0.95$ and $c=0.05$, which depends only on $c_m(10,q;3)$. This would correspond to the slightly larger range of $q=0.107$ to $0.571$.

Consider what happens if we observe $k_0=0$ out of $n=10$ events. If $q=0$ is assumed (the lowest value possible), the cumulative value at $k_0=0$ will be $c_m(10,0;0)=0.5$ and can take on no larger value. By retaining the cumulative mapping, this $k_0=0$ example merely demonstrates that at the lower boundary, only extreme cumulative values (intervals) at $k_0$ at the \textit{lower} range of the confidence interval are necessary to consider (for increasingly larger $q>0$). This important asymmetry of the above-constructed confidence intervals is missing in the Clopper-E.~Pearson, who attempt to define a central 90\% confidence level no matter the value of $k_0$, which will clearly fail at the boundaries, as pointed out by Wilson\cite{wilson_confidence_1942}. These boundary concerns can be neglected for continuous data, as data points that would fall directly at the boundaries of a continuous interval would occupy a set of measure zero.
Wilson's preferred solution\cite{wilson_probable_1927} is to consider for an hypothesized distribution all outcomes with probabilities greater than or equal to the probability to obtain $k_0$. This is akin to a likelihood approach. Because this is carried out on discrete data, the likelihood has an invariant meaning, as coordinate changes, which affect the likelihood value for continuous data, do not exist. The extension of Wilson's approach to continuous distributions would therefore correspond to a likelihood ratio analysis, which should be rejected for the reasons given in \S\ref{sec:intro} for the \textit{likelihood fallacy}. Wilson's approach for this discrete situation is in general not as ideal or as intuitive as use of the cumulative interval, as different distributions have different shapes, with the likelihood being a better or worse discriminating statistic depending critically on these shape considerations. Consider the following extreme example of a flat (and therefore non-Bernoulli) hypothesized probability distribution over a cluster of $k$ values that can be shifted to cover $k_0$. Wilson's approach would give a similar concordance if $k_0$ lies anywhere within the flat distribution, whereas the cumulative method outlined above would give a preference for $k_0$ to lie at the center of the flat distribution (representing a better estimate of this shift parameter), with an increasingly worse concordance as $k_0$ is shifted to either extremity. This same criticism would apply to maximum likelihood on continuous distributions. The shape-dependent discriminatory power of maximum likelihood is therefore not as general as the fidelity, for which model discrimination is based on the universal nature of the cumulative distribution.

An example of a binary process is Peirce's consideration of the occurrence of the number 5 in the sequence of $\pi$ (Peirce 7.121-122)\cite{peirce_collected_1974}. He found that 5 occurred 33 times in the first 350 digits of $\pi$ and 28 times in the next 350 digits, with both values within the statistical expectation for a binomial process with $q=0.1$ for the occurrence of $5$ (success) versus the other 9 digits. Consider the even simpler case of what the fidelity tells us about the occurrence of 1's in the expression of $\pi$ in binary units, $\pi=11.00100100001111110110101010001000\ldots$ Within the first 1000 digits of $\pi$, 1's appear 489 times. The cumulative value for this outcome assuming a binomial process with $q=0.5$ is $c_m(1000,0.5;489)=0.243$. The value of $q$ that centralizes the cumulative distribution at $k_0=489$ is, unsurprisingly, $q\simeq0.489$. More interestingly, the lower and upper limits for $q$ for which the midpoint cumulative value attains either $c=0.95$ or $c=0.05$ based on $c_m(1000,q;489)$ are, respectively, $q\simeq0.463$ and $q\simeq0.515$.

All of the above pertains to the assumed model of a binomial process of independent events with fixed probability $q$ for each event (Bernoulli process). However, a given observed sequence of binary events need not be generated by such a process, as Wilson\cite{wilson_probable_1927,wilson_comparative_1964} pointed out with regard to the ``Lexian ratio'' (discussed in \S\ref{sec:intro}). Consider the following extreme example. Assume that the experiment of $n=10$ events is repeated multiple times, with the result always being $k_0=3$. One simple way this could happen is if the $n=10$ events are blue (success) or red (failure) balls drawn from a bag (without replacement) that itself contains only $n=10$ balls (3 blue, 7 red). If only one such sequence exists (and absent any knowledge about the total number of balls in the bag), a Bernoulli process for the events should be given similar consideration to any other model that would roughly centralize the cumulative distribution at the observed number of successes (absent any other information about the events). The often observed deviations of real-world data from a Bernoulli process is due to similar additional constraints on the system that may not be clear from the outset. This again reinforces Wilson's expressed hesitation to assign a definite probability for the Bernoulli parameter to take on any particular range of $q$ values defined by his interval\cite{wilson_confidence_1942}.

%%%%%%%%%%%%%%%%%%%%%%%%%%%%%%%%%%%%%%%%%%%%%%%%%%%%%%%%%%%%%%
\section{Multidimensional Data}\label{sec:multidim}
%%%%%%%%%%%%%%%%%%%%%%%%%%%%%%%%%%%%%%%%%%%%%%%%%%%%%%%%%%%%%%

For multidimensional data, the situation is more complex, as there is no way to uniquely define the cumulative distribution in higher dimensional spaces (i.e. there is no natural topological ordering of points in dimensions higher than 1D). Nevertheless, in the following, I will present a method based on ``inverse Monte Carlo'' that allows extension of the fidelity to 2D (and higher) distributions for both parameter estimation and concordance determination. Significantly, the fidelity and associated concordance values obtained below are coordinate \textit{independent}. It should be noted that Ranneby and colleagues have proposed a method for extending maximum spacings to higher dimensions based on tesselations of the coordinate space representation of the data\cite{ranneby_maximum_2005}, which unfortunately makes their approach dependent on the exact choice of coordinate system (\textit{coordinate fallacy}).

The approach proposed below views inductive inference in higher dimensional spaces as a form of ``inverse Monte Carlo,'' which is already the implicit basis of the fidelity in 1D. The mapping in Fig.~\ref{fig:fidover} from the data coordinate space to the model-based cumulative interval is simply an inversion of the normal process of Monte Carlo, for which points are randomly chosen on the unit interval ($y$-axis) and then mapped onto the coordinate space ($x$-axis) via the cumulative distribution. The fidelity can be thought of as a seemingly optimal measure of how \textit{random} the inverse mapping shown in Fig.~\ref{fig:fidover} appears on the unit interval for each hypothesized model distribution, fulfilling the task that Pearson described in 1933\cite{pearson_method_1933}.

Now consider the generalization to 2D. Specifically, assume that one wishes to simulate a 2D Gaussian distribution of data points. This problem can be approached in different ways. One could discretize the space, order the bins, and then assign each event to one of the ordered bins, which are individually weighted according to the probability distribution. However, due to the discretization, this is merely an approximate simulation. A better way to simulate the data without sacrificing resolution is as follows. Events can be simulated separately in the $x$ and $y$ directions over the unit square. These events can then be stretched and rotated to conform to the desired 2D Gaussian distribution. Another method could be to generate random events in $r$ and $\theta$ over the unit sphere (circular disc). These events can then be stretched and rotated to conform to the 2D Gaussian distribution. The methods for inductive inference demonstrated below draw on these symmetries and can be considered a form of 2D ``inverse Monte Carlo'', with the task of finding the distribution that best concords with the observed data points in both dimensions.

\begin{figure}[!t]
   \begin{center}
     \vspace{-2cm}
     \includegraphics*[width=6.5in]{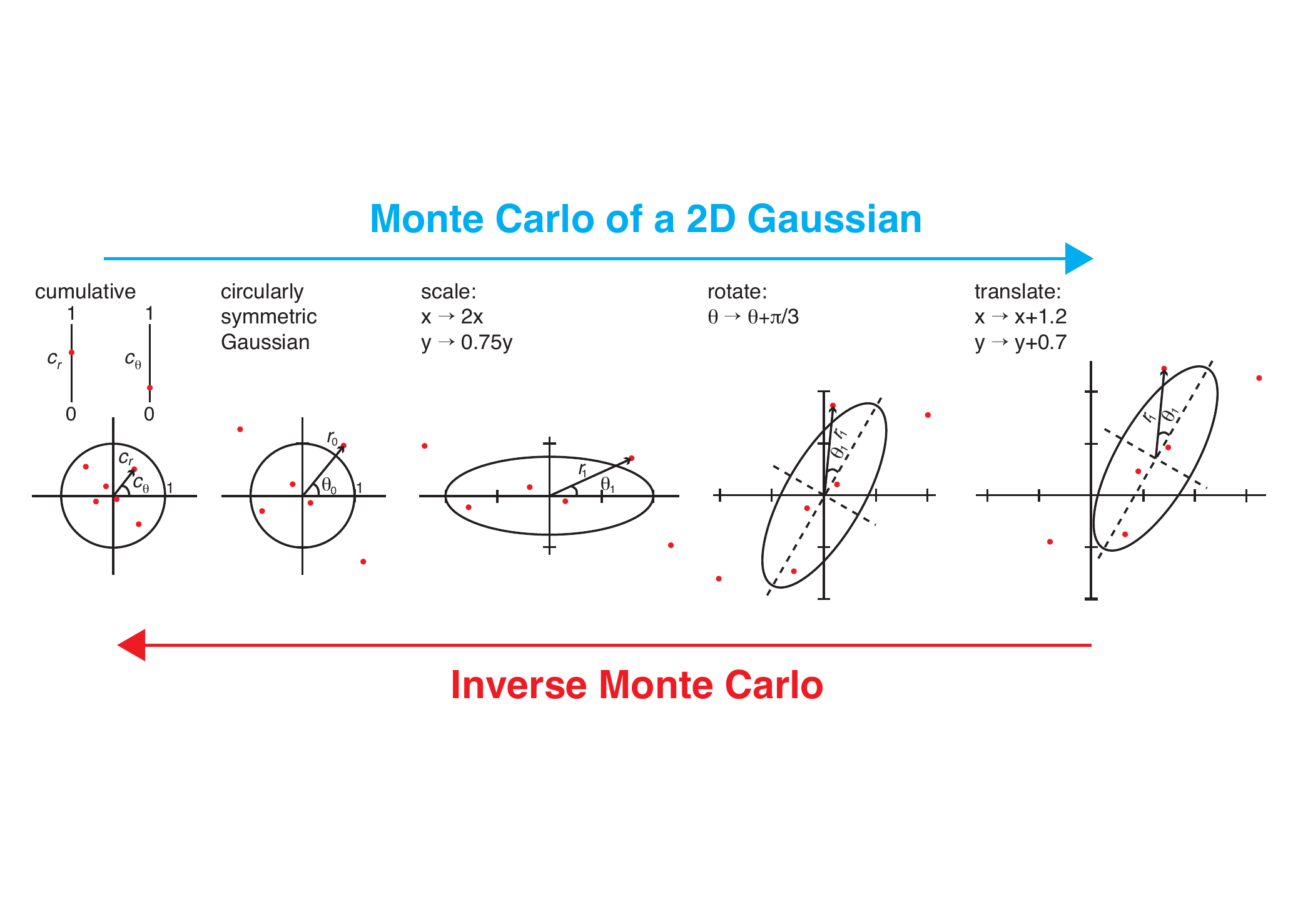}
     \vspace{-3.cm}
      \caption{Monte Carlo and inverse Monte Carlo for a 2D Gaussian. Consider first the Monte Carlo simulation of the data. First, cumulative values $c_r$ and $c_{\theta}$ are chosen from the unit interval (also shown within the unit circle, with $c_{\theta}$ ranging from 0 to 1 corresponding to a full $2\pi$ rotation). Second, these data points are converted to $r$ and $\theta$ values for a 2D Gaussian by Eq.~\ref{eq:convert}. Third, coordinates for the data points are stretched along the major/minor axes of the Gaussian, which is aligned along the $x$ and $y$ axes. Fourth, rotation by $\pi/3$ is carried out. Fifth, the points are translated to their final coordinates. Inverse Monte Carlo corresponds to mapping the observed data points to the intervals $c_r$ and $c_{\theta}$ by a reversal of the above transformations. This enables calculation of the fidelity and the corresponding $p$ value for an hypothesized 2D Gaussian model.}
      \label{fig:2D_cartoon}
   \end{center}
\end{figure}

The essential notion of 2D inverse Monte Carlo is displayed in Fig.~\ref{fig:2D_cartoon}. In this figure, the Monte Carlo-based simulation of data points drawn from a 2D Gaussian distribution is shown. Due to the symmetry of the distribution function, the following approach can be taken. First, two random numbers between 0 and 1 are drawn from a uniform distribution on the unit interval and designated $c_r$ and $c_{\theta}$. These points can be displayed on two separate unit intervals or within the unit circle (with $c_{\theta}$ going from 0 to 1 representing a $2\pi$ rotation). For each point, we can convert $c_r$ and $c_{\theta}$ to the coordinates $r$ and $\theta$ of a circularly symmetric Gaussian (second plot) according to the following equations:
\begin{eqnarray}\label{eq:convert}
r&=&\sqrt{-2\log{(1-c_r)}}\nonumber\\
\theta &=& 2\pi c_{\theta}.
\end{eqnarray}
Scaling of this distribution along the major/minor axes (aligned in the $x$ and $y$ directions) is then carried out (third plot), followed by rotation (fourth plot), and translation to reach the desired physical coordinates of the system under consideration (fifth plot). All steps in this process are of course reversible, which permits the inverse process of mapping the observed data points to the unique values $c_r$ and $c_{\theta}$ for a particular hypothesized Gaussian. With these values in hand, the fidelity and associated $p$ values can be calculated immediately using joint analysis (see \S\ref{sec:joint}) of the $n$ coordinates $c_r$ on the line and the $n$ coordinates $c_{\theta}$ on the circle. I will refer to this overall approach as the $r$-$\theta$ transform. 

\begin{figure}[!t]
   \begin{center}
     \vspace{-0.cm}
     \includegraphics*[width=3.3in]{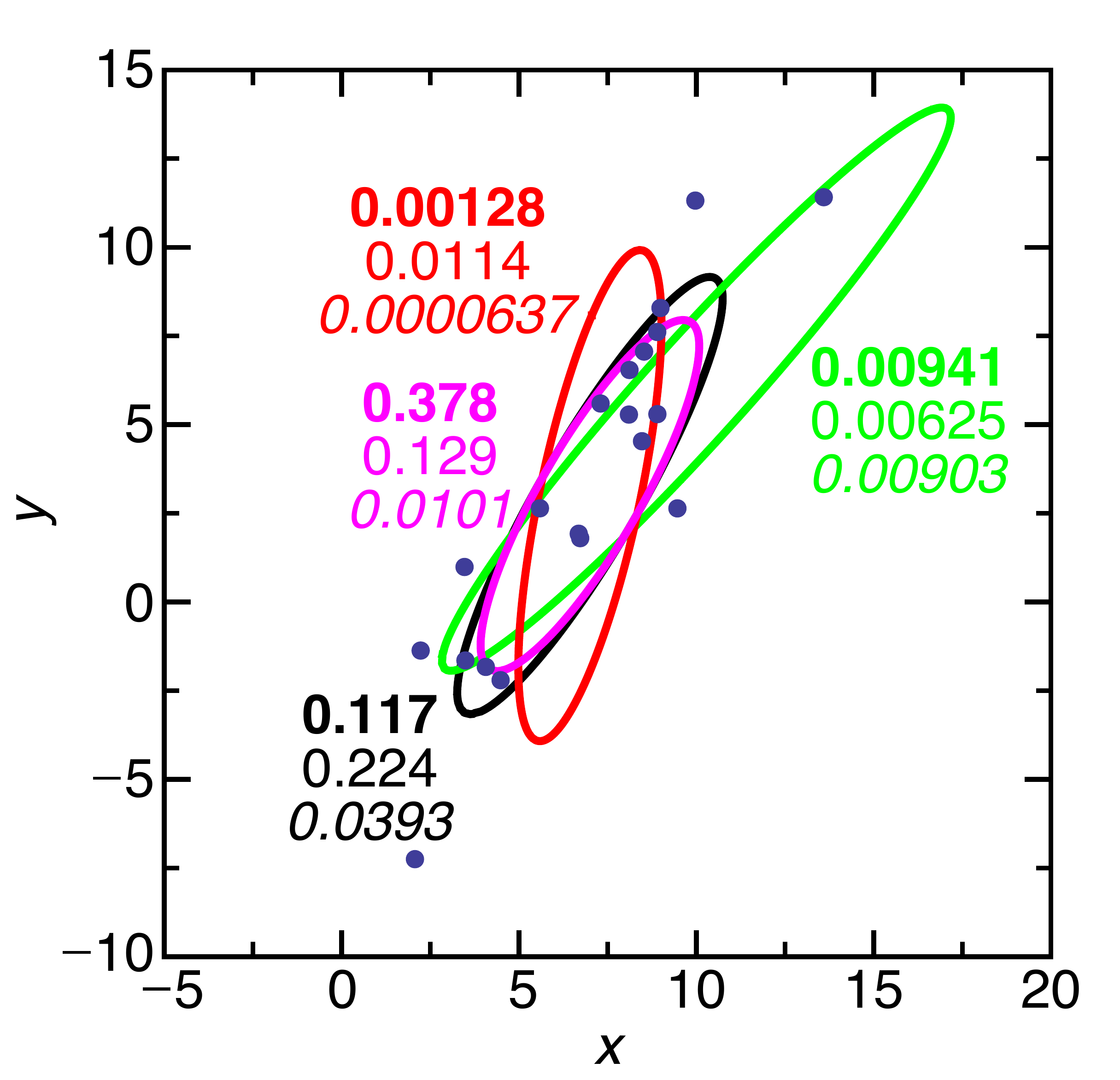}
     \vspace{-0.cm}
      \caption{2D Gaussian fitting of data drawn from a 2D Gaussian. $p$ values are shown for the $r$-$\theta$ transform (bold font), the model-based $x$-$y$ transform (regular font), and the coordinate-based $x$-$y$ transform (italic font) for each of the displayed trial models: the true model in black ($x_0=7$, $y_0=3$, $a_0=3$, $b_0=2$, $\phi_0=\pi/3$), the magenta model ($a_1=0.8 a_0$), the red model ($\phi_2=\phi_0+\pi/10$), and the green model ($x_3=x_0+3$, $y_3=y_0+3$, $a_3=1.5a_0$, $\phi_3=\phi_0-\pi/15$).}
      \label{fig:2D_Gauss}
   \end{center}
\end{figure}
\begin{figure}[]
   \begin{center}
     \vspace{0cm}
     \includegraphics*[width=3.3in]{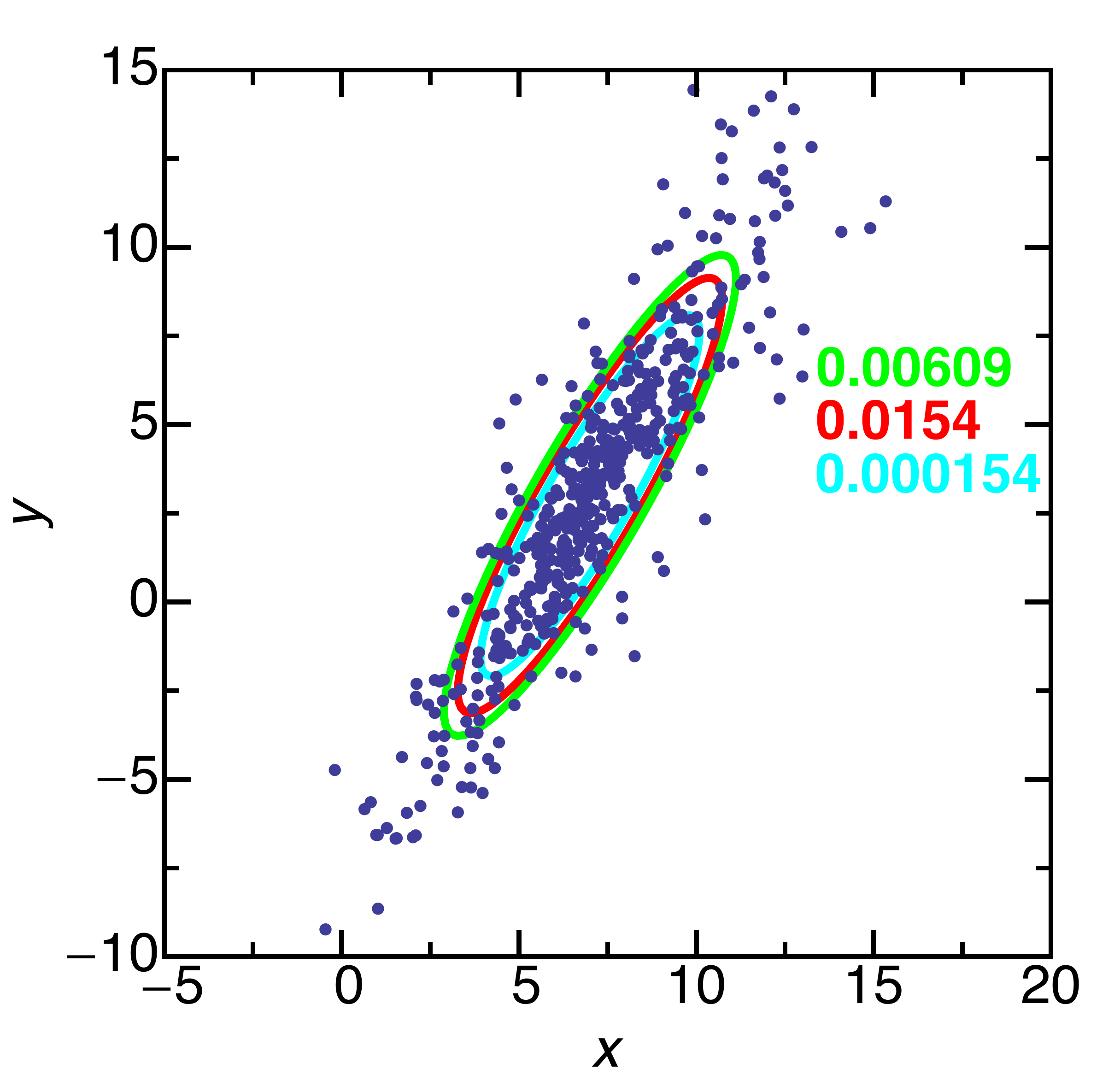}
     \vspace{-0.cm}
      \caption{2D Gaussian fitting to data drawn from a 2D Exponential ($n=500$) using the $r$-$\theta$-transform approach. The 2D Exponential model parameters are $x_0=7$, $y_0=3$, $a=5$, $b=1$, $\phi=\pi/3$ (where $a$ and $b$ give the distance to the 1/e surface of the 2D Exponential). The $p$ values for the various 2D Gaussains are listed in the plot. These three models have the same major-minor axis ratio and alignment as the underlying 2D Exponential and differ from each other only by an overall scale factor $s$ (ratio of the distance to the 1$\sigma$ surface of the 2D Gaussian to the distance to the 1/e surface of the 2D Exponential along the major axis): cyan ($s=1.5$), red ($s=1.81$), green ($s=2$). Maximization of the concordance ($p$ value) with respect to this scale factor leads to the red model with $s=1.81$ and $p=0.0154$ (the $p$ value of the correct exponential model is $p=0.328$). None of the models is therefore sufficiently concordant with the data (in the sense of $p>0.05$). }
      \label{fig:2D_exp}
   \end{center}
\end{figure}

This is not the only possible inverse Monte Carlo approach. We could have also started from the symmetric Gaussian displayed in the second panel in Fig.~\ref{fig:2D_cartoon} to separately Monte Carlo the $x$ and $y$ coordinates of the data. Due to the symmetry of the Gaussian, the $x$ and $y$ values can be independently generated, then stretched, rotated, and translated to the final physical coordinate system. This can also, of course, be reversed to yield a second inverse Monte Carlo approach based on the $c_x$ and $c_y$ values on the line. I will refer to this as the model-based $x$-$y$ transform. 

A final approach for generating cumulative values, which does not involve inverse Monte Carlo, is to simply project the observed data onto the $x$ and $y$ axes and then compute the separate cumulative projections of the hypothesized model on each axis. The observed data points can then be mapped by the projected cumulative distributions to determine their $c_x$ and $c_y$ values (which will differ from the $x$-$y$ transform). I will refer to this as the coordinate-based $x$-$y$ transform below.

In Fig.~\ref{fig:2D_Gauss}, different 2D Gaussian models are used to fit data drawn from a 2D Gaussian (black contour at $1\sigma$) with three different transform methods used for calculating the fidelity/concordance: the recommended $r$-$\theta$ transform (bold font), the model-based $x$-$y$ transform (normal font), and the coordinate-based $x$-$y$ transform (italic font). The $r$-$\theta$ transform is recommended over the other two for the following reason. One can easily consider a model that is symmetric with respect to the original physical coordinate systems (indicating no correlation), but that would project onto the physical $x$ and $y$ axes (coordinate-based $x$-$y$ transform) in a way that would mimic the projection of the Gaussian distribution shown in the last panel of Fig.~\ref{fig:2D_cartoon}. This weakness would affect the model-based $x$-$y$ transform in a similar fashion as well. Models that would project onto the $r$ direction and the $\theta$ direction in a way that mimicked a rotated Gaussian are more \textit{baroque} in their construction (e.g. a spiral pattern or a pattern with holes at various $r$ and $\theta$ that are filled in by displaced densities at different $r$ but the same $\theta$ to give a similar projection of the density across $\theta$), though these possibilities should nevertheless be kept in mind. 

Unlike the tracking of the fidelity in 1D upon acquisition of more and more data, taking more data in higher dimensions will not therefore by itself be guaranteed to conform to ``a proceeding which must in the long run approximate to the truth'' (Peirce 2.780)\cite{peirce_collected_1974}. This overall approach, however, extracts more information about the model fit than available with any other approach --- in particular, in the form of an absolute, coordinate-independent concordance measure that can be used not only to find the optimal central point of the distribution and correlation (similar to principal component analysis) but also to determine if a multidimensional Gaussian is even an appropriate model. For example, different 2D Gaussian models were used to fit data drawn from a 2D Exponential distribution in Fig.~\ref{fig:2D_exp}. 
All of these models provide inadequate fits based on the $r$-$\theta$ transform method for determining the fidelity/concordance. While all the hypothesized models used in Figs.~\ref{fig:2D_Gauss} and \ref{fig:2D_exp} were 2D Gaussians, this approach can be applied to any class of 2D models. However, symmetric models (like the 2D Gaussian or the 2D Exponential) are mathematically more tractable as the conventions for mathematical transformation of these symmetric distributions take on an obvious and (for the 2D case) unique form.

While the $r$-$\theta$ transform (based on the symmetry of the Gaussian) is uniquely defined for 2D Gaussians, for 3D Gaussians (or other similarly symmetric distributions), a coordinate convention must be chosen to transform the observed data points using a spherically symmetric coordinate system ($r$, $\theta$, $\phi$). The obvious choice is to align the longest axis of the hypothesized Gaussian along the azimuthal axis $\phi=0$. Inverse Monte Carlo would then map the observed points to cumulative values for the $r$ direction ($c_r$ values on the line), the $\theta$ direction ($c_{\theta}$ values on the circle), and the $\phi$ direction ($c_{\phi}$ values on the line). Computation of the fidelity and associate $p$ value is then straightforward. For even higher dimensional Gaussians, corresponding analogues of the Euler angles can be used with more choices of the coordinate orientation convention required. Any arbitrary $D$-dimensional distribution can be converted to the hypersphere (as in Fig.~\ref{fig:2D_cartoon}) upon a complicated enough transformation function; however, models that can easily be converted to $n$-spheres certainly simplify the process of inverse Monte Carlo and they importantly allow for universal conventions to be established (restricting the ability of the researcher to fool oneself or to fool others through statistical manipulation). As symmetric models like multidimensional Gaussians are typically invoked to explain higher dimensional data, the fidelity presents itself as a uniquely powerful tool for coordinate-independent assessment of model concordance in higher dimensional spaces.
 
%%%%%%%%%%%%%%%%%%%%%%%%%%%%%%%%%%%%%%%%%%%%%%%%%%%%%%%%%%%%%%
\section{Nonparametric Comparisons}\label{sec:nonparam}
%%%%%%%%%%%%%%%%%%%%%%%%%%%%%%%%%%%%%%%%%%%%%%%%%%%%%%%%%%%%%%

\begin{figure}[]
   \begin{center}
     \vspace{-1cm}
     \includegraphics*[width=3.5in]{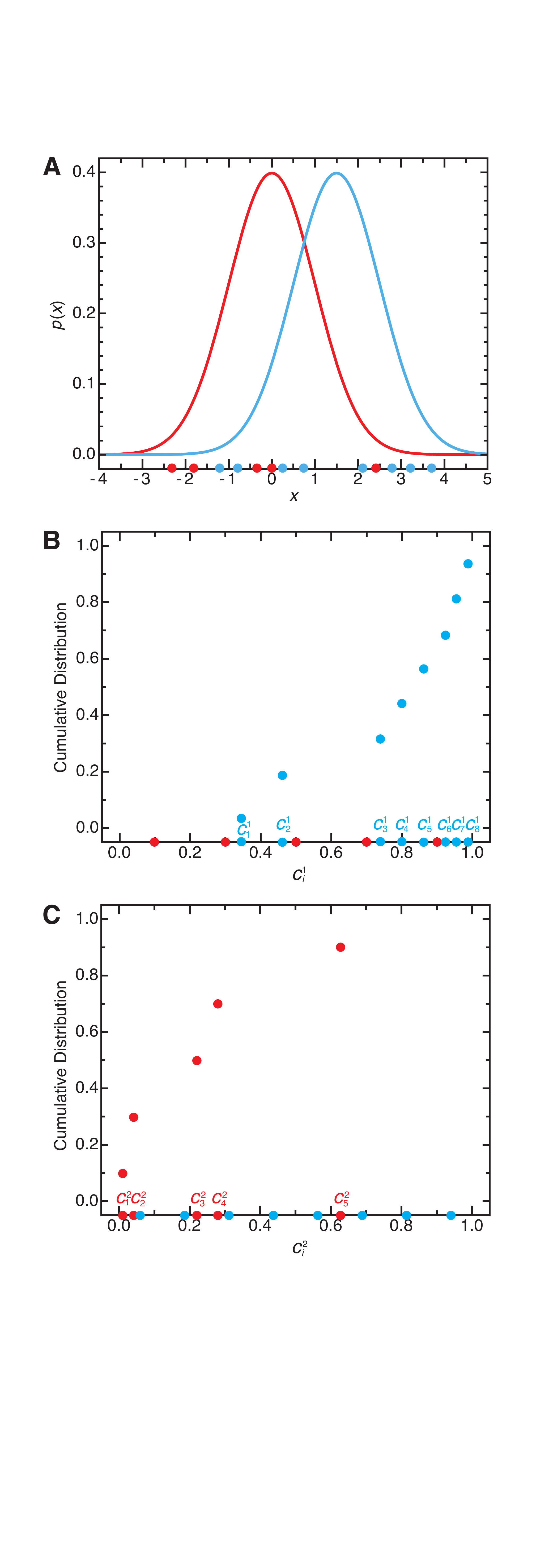}
     \vspace{-5.cm}
      \caption{Fidelity-based approach for testing whether two data sets were drawn from the same underlying model. (A) A Gaussian with $\mu_1=1.5$ and $\sigma_1=1$ (blue) and a second Gaussian with $\mu_2=0$ and $\sigma_2=1$ (red) are displayed with a representative set of data points drawn from each distribution ($n_1=8$ and $n_2=5$). To calculate a representative value for the fidelity, we consider the fidelity in both directions. (B) First, the fidelity of the blue data points are calculated with respect to the assumption of a distribution constructed from the red points. To do this, we distribute the red data points evenly across the cumulative interval ($x$-axis). We then place each blue data point at the midpoint of the particular interval in which it is located. If two or more blue points fall in the same interval, they are spread evenly over the interval in a way that maximizes their fidelity within the interval (as in Eq.~\ref{eq:fidspace}). The fidelity, $f_1$, can then be calculated according to Eq.~\ref{eq:nonparam}. (C) Here, the blue points are evenly distributed across the cumulative interval ($x$-axis) with the red points interspersed to enable calculation of $f_2$ in Eq.~\ref{eq:nonparam}. The final symmetrized fidelity is $f=(f_1+f_2)/2$ (Eq.~\ref{eq:nonparam}).} 
      \label{fig:nonparam}
   \end{center}
\end{figure}

Whether two observed data sets are drawn from the same underlying distribution is an important problem in statistics.
A fidelity-based method for comparing two empirical distributions, derived from two observed data sets with $n_1$ and $n_2$ total events, is presented in Fig.~\ref{fig:nonparam}. 
This test is compared with other non-parametric tests, specifically, Student's $t$ test, the Wilcoxon-Mann-Whitney test\cite{wilcoxon_individual_1945,mann_test_1947}, and the 2-sample Kolmogorov-Smirnov test\cite{kolmogorov_sulla_1933,smirnov_table_1948}. 

Construction of the nonparametric or model-independent form of the fidelity statistic is shown in Fig.~\ref{fig:nonparam}. In Fig.~\ref{fig:nonparam}A, eight data points are shown drawn from a Gaussian distribution with $\mu=1.5$ and $\sigma=1$ (blue) and five data points are shown drawn from a second Gaussian distribution with $\mu=0$ and $\sigma=1$ (red). To construct the relevant statistic, we assume the second distribution is known (defined by its fidelity-maximizing spacing across the cumulative interval), which then permits calculation of the fidelity of the first distribution relative to this reference (Fig.~\ref{fig:nonparam}B). Each data point from the first data set is placed in the bins defined by the second distribution in a way that maximises the fidelity within each bin (as in Eq.~\ref{eq:fidspace}). Upon switching the roles of each data set, we can determine a fidelity in the other direction as well (Fig.~\ref{fig:nonparam}C). We combine these measures into a single fidelity statistic by averaging the individual fidelities:
\begin{eqnarray}
f_1&=&\frac{1}{2n_1}\log{\left[2 n_1\left(c^1_1\right)\right]}+\frac{1}{2n_1}\log{\left[2 n_1\left(1-c^1_{n_1}\right)\right]}+\frac{1}{n_1}\sum_{i=1}^{n_1-1}\log{\left[n_1\left(c^1_{i+1}-c^1_i\right)\right]}\nonumber\\
f_2&=&\frac{1}{2n_2}\log{\left[2 n_2\left(c^2_1\right)\right]}+\frac{1}{2n_2}\log{\left[2 n_2\left(1-c^2_{n_2}\right)\right]}+\frac{1}{n_2}\sum_{i=1}^{n_2-1}\log{\left[n_2\left(c^2_{i+1}-c^2_i\right)\right]}\nonumber\\
f&=&(f_1+f_2)/2.
\label{eq:nonparam}
\end{eqnarray}
In the last line, the final combined fidelity is simply the straight average of $f_1$ and $f_2$. Other conventions for combining $f_1$ and $f_2$ could be made, such as $f'=n_1f_1+n_2f_2$, but the average of $f_1$ and $f_2$ appears to be the most reasonable, as it balances the contributions from each direction. This can be seen by examining $f'=n_1f_1+n_2f_2$. In $f'$, each data point contributes a comparable amount to the overall fidelity, as $f_1$ and $f_2$ are averages that asyptotically approach Euler's constant $\gamma$ (under the null hypothesis). However, $f_1$ and $f_2$ were obtained by fixing the reference distribution, constructed from the opposite data set, to an idealized approximation (equally spaced data points on the cumulative interval in Figs.~\ref{fig:nonparam}B and C), the reliability of which will depend on the number of data points in the opposite data set. The simplest and likely most fundamental way to incorporate this information-based \textit{reliability index} is by multiplying each term by the number of data points from the opposite data set to obtain $f''=n_2n_1f_1+n_1n_2f_2$, but this is just $2n_1n_2$ times $f$ defined above in Eq.~\ref{eq:nonparam}.

\begin{figure}[]
   \begin{center}
     \vspace{0cm}
     \includegraphics*[width=5in]{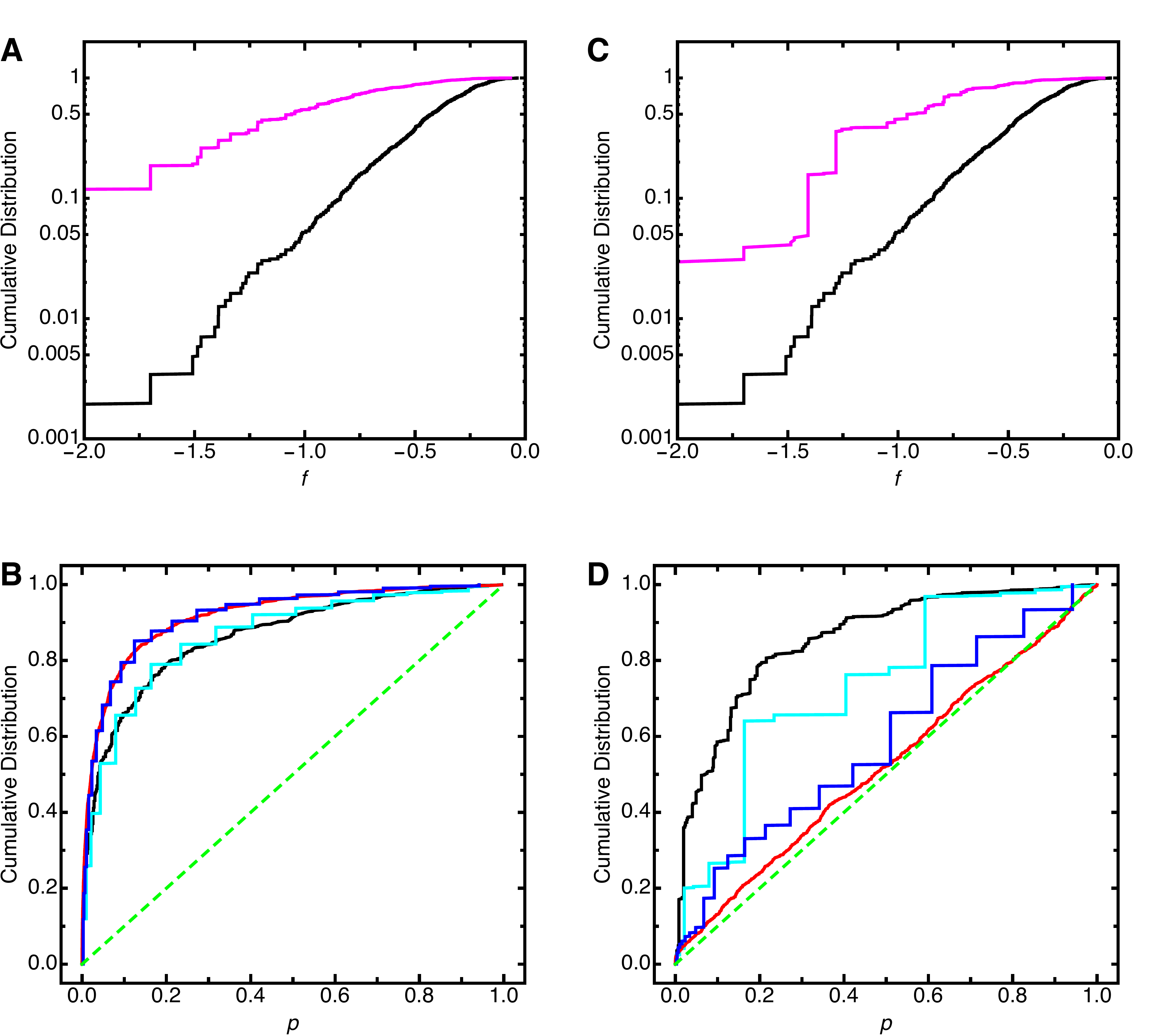}
     \vspace{0cm}
      \caption{Comparison of nonparametric tests for discriminating Gaussians. (A) The null hypothesis: $n_1=8$ data points and $n_2=5$ data points were repeatedly drawn (20000 realizations) from the uniform distribution on the unit interval to compute the cumulative fidelity distribution (black curve). The Gaussian ``location'' test: $n_1=8$ data points from a Gaussian with $\mu_1=0$ and $\sigma_1=1$ and $n_2=5$ data points from a Gaussian with $\mu_2=1.5$ and $\sigma_2=1$ were repeatably drawn (1000 realizations) to generate the cumulative distribution of the fidelity (magenta curve). (B) Through use of the null distribution (black curve in A), the fidelity distribution (magenta curve in A) was converted to the displayed cumulative distribution of $p$ values (black curve). Similar conversions were used to compute the two-sample Kolmogorov-Smirnov test (cyan) on the same data. Routines in Mathematica\textsuperscript{\textregistered} allowed immediate computation of the $p$ values for Student's $t$ test (red) and the Wilcoxon-Mann-Whitney test (blue) on the same data. The dashed green line indicates the null distribution. (C) The null hypothesis: same as in panel A for $n_1=8$ and $n_2=5$ data points (black curve). The Gaussian width test: $n_1=8$ data points were repeatedly drawn (1000 realizations) from a Gaussian with $\mu_1=0$ and $\sigma_1=1$, and $n_2=5$ data points were repeatedly drawn from a Gaussian with $\mu_2=0$ and $\sigma_2=5$ (magenta curve). (D) Through use of the null distribution (black curve in C), the fidelity distribution (magenta curve in C) was converted to the displayed cumulative distribution of $p$ values (black curve). Similar conversions were used to compute the two-sample Kolmogorov-Smirnov test (cyan) on the same data. Routines in Mathematica\textsuperscript{\textregistered} allowed immediate computation of the $p$ values for Student's $t$ test (red) and the Wilcoxon-Mann-Whitney test (blue) on the same data.  The dashed green line indicates the null distribution. The raggedness of the displayed distributions is not due to a calculational approximation (calculational uncertainty is \textit{typically} much smaller than these features), but is a fundamental aspect of the tests (i.e. the discrete nature of the computed statistics).}
      \label{fig:nonparam_gauss_gauss}
   \end{center}
\end{figure}

\begin{figure}[]
   \begin{center}
     \vspace{0.5cm}
     \includegraphics*[width=5in]{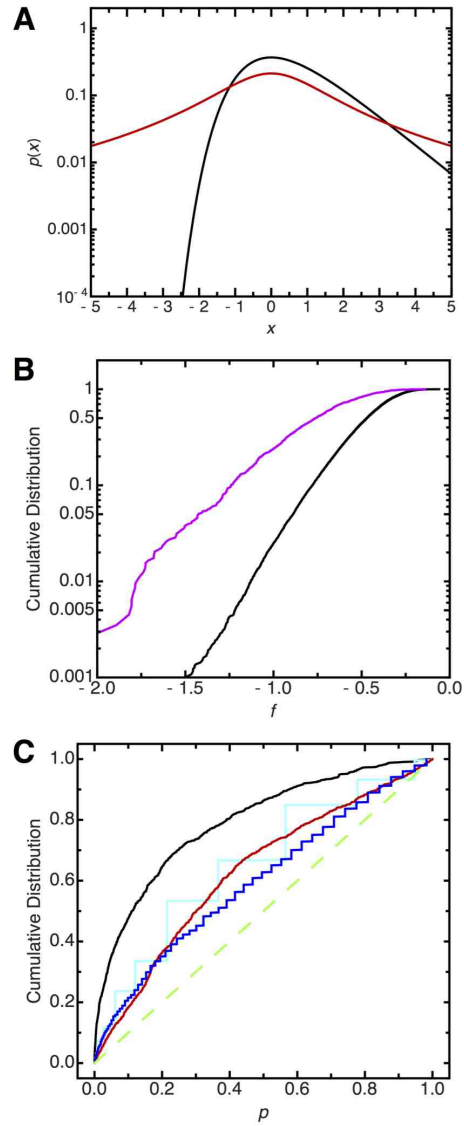}
     \vspace{-2.7cm}
      \caption{Comparison of non-parametric tests for discriminating more general distributions. (A) $n_1=20$ data points were drawn from an Extreme Value distribution with $\beta_1=0$ and $\alpha_1=1$ (black), and $n_2=10$ data points were drawn from a Cauchy distribution with $\beta_2=0$ and $\alpha_2=1.5$ (red) (see Table~\ref{tab:line_dist}). (B) The null distribution of the fidelity for $n_1=20$ and $n_2=10$ (black, constructed from 20000 Monte Carlo realizations) is shown along with the cumulative distribution of the fidelity obtained for the particular pair of distributions (magenta, 1000 realizations). (C)  Through use of the null distribution (black curve in B), the fidelity distribution (magenta curve in B) was converted to the displayed cumulative distribution of $p$ values (black curve) Also shown are the cumulative $p$ value distributions obtained from Student's $t$ test (red), the Wilcoxon-Mann-Whitney test (blue), and the two-sample Kolmogorov-Smirnov test (cyan) on the same data (see Fig.~\ref{fig:nonparam_gauss_gauss} for more details). The dashed green line indicates the null distribution.}
      \label{fig:nonparam_extreme_cauchy}
   \end{center}
\end{figure}

In Fig.~\ref{fig:nonparam_gauss_gauss}, the results of applying this empirical test based on the fidelity to Gaussians that differ only in location (left panels) and Gaussians that differ only in width (right panels) is shown in comparison with other tests: Student's $t$ test (red curve), the Wilcoxon-Mann-Whitney test (blue curve), and the two-sample Kolmogorov-Smirnov test (cyan curve). For ``location'' testing, both the fidelity and the two-sample Kolmogorov-Smirnov test have slightly less power than the other two tests. However, the fidelity is clearly superior in discriminating differences in width compared to these other tests (in fact, Student's $t$ test and the Wilcoxon-Mann-Whitney test have no power to discriminate width differences due to their mathematical definitions). That the fidelity represents a general test with discriminatory power for arbitrary distributions is shown in Fig.~\ref{fig:nonparam_extreme_cauchy} for discrimination of an Extreme Value distribution from a Cauchy distribution. The fidelity is clearly superior to the other three tests.

While showing a slight decrease in sensitivity for ``location'' testing (as compared to the exclusive ``location'' tests corresponding to Student's $t$ test and the Wilcoxon-Mann-Whitney test), the fidelity more than makes up for this by its ability to test more general differences between distributions (different widths or shapes). Perhaps of most significance, it equals or surpasses the two-sample Kolmogorov-Smirnov test, which represents the ``gold standard'' for general nonparametric testing. The fidelity also does not suffer from the asymmetric discrimination problems of the Kolmogorov-Smirnov test (better discrimination for differences near the median than at the boundaries), as it treats all discrepancies between the distributions equally, no matter their location on the cumulative intervals (Figs.~\ref{fig:nonparam}B and C). Similar results and behavior are expected for nonparametric testing on the circle as well (upon the obvious modification of Fig.~\ref{fig:nonparam} at the boundaries).

\newpage

%%%%%%%%%%%%%%%%%%%%%%%%%%%%%%%%%%%%%%%%%%%%%%%%%%%%%%%%%%%%%%
\section{Discussion}\label{sec:discussion}
%%%%%%%%%%%%%%%%%%%%%%%%%%%%%%%%%%%%%%%%%%%%%%%%%%%%%%%%%%%%%%

\begin{quote}
``Thus it is that inquiry of every type, fully carried out, has the vital power of self-correction and of growth. This is a property so deeply saturating its inmost nature that it may truly be said that there is but one thing needful for learning the truth, and that is a hearty and active desire to learn what is true. If you really want to learn the truth, you will, by however devious a path, be surely led into the way of truth, at last. No matter how erroneous your ideas of the method may be at first, you will be forced at length to correct them so long as your activity is moved by that sincere desire. Nay, no matter if you only half desire it, at first, that desire would at length conquer all others, could experience continue long enough. But the more veraciously truth is described at the outset, the shorter by centuries will the road to it be.'' (Peirce 5.582)\cite{peirce_collected_1974}
\end{quote}
Conventional statistical approaches by ``however devious a path'' will generally allow one to arrive at the ``truth'', as long as one keeps in mind the particular limitations of the assumed approaches. However, as I hope to have demonstrated in this manuscript, maximum fidelity represents a particularly ``veracious'' method for identifying the model (or set of models) that best concord with the data, serving as a highly efficient and seemingly optimal basis for statistical inference.

I have argued that maximum fidelity is superior to all other methods, including maximum likelihood, for parameter estimation (see \S\ref{sec:estimation}). The likelihood was previously considered as \textit{fundamental} by Fisher and by many others throughout the twentieth century, however its spectacular failure for parameter estimation within certain distribution families tarnished its status\cite{stigler_epic_2007} (see also the references in \S\ref{sec:intro}). The superiority of parameter estimation by maximization of the fidelity across all of the diverse cases tested in \S\ref{sec:estimation}, as well as the fact that it never fails no matter the distribution family (due to its boundedness from above\cite{ranneby_maximum_1984}), together demonstrate that the fidelity is more \textit{fundamental} than the likelihood. In \S\ref{sec:fidlike}, I have shown that the likelihood can be derived as an asymptotic approximation to the fidelity on the circle and on the line. Whatever good qualities the likelihood possesses for the analysis of univariate data can therefore be attributed to its asymptotic approximation of the fidelity.

The fidelity is also superior to the spacings statistic for parameter estimation. Ranneby claimed that the spacings statistic represented a better approximation than the likelihood for the Kullback-Leibler divergence\cite{ranneby_maximum_1984}; however, he wisely equivocated on the use of the maximum spacings estimage (MSP) versus the maximum likelihood estimate (MLE) for small sample parameter estimation:
\begin{quote}
``To give rules for choosing between the MSP estimate and the MLE when they are asymptotically equivalent we have to know more, especially for small samples.''\cite{ranneby_maximum_1984}
\end{quote}
Small-sample comparisons of the estimates obtained from maximum spacings and from maximum likelihood have not successfully demonstrated the superiority of either statistic\cite{ekstrom_alternatives_2008}. In light of this and, more significantly, of the superiority of the fidelity for general parameter estimation (\S\ref{sec:estimation}), it is clear that the spacings statistic should also be regarded as an approximation of the more fundamental fidelity statistic, which provides the best representation of the \textit{fidelity} (Kullback-Leilber divergence) of the model with the observed data. As argued in \S\ref{sec:fidelity}, the success of the fidelity statistic over the spacings statistic can be attributed to the former's better respect for the density/symmetry of the contribution of each data point  (and not the \textit{spacings} they create) across the cumulative interval.

Complete avoidance of the \textit{probability fallacy} and the \textit{parameter fallacy} prohibits the construction of a quantitative \textit{universal} measure for assessing the ``optimality'' of a parameter estimation method that would be valid across all distribution families. What the fidelity accomplishes  is the determination of the model distribution (often from a restricted set of distributions) that best summarizes the \textit{local} information represented by the data. That the fidelity's consideration of only \textit{local} information turns out to be ``optimal'' for parameter estimation --- in the loose sense of the term as used in \S\ref{sec:estimation} to refer to low median bias and a narrow distribution about the true value --- is a highly interesting byproduct. What we are driven to conclude is that we should have always been on the search for the distribution that maximizes the fidelity. The fidelity is the quantity on which it makes the most sense to base statistical inference. 

The best justifications for the fidelity, therefore, lie in its empirically determined estimation optimality, its reliable coordinate-independent basis on the cumulative distribution (unlike maximum likelihood, which cannot be compared across different distribution families and even fails for parameter estimation within certain distribution families), and its intuitive and \textit{balanced} consideration of the local density contributed by each data point across the cumulative interval (in contrast to maximum spacings, see Fig.~\ref{fig:symmetry}).

I have also argued that the concordance value associated with the fidelity is the most \textit{generally} optimal approach for discriminating incorrect models for data on the circle and on the line (see \S\ref{sec:concordance}). Due to its basis on the cumulative distribution, both the fidelity and its associated concordance value have an absolute meaning that is invariant to the choice of coordinate system for the data (with the only restriction being to coordinate systems that share the same topological ordering of the data points, see \S\ref{sec:neyman}). As the fidelity groups together distributions based on their ``frequency of occurrence'' (according to the data ``density'' assumption, see \S\ref{sec:fidelity}), it appears to be a fundamental measure of the directional distance of an hypothesized distribution from the data-inferred estimate of the underlying distribution. It is important to note that the fidelity is only based on the logarithmic sum of the distribution of the sizes of the cumulative intervals, not on their particular order. The fidelity tests only \textit{zeroth order} information and is therefore completely insensitive to the \textit{location} of multiple shorter-than-average intervals across the cumulative interval, which might be dispersed randomly over the cumulative interval, might lie directly next to each other (as could be revealed by a more sensitive but less general ``location'' or ``clustering'' test), or might be distributed in a regular, evenly-spaced fashion, etc. This explains the slight reduction in sensitivity of the fidelity upon the (\textit{artificial}) restriction of testing to the ``location'' of a model distribution of fixed shape (see Fig.~\ref{fig:good_Gauss}B). However, the fidelity exhibits much greater power for discriminating more general differences between the distributions. The fidelity might assign equal discrepancies to a ``location'' difference as to a local concentration difference, with a clear meaning of the degree of this discrepancy in both cases from the definition of the fidelity of the distribution in terms of its ``frequency of occurrence'' (see \S\ref{sec:fidelity}). The extension of maximum fidelity to the model-independent assessment of whether two data sets were drawn from the same underlying distribution in \S\ref{sec:nonparam} demonstrates that the so-defined fidelity possesses a similar sensitivity to more general distribution discrepancies at the sacrifice of some power in the testing of distribution ``location'' (compared, in this case, to the \textit{pure} ``location'' tests represented by Student's $t$ test and the Wilcoxon-Mann-Whitney test).

That the fidelity represents the zeroth order step in an analysis certainly does not preclude examination of the data at a higher order. It is nevertheless remarkable how much can already be obtained upon this zeroth order consideration (optimal parameter estimation and \textit{generally} optimal concordance assessment, both in a coordinate-independent manner). The different orders of data analysis are represented well in Peirce's examination of the randomness of the digits of $\pi$ (see also \S\ref{sec:binary}):
\begin{quote}
``In order to illustrate this mode of induction, I have made a few observations on the calculated number. There ought to be, in 350 successive figures, about 35 fives. The odds are about 2 to 1 that there will be 30-39 [and] 3 to 1 that there will be 29-41. Now I find in the first 350 figures 33 fives, and in the second 350, 28 fives, which is not particularly unlikely under the supposition of a chance distribution. During the process of counting these 5's, it occurred to me that as the expression of a rational fraction in decimals takes the form of a circulating decimal in which the figures recur with perfect regularity, so in the expression of a quantity like $\pi$, it was naturally to be expected that the 5's, or any other figure, should recur with some approach to regularity. In order to find out whether anything of this kind was discernible I counted the fives in 70 successive sets of 10 successive figures each. Now were there no regularity at all in the recurrence of the 5's, there ought among these 70 sets of ten numbers each to be 27 that contained just one five each; and the odds against there being more than 32 of the seventy sets that contain just one five each is about 5 to 1. Now it turns out upon examination that there are 33 of the sets of ten figures which contain just one 5. It thus seems as if my surmise were right that the figures will be a little more regularly distributed than they would be if they were entirely independent of one another. But there is not much certainty about it. This will serve to illustrate what this kind of induction is like, in which the question to be decided is how far a given succession of occurrences are independent of one another and if they are not independent what the nature of the law of their succession is.'' (Peirce 7.121)\cite{peirce_collected_1974}
\end{quote}
When presented with the digits of $\pi$, the first natural question one should ask is ``Does each digit appear roughly the same amount of times as every other digit?" This zeroth order question is a question concerning \textit{local} information, in the same sense as the fidelity. No correlations with other digits are considered. One simply counts the occurrence of each digit in the first $N$ digits and compares with the statistical expectation (for the binary expression of $\pi$, this would involve a comparison with a Bernoulli process having $q=0.5$, as discussed in \S\ref{sec:binary}). In Peirce's example, he finds the first 350 digits of $\pi$ contain 33 fives, in accord with statistical expectations for a Bernoulli process with $q=0.1$. However, he then searches for higher order discrepancies, which may still be present. He considers the mean number of 5's in every 10 digits of $\pi$ in order to see if there is any regularity on this 10 digit ``length scale''. He finds marginal evidence for a discrepancy from a pure Bernoulli process, but the evidence is not statistically significant. Of course, the choice of every 10 digits was arbitrary. One might also look for correlations of digits that are displaced from one another (e.g. examine the distribution of 5's in every other digit or by skipping every 2 digits, etc.). Due to the many possible correlations that can be defined, and might actually be present in any data set, these higher order investigations require a degree of caution that is unnecessary at the zeroth order (at least for the analysis of the digits of $\pi$ or, as throughout this manuscript, for the testing of model concordance with univariate data through use of the fidelity). Note that ``location'' or ``clustering'' tests essentially skip this zeroth order step and proceed directly to a higher order of analysis, with the higher order test often based on assumptions that may already be clearly violated at the zeroth order. For example, one might incorrectly apply Student's $t$ test to two data sets that are not drawn from Gaussians; by contrast, the fidelity-based generalization of Student's $t$ test shown in Figs.~\ref{fig:gaussian_t_test} and \ref{fig:extreme_value_t_test} allows examination of the concordance of each model distribution with the data set it is supposed to represent ($p_1$ and $p_2$) as well as the absolute concordance of the joint model fit ($p$). 

The extension of maximum fidelity to higher dimensional data is possible in a coordinate-independent fashion but such an extension is not unique due to the lack of a unique way to define the cumulative distribution in higher dimensions (see \S\ref{sec:multidim}). Whether the extension to higher dimensional data is straightforward or not depends on the symmetry properties of the coordinate system along with the class of models under consideration. This is not a \textit{specific} weakness of maximum fidelity as compared to other approaches (as few techniques in general exist for higher dimensional data), but rather a general recognition of the significantly more challenging problems presented by higher dimensional spaces. 

\begin{figure}[]
   \begin{center}
     \vspace{-4cm}
     \includegraphics*[width=5in]{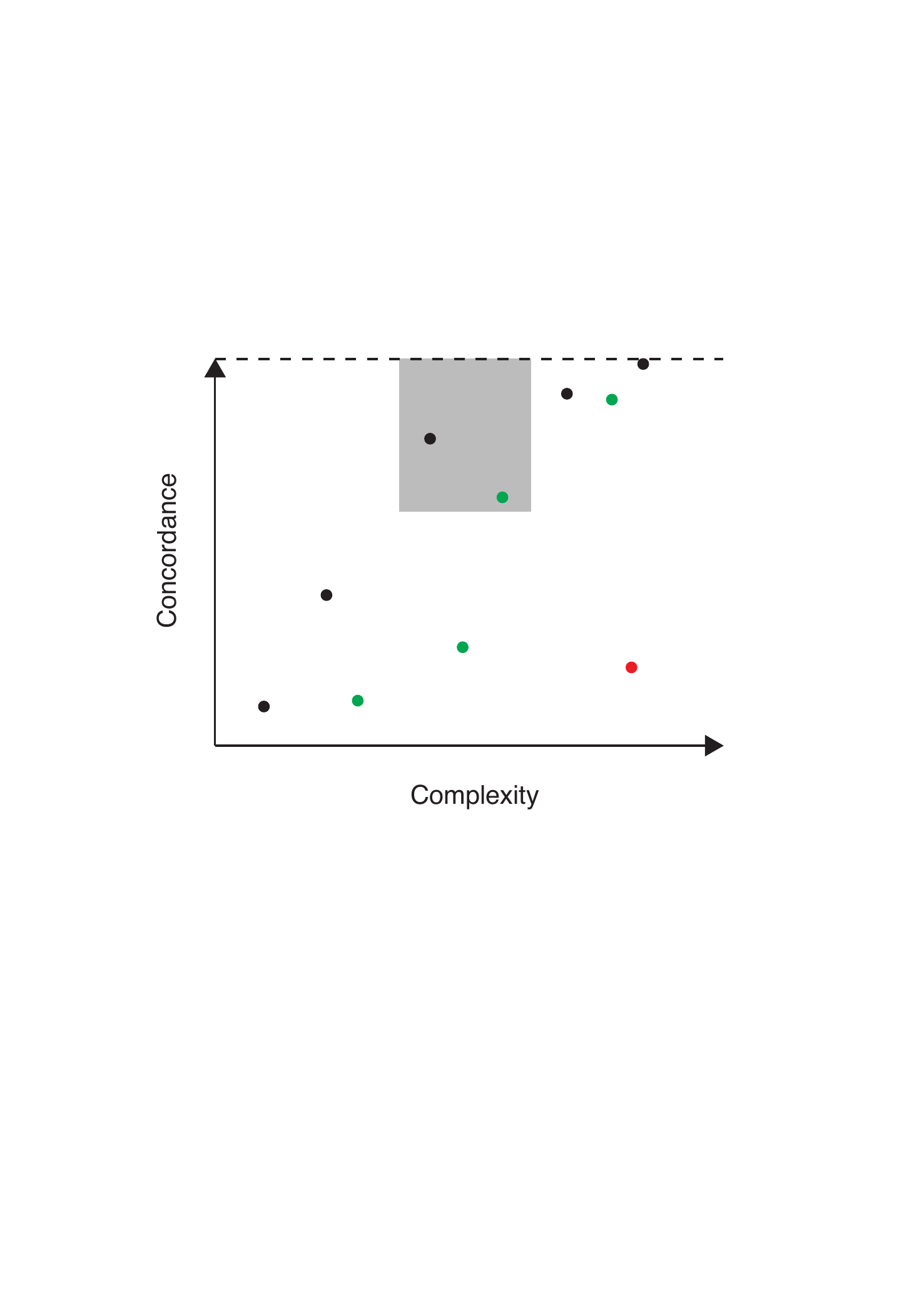}
     \vspace{-6.5cm}
      \caption{Efficient scientific inference requires a balance of model concordance with model simplicity. Each point represents a particular theory's concordance with the data versus its complexity. The dashed line indicates perfect concordance of the model with the data (e.g. perfect fidelity with $p=1$). The different colors represent different classes of theories, with the progression from left to right indicating the inclusion of additional parameters to bring the model into concordance with the data. Of course, complex models can also be highly discrepant with the data (red point). The grey zone indicates the region of optimal scientific inference, corresponding to models with sufficient concordance that are not too complex (if the model is too complex, it would provide little insight). Models that are too simple are not worth considering as they will generally be highly discrepant with the data (note the absence of points in the upper left corner). A ``good model'' (like the black or green dots in the grey zone) may not be perfectly concordant with the data due to the remaining statistical uncertainty.}
      \label{fig:ockham}
   \end{center}
\end{figure}

If scientific inference were based only on finding the most concordant model to the data, then we should simply seek out models with as many parameters necessary to fit each data set, which would amount to an intellectually fallow \textit{descriptive empiricism}. But science has successfully been shown (often enough anyway!) to be based on simple theories that can summarize large amounts of data. To be efficient, scientific inference must be based on Ockham's razor, as Peirce also recognized, but in a \textit{confused} manner that I will address below:
\begin{quote}
``Parsimony (law of): Ockham's razor, i.e. the maxim `\textit{Entia non sunt multiplicanda praeter necessitatem}.' The meaning is, that it is bad scientific method to introduce, at once, independent hypotheses to explain the same facts of observation.

Though the maxim was first put forward by nominalists, its validity must be admitted on all hands, with one limitation; namely, it may happen that there are two theories which, so far as can be seen, without further investigation, seem to account for a certain order of facts. One of these theories has the merit of superior simplicity. The other, though less simple, is on the whole more likely. But this second one cannot be thoroughly tested by a deeper penetration into the facts without doing almost all the work that would be required to test the former. In that case, although it is good scientific method to adopt the simpler hypothesis to guide systematic observations, yet it may be better judgment, in advance of more thorough knowledge, to suppose the more complex hypothesis to be true. For example, I know that men's motives are generally mixed. If, then, I see a man pursuing a line of conduct which apparently might be explained as thoroughly selfish, and yet might be explained as partly selfish and partly benevolent, then, since absolutely selfish characters are somewhat rare, it will be safer for me in my dealings with the man to assume the more complex hypothesis to be true; although were I to undertake an elaborate examination of the question, I ought to begin by ascertaining whether the hypothesis of pure selfishness would quite account for all he does.'' (Peirce 7.92-93)\cite{peirce_collected_1974}
\end{quote}
There is actually no evidence that Ockham ever wrote ``\textit{Entia non sunt multiplicanda praeter necessitatem}'' (``Entities should not be multiplied without necessity''), but what he did write in a similar vein was ``\textit{Numquam ponenda est pluralitas sine necessitate}'' (``Plurality must never be posited without necessity'') and ``\textit{Frustra fit per plura quod potest fieri per pauciora}'' (``It is futile to do with more things that which can be done with fewer'')\cite{thorburn_myth_1918}. Ockham, along with many other medieval scholars, was interested in categorizing \textit{all} of existence into the simplest possible taxonomic tree containing the fewest branches. Insertion of additional categories, when fewer would suffice, was considered a roadblock for understanding the fundamental organization of the entire universe. Elsewhere, Peirce arrives closer to the true spirit of Ockham's razor:
\begin{quote}
``Science ought to try the simplest hypothesis first, with little regard to its probability or improbability, although regard ought to be paid to its consonance with other hypotheses, already accepted.'' (Peirce 4.1)\cite{peirce_collected_1974}
\end{quote}
``Consonance with other hypotheses, already accepted'' is absolutely key to the very definition of simplicity, a point which Peirce failed to fully appreciate in his own example above. Simplicity should never be examined in a vacuum, but always in the context of one's theory of the entire universe, in line with Peirce's central maxim of his philosophy of \textit{Pragmatism}:
\begin{quote}
``Consider what effects, that might conceivably have practical bearings, we conceive the object of our conception to have. Then, our conception of these effects is the whole of our conception of the object.'' (Peirce 5.2)\cite{peirce_collected_1974}
\end{quote}
Any theory we posit for a given data set must not unnecessarily overcomplicate our model for the \textit{entire universe}.  For Peirce's example above, assume for simplicity that ``purely selfish'' acts are not only rare but not known to exist, and also assume that we have no other knowledge about the act under consideration. Therefore, to assume the man is acting with ``purely selfish'' intentions would not be the simplest assumption, as it would create a new category not present in the known universe. In this case, Peirce errs in considering the ``partly selfish and partly benevolent'' hypothesis as the more complex one. In fact, it is the simpler hypothesis (it leads to the simplest extension of our universe model), perfectly in line with Ockham's razor. The scientific method is based on positing the simplest theories that are \textit{sufficiently} concordant with previous knowledge for further empirical testing. Maximum fidelity allows a universal notion of concordance that should help assist in testing a range of models with different degrees of complexity for further testing (Fig.~\ref{fig:ockham}). The simplest models pose questions that allow for the fastest gain in knowledge, in accord with Peirce's revealing analogy of scientific inference to the game of ``twenty questions'':
\begin{quote}
``The qualities which these considerations induce us to value in a hypothesis are three, which I may entitle Caution, Breadth, and Incomplexity. In respect to caution, the game of twenty questions is instructive. In this game, one party thinks of some individual object, real or fictitious, which is well-known to all educated people. The other party is entitled to answers to any twenty interrogatories they propound which can be answered by Yes or No, and are then to guess what was thought of, if they can. If the questioning is skillful, the object will invariably be guessed; but if the questioners allow themselves to be led astray by the will-o-the-wisp of any prepossession, they will almost as infallibly come to grief. The uniform success of good questioners is based upon the circumstance that the entire collection of individual objects well-known to all the world does not amount to a million. If, therefore, each question could exactly bisect the possibilities, so that yes and no were equally probable, the right object would be identified among a collection numbering $2^{20}$ \ldots or over one million and forty-seven thousand, or more than the entire number of objects from which the selection has been made. Thus, twenty skillful hypotheses will ascertain what two hundred thousand stupid ones might fail to do. The secret of the business lies in the caution which breaks a hypothesis up into its smallest logical components, and only risks one of them at a time. What a world of futile controversy and of confused experimentation might have been saved if this principle had guided investigations into the theory of light! The ancient and medieval notion was that sight starts from the eye, is shot to the object from which it is reflected, and returned to the eye. This idea had, no doubt, been entirely given up before R\"omer showed that it took light a quarter of an hour to traverse the earth's orbit, a discovery which would have refuted it by the experiment of opening the closed eyes and looking at the stars. The next point in order was to ascertain of what the ray of light consisted. But this not being answerable by yes or no, the first question should have been `Is the ray homogeneous along its length?' Diffraction showed that it was not so. That being established, the next question should have been `Is the ray homogeneous on all sides?' Had that question been put to experiment, polarization must have been speedily discovered; and the same sort of procedure would have developed the whole theory with a gain of half a century.'' (Peirce 7.220)\cite{peirce_collected_1974}
\end{quote}
The absolute \textit{simplicity} or generality of a given scientific hypothesis or question is of course impossible to quantify. Even if it could be quantified, many hypotheses might still have a similar simplicity (even in an exact mathematical sense). One should, however, not make the mistake of relying on mathematical and/or statistical reasoning alone to define what is \textit{simple} (e.g. by simply counting model parameters). The notion of simplicity is often highly controversial (indeed, many scientific debates are ultimately based on the question of which hypothesis or explanation is the simplest or the most worthwhile for further testing), but the process elegantly described above by Peirce nevertheless plays an essential guiding principle for truly efficient scientific inference.

There are a few obvious open problems immediately suggested from the results of this manuscript. A mathematical proof that can account for the exact overlap of the parameter estimate distribution derived from maximum fidelity and that derived from the standard deviation for the Gaussian width $\sigma$ (Fig.~\ref{fig:gauss_mean_sigma}) would be of great interest. Also of interest would be a \textit{rigorous} demonstration of the asymptotic equivalence of maximum fidelity and minimum $\chi^2$ for binned data (see \S\ref{sec:binned}). Universal conventions would be worthwhile to define for the analysis of 2D distributions on the surface of the sphere or on the torus, for higher dimensional Gaussians (ordering of the Euler angle axes), and for other general higher dimensional distributions. A rather deeper line of enquiry would be whether there exists a \textit{generic} foundation for taking into account the ``location'' or ``clustering'' of shorter-than-average vs. longer-than-average cumulative intervals along the entire cumulative distribution. Discovery of a general approach (or approaches) for analyzing successively higher orders of information would of course be worthwhile and would help to ground the current \textit{ad hoc} nature of popular ``location'' or ``clustering'' tests. The best way to mathematically integrate the asymptotic notion of ``degrees of freedom'' to ``correct'' the absolute $p$ value (obtained from maximization of the fidelity over some or all of the model parameters) based on the work of Cheng \& Stephens for goodness-of-fit using the spacings statistic might also be of interest to many researchers\cite{cheng_goodness--fit_1989}; such ``corrections'' should nevertheless be interpreted very carefully (as should the absolute $p$ value for that matter) for the reasons given in \S\ref{sec:concordance}.

Maximum fidelity can be applied in a unique and universal fashion for the testing of arbitrary models against the complete \textit{local} information present in a univariate data set (or across multiple univariate data sets). Upon a generic choice of convention, it can also be used to efficiently extract information from multivariate data, including a coordinate-independent assessment of model concordance not possible with any other technique. The results in this manuscript prove that statistical inference can \textit{indeed} be based completely on the fundamental notion of model concordance, with maximum fidelity representing a ``proper form of probable inference''\cite{wilson_probable_1927} through its complete avoidance of the \textit{probability} and \textit{parameter fallacies} required for Bayesian and most frequentist approaches to statistics. Most intriguingly, the results obtained herein suggest that the most optimal and efficient basis for statistical inference is through the assessment of the information contained in the \textit{cumulative distribution}, establishing a specific fundamental connection between optimal statistical inference and information theory that had either not been so envisioned (e.g. the work of Jaynes\cite{jaynes_information_1957}),  or even if envisioned, as for the spacings statistic\cite{kale_test_1967,gebert_goodness_1969,kale_unified_1969,cheng_estimating_1983,ranneby_maximum_1984,ekstrom_alternatives_2008}, never convincingly established.

\section*{Acknowledgements}

This manuscript was written using \LaTeX{} and \BibTeX\  within the TeXShop environment.  All analysis was performed in Mathematica\textsuperscript{\textregistered} with additional use of routines for Stineman interpolation (S. Wagon, http://library.wolfram.com/infocenter/Articles/2176/) and plot tick formatting (M. Caprio, 
http://library.wolfram.com/infocenter/MathSource/5727/).  Plots were created in Mathematica\textsuperscript{\textregistered} with additional formatting and annotation created within Adobe\textsuperscript{\textregistered} Illustrator\textsuperscript{\textregistered}. Referencing of primary sources would not have been possible without the extensive digital archives maintained by the cited journals, JSTOR\textsuperscript{\textregistered}, and Project Euclid. Many helpful Wikipedia\textsuperscript{\textregistered} webpages on statistical approaches and their history were consulted in the preparation of this manuscript.

\bibliography{fidelity}{}
\bibliographystyle{ieeetr}

\end{document}